\documentclass[psamsfonts]{amsart}
\usepackage[T1]{fontenc}
\usepackage{geometry}
\usepackage{graphicx}
\usepackage{amssymb}

\makeatletter
 \theoremstyle{plain}
\newtheorem{thm}{Theorem}[section]
  \theoremstyle{definition}
  \newtheorem{defn}[thm]{Definition}
  \theoremstyle{plain}
  \newtheorem{prop}[thm]{Proposition}
  \theoremstyle{remark}
  \newtheorem{rem}[thm]{Remark}
\newenvironment{lyxlist}[1]
{\begin{list}{}
{\settowidth{\labelwidth}{#1}
 \setlength{\leftmargin}{\labelwidth}
 \addtolength{\leftmargin}{\labelsep}
 }}
{\end{list}}
 \theoremstyle{definition}
  \newtheorem{example}[thm]{Example}
  \theoremstyle{remark}
  \newtheorem*{acknowledgement*}{Acknowledgement}

\usepackage{graphicx}
\usepackage{color}
\definecolor{rltblue}{rgb}{0,0,0.75}
\usepackage{caption}
\usepackage{euscript}
\usepackage{eufrak}
\usepackage{threeparttable}
\usepackage{dcolumn}

\theoremstyle{definition}

\theoremstyle{remark}

\numberwithin{equation}{section}

\renewcommand{\subsection}{\@startsection%
 {subsection}%
 {2}%
 {0mm}%
 {0mm}%
 {0.1\baselineskip}%
 {\normalfont\bfseries}}
\newcolumntype{d}[1]{D{.}{.}{#1}}
\newcolumntype{f}[1]{D{+}{+}{#1}}

\newcommand{\FD}{\ensuremath{\EuScript{F}}}

\renewcommand{\H}{\EuScript{H}}

\newcommand{\Arg}{\mathbf{Arg}}

\newcommand{\SLZ}{SL(\ensuremath{2,\mathbb{Z}})}
\newcommand{\PSLR}{PSL(\ensuremath{2,\mathbb{R}})}
\newcommand{\PSLZ}{PSL(\ensuremath{2,\mathbb{Z}})}

\newcommand{\SLR}{SL(\ensuremath{2,\mathbb{R}})}

\newcommand{\PGLR}{PGL(\ensuremath{2,\mathbb{R}})}
\newcommand{\GLR}{GL(\ensuremath{2,\mathbb{R}})}

\newcommand{\Id}{\text{Id}}
\newcommand{\sgn}{\text{sgn}}
\newcommand{\MAS}{\ensuremath{\EuScript{M}}}
\newcommand{\Mas}[1]{\EuScript{M}(#1)}

\newcommand{\Tm}{\ensuremath{S}}
\newcommand{\E}{\ensuremath{T}}

\renewcommand{\mod}{\  \text{mod} \ }
\newcommand{\T}[2][v]{T_{#2,k}^{#1}}

\newcommand{\Th}[2][v]{\Theta_{#2,k}^{#1}}

\newcommand\prefixtext[1]{%
\ifvmode\else\\\@empty\fi
\noalign{%
\penalty0%
\vbox{\mathstrut}%
\penalty10000%
\vskip-\baselineskip
\penalty10000%
\vbox to 0pt{%
\normalbaselines
\ifdim\linewidth=\columnwidth
\else
\parshape\@ne
\@totalleftmargin\linewidth
\fi
\vss
\noindent#1\par}%
\penalty10000%
\vskip-\baselineskip}%
\penalty10000}

\def\mytoday{\number\year --\number\month--\number\day}
\def\blfootnote{\xdef\@thefnmark{}\@footnotetext} 

\author{Fredrik Str\"omberg}

\address{
Institut f\"ur Theoretische Physik \\
TU Clausthal   \\
Abteilung Statistische Physik und Nichtlineare Dynamik \\
Arnold-Sommerfeld-Stra{\ss}e-6\\
38678 Clausthal-Zellerfeld\\
Germany}
\email{fredrik.stroemberg@tu-clausthal.de}
\date{\mytoday}


\makeatother
\begin{document}

\title{Computation of Maass Waveforms with Non-Trivial Multiplier Systems}

\date{\mytoday}

\begin{abstract}
The aim of this paper is to describe efficient algorithms for computing
Maass waveforms on subgroups of the modular group $\PSLZ$ with general
multiplier systems and real weight. A selection of numerical results
obtained with these algorithms is also presented. Certain operators
acting on the spaces of interest are also discussed. The specific
phenomena that were investigated include the Shimura correspondence
for Maass waveforms and the behavior of the weight-$k$ Laplace spectra
for the modular surface as the weight approaches $0$. 
\end{abstract}
\maketitle

\section{Introduction and notation}

The purpose of the present paper is to present computational methods
and experimental results for Maass waveforms with general real weight
and general multiplier systems.

The classical theory of holomorphic automorphic forms was developed
in the setting of (even) integer weight. This was motivated both from
a geometrical point of view and by number theoretical applications
(e.g. the study of modular forms related to the modular invariant
$j$ and the discriminant $\Delta$). The problem of finding the number
of representations of an integer as a sum of a fixed number of squares
(cf.~e.g.~\cite{mordell:representations} and \cite{hardy:on_toe_representaation})
was successfully treated by using the theory of half-integral weight
forms (e.g.~$\theta$-series). This motivated Petersson to develop
a theory of automorphic forms and multiplier systems of arbitrary
real weight \cite{petersson:30} and later also complex weight \cite[I-IV]{MR0028964}. 

Many applications of modular forms use Hecke operators as a principal
tool and Wohlfahrt developed a theory of Hecke-like operators in arbitrary
real weight \cite{MR0106888} (cf. also \cite{van-Lint:HeckeOperators:MR0090616})
. 

Now recall the definition of a Maass waveform as a real-analytic square-integrable
eigenfunction of the Laplacian on a Riemann surface of finite volume
with constant negative curvature $-1$. These were introduced by Maass
for zero weight in \cite{maass:49}. In \cite{MR0065583} the theory
of waveforms for general weights was developed by using the so called
lowering and raising operators, which send a waveform of a given weight
to one with a smaller or larger weight. 

In a slightly different setting Selberg \cite[pp.\ 82--83]{MR0088511}
observed that considering the invariant differential operators on
the space $\H\times S^{1}$ with representation $\chi$ and separating
variables lead to real analytic eigenfunctions of the form $f\left(z,\phi\right)=y^{\frac{k}{2}}F\left(z\right)e^{-ik\phi}$
where $F\left(z\right)$ is a holomorphic modular form of integer
weight $k$ and character $\chi$. 

For an overview of the spectral theory of real analytic modular forms
with arbitrary real weights and multiplier systems see e.g. Maass
\cite{maass:modular_functions} and Roelcke \cite{MR0243062}. More
recent work can be found in e.g.~\cite{MR0384702,MR0419364,MR0491511},
\cite{hejhal:lnm1001,MR803365} and \cite{MR808915,MR840831,bruggeman:varying_weightIII}. 

It is also worth mentioning that the recent interest in Maass waveforms
as representative objects for studying quantum chaos also applies
to real weights. If a weight zero waveform corresponds to a quantum
mechanical particle moving freely on a Riemann surface then a weight
$k$ waveform represents a similar particle moving in a constant magnetic
field with field strength proportional to $k$.

\subsection{Algorithms}

Previously published algorithms for computing Maass waveforms on cofinite
Fuchsian groups have been restricted to groups with one cusp, e.g.
the full modular group or Hecke triangle groups (see e.g. \cite{MR803370,MR936998,hejhal:92,hejhal:99_eigenf,hejhal:calc_of_maass_cusp_forms,muhlenbruch:03}).
By adjoining certain elements it is also possible to bring certain
Hecke congruence subgroups to the one-cusp case, cf.~e.g.~\cite{farmer-lemurell,farmer-lemurell:2,helen:deform_published}. 

In addition, with the exception of the (somewhat crude) computations
in \cite{muhlenbruch:03}, only trivial multiplier system and zero
(or even integer for the holomorphic case in \cite{hejhal:99_eigenf})
weight has been considered. 

The most stable of the algorithms cited above is the one based on
{}``implicit automorphy'' by Hejhal (as detailed in e.g.~\cite{hejhal:calc_of_maass_cusp_forms})
which admits generalizations first of all to remarkably large spectral
parameter (cf.~\cite{then:02}) and a further advantage is that it
does not depend on any underlying arithmetical properties (i.e. Hecke
operators). 

In \cite{stromberg:04:1} this algorithm was generalized to groups
with several cusps, e.g.~Hecke congruence subgroups $\Gamma_{0}(N)$
with non-trivial Dirichlet characters and in and \cite[ch.\ 3]{stromberg:thesis}
we also considered general subgroups of the modular group. Recently,
in \cite{andreas:effective_comp}, this algorithm in combination with
other theoretical methods was used to show the existence of certain
Maass waveforms close to the tentative waveforms produced by the algorithm. 

The aim of the present paper (which is based on \cite[ch.\ 2]{stromberg:thesis})
is to demonstrate how to extend the algorithm to general multiplier
systems and arbitrary real weights. 

The first section contains a brief review of the basic theory of multiplier
systems and then we will introduce the notion of Maass waveforms in
this context. We will also provide some details on the different operators
that act on the space of Maass waveforms with non-trivial multiplier
system. While being of interest in themselves, these operators can
also be used in combination with other tests of reliability and accuracy
of the algorithm. 

In section \ref{sec:Some-Computational-Remarks} we will give the
specifics of how the algorithm is modified and in the last section
we present a selection of the results which has been obtained with
the described method.

\subsection{\label{sub:Summary-of-Notation}Summary of notation}

We will use the notation $e(x)=e^{2\pi ix}$ and for a complex number
$z$ we always use the principal branch of the argument, $-\pi<\Arg z\le\pi$.

Let $\H=\left\{ z=x+iy\,|\, y>0\right\} $ be the upper half-plane
equipped with the hyperbolic line- and area-elements $ds^{2}=\frac{|dz|^{2}}{y^{2}}$
and $d\mu=\frac{dxdy}{y^{2}}$ respectively. The boundary of $\H$
is $\partial\H=\mathbb{R}\cup\left\{ \infty\right\} $. The isometry
group of $\H$ is identified with $\PGLR=\GLR/\left\{ \pm\Id\right\} $,
where $\GLR$ is the group of invertible two-by-two matrices with
real elements and $\Id=\left(\begin{smallmatrix}1 & 0\\
0 & 1\end{smallmatrix}\right)$. For $\gamma=\left(\begin{smallmatrix}a & b\\
c & d\end{smallmatrix}\right)\in\GLR$ and $z\in\H$ we define an action by \[
\gamma z=\left\{ \begin{array}{ccc}
\frac{az+b}{cz+d}, & \textrm{if} & ad-bc>0,\\
\\\overline{\frac{az+b}{cz+d}}, & \textrm{if} & ad-bc<0.\end{array}\right.\]
The subgroup of orientation-preserving isometries of $\H$ is given
by $\PSLR=\SLR/\left\{ \pm\Id\right\} $ where $\SLR$ is the subgroup
of $\GLR$ consisting of matrices with determinant $1$. For any subgroup
$\Gamma\subseteq\textrm{\PGLR,}$ we use $\overline{\Gamma}$ to denote
the inverse image of $\Gamma$ in $\GLR$. Note that this forces $-\Id\in\overline{\Gamma}$.
We are mainly interested in \emph{Fuchsian groups,} i.e.~discrete
subgroups of $\PSLR$. Of particular interest is the subgroup consisting
of matrices with integer entries, the \emph{modular group,} $\PSLZ$.
We are also interested in the so-called Hecke congruence subgroups,
$\Gamma_{0}(N)=\left\{ \left(\begin{smallmatrix}a & b\\
c & d\end{smallmatrix}\right)\in\PSLZ\,|\, c\equiv0\mod N\right\} ,$ defined for any positive integer $N$ (note that $\Gamma_{0}(1)=\PSLZ$). 

We say that an element $\gamma$ of $\PSLR$ is \emph{elliptic}, \emph{parabolic}
or \emph{hyperbolic} if the absolute value of the trace of the associated
matrix is smaller than, equal to or greater than $2$ respectively,
or equivalently, if $\gamma$ has one fixed point in $\H$, one (double)
fixed point in $\partial\H$ or two (different) fixed points in $\partial\H$.
Fixed points of parabolic elements are called \emph{cusps}. 

If $\Gamma\subset\PSLR$ is a finitely generated Fuchsian group we
identify the set of $\Gamma$-orbits with a connected subset of $\H,$
$\mathcal{F=}\Gamma\backslash\H$, a \emph{fundamental domain} of
$\Gamma$. If $\Gamma$ has a set $p_{1},\ldots,p_{\kappa}$ of inequivalent
cusps then the set $\mathcal{F}$ will meet $\partial\H$ at $\kappa$
inequivalent points which we will also denote by $p_{1},\ldots,p_{\kappa}$
and we usually abuse the notation and call these points the cusps
of $\Gamma$ (or of $\mathcal{F}$). By conjugation we may always
assume $p_{1}=i\infty$. Corresponding to a cusp $p_{j}$ of $\Gamma$
we use $\Tm_{j}$ to denote a parabolic generator of $\Gamma_{p_{j}}$
- the subgroup of $\Gamma$ which fixes $p_{j}$. We also choose a
\emph{cusp normalizer,} $\sigma_{j}\in\PSLR$, with the property that
$\sigma_{j}\left(\infty\right)=p_{j}$ and $\sigma_{j}\Tm\sigma_{j}^{-1}=\Tm_{j},$
where $\Tm=\left(\begin{smallmatrix}1 & 1\\
0 & 1\end{smallmatrix}\right)$ is the parabolic generator of $\PSLZ$ (the other generator being
$\E=\left(\begin{smallmatrix}0 & -1\\
1 & 0\end{smallmatrix}\right)$. The map $\sigma_{j}$ is uniquely determined up to a translation
$\Tm^{k}$. If $\mathcal{F}$ meet $\partial\H$ at the points $q_{i},$
$1\le j\le\kappa_{0}$ we fix a set of maps $U_{i}\in\Gamma$ such
that $U_{i}q_{i}=p_{j}$ where $p_{j}$ is the unique cusp equivalent
to $q_{i}$.

\section{Multiplier systems}

\subsection{Introduction}

We will give a brief introduction to multiplier systems, for more
extensive treatments see \cite[pp.\ 331-338]{hejhal:lnm1001}, \cite{petersson:37:analytischen},
or \cite[pp.\ 70-87]{rankin:mod}. 

Let $\Gamma$ be a Fuchsian group and $m$ an even integer. Classically,
a function $\varphi,$ meromorphic on the upper half-plane $\H$,
which satisfies \begin{equation}
\varphi(Az)=\Theta_{A}(z;m)\varphi(z)=(cz+d)^{m}\varphi(z),\quad\forall A=\left(\begin{array}{cc}
a & b\\
c & d\end{array}\right)\in\Gamma,\label{eq:autom_form_def}\end{equation}
is called an \emph{automorphic form} of weight $m$ for $\Gamma$.
The function \[
\Theta_{A}(z;m)=\left(cz+d\right)^{m},\, A=\left(\begin{array}{cc}
a & b\\
c & d\end{array}\right)\in\Gamma\]
is said to be an \emph{automorphy factor} on $\Gamma$\emph{.} The
classical theory of automorphic forms is well-known; for instance,
if $m=2,$ then the automorphic forms can be identified with the meromorphic
differential forms of degree $1$ on the orbifold (classical Riemann
surface) $\Gamma\backslash\H$. 

We observe that, for even $m,$ the number $(cz+d)^{m}$ is uniquely
defined and the automorphy factor $\Theta_{A}(z;m)$ in (\ref{eq:autom_form_def})
clearly satisfies \[
\Theta_{A}(Bz;m)\Theta_{B}(z;m)=\Theta_{AB}(z;m).\tag{*}\]

To generalize these notions to arbitrary real $m$, there needs to
be a choice of branch of the argument, and to make certain everything
is well-defined, we have to introduce the notion of a multiplier system.

\begin{defn}
For any real number $m$ define \[
j_{A}(z;m)=e^{im\Arg(cz+d)}=\frac{(cz+d)^{m}}{|cz+d|^{m}}=\left(\frac{cz+d}{c\overline{z}+d}\right)^{\frac{m}{2}},\, A=\left(\begin{array}{cc}
a & b\\
c & d\end{array}\right)\in\SLR.\]
To adapt the relation ({*}), we also write\begin{equation}
\sigma_{m}(A,B)=j_{A}(Bz;m)j_{B}(z;m)j_{AB}(z;m)^{-1}.\label{eq:chain_rule_for_j_incl_sigma}\end{equation}
 It is clear that for integer $m$, $\sigma_{m}(A,B)=1,$ but it can
also be shown (cf.~\cite[\S2, pp.\ 42--50]{petersson:37:analytischen})
that the only values which $\sigma_{m}$ can take are $1$ and $e^{\pm2\pi im}$.
\end{defn}

\begin{defn}
\label{def:multiplier-system} $v:\overline{\Gamma}\rightarrow S^{1}=\left\{ z\,|\,\left|z\right|=1\right\} $
is said to be a \emph{multiplier system} of weight $m$ on $\overline{\Gamma}$
if
\begin{itemize}
\item $v(-I)=e^{-\pi im},$ and
\item $v(AB)=\sigma_{m}(A,B)v(A)v(B),$ $\forall A,B\in\overline{\Gamma}.$ 
\end{itemize}
\end{defn}
Observe that $v$ can be regarded equally well as a multiplier system
of any weight $m'\equiv m\mod2$. The question of whether there exist
multiplier systems of a given weight and on a given group is most
easily answered by the following proposition (cf. \cite[ Prop.\ 2.1, p.\ 333]{hejhal:lnm1001}). 

\begin{prop}
\label{pro:multiplier_exists}Given $v:\overline{\Gamma}\rightarrow S^{1}$
and $m\in\mathbb{R}.$ The following are equivalent: 
\begin{itemize}
\item $v(T)$ is a multiplier system of weight $m$ on $\overline{\Gamma}.$
\item There exists a function $\varphi\not\equiv0$ on $\H$ which is \emph{}either
$C^{\infty}$ \emph{}or \emph{}meromorphic such that \[
\varphi(Az)=v(A)\varphi(z)(cz+d)^{m},\,\forall A=\left(\begin{array}{cc}
a & b\\
c & d\end{array}\right)\in\overline{\Gamma}.\]
 
\end{itemize}
\end{prop}
A $\varphi$ as above is said to be an automorphic form of weight
$m$ and multiplier system $v$ on $\Gamma$. 

\begin{defn}
\label{def:conjugate-multiplier}Given a multiplier system $v$ on
$\overline{\Gamma}$ and an element $\alpha\in\GLR$ we define a multiplier
system, $v^{\alpha},$ on the group $\alpha^{-1}\overline{\Gamma}\alpha$
by \[
v^{\alpha}(A)=v\left(\alpha A\alpha^{-1}\right)\frac{\sigma_{m}\left(\alpha A\alpha^{-1},\alpha\right)}{\sigma_{m}\left(\alpha,A\right)},\, A\in\alpha^{-1}\overline{\Gamma}\alpha.\]

That this indeed gives a multiplier system on $\alpha^{-1}\overline{\Gamma}\alpha$
is shown in \cite[p.\ 138]{maass:modular_functions}.
\end{defn}
Using Prop. \ref{pro:multiplier_exists}, we will construct the two
most widely used multiplier systems in the following sections. Compare:
\cite[pp.\ 334-337]{hejhal:lnm1001}.

\subsection{The $\eta$ multiplier system\label{sub:The-eta-multiplier}}

\subsubsection{The $\eta$-function}

The Dedekind $\eta-$function is a holomorphic function on $\H,$
defined by \[
\eta(z)=e\left(\frac{z}{24}\right)\prod_{n=1}^{\infty}\left(1-e(nz)\right).\]
 It is clear from the definition that $\eta(z)\neq0$ for $z\in\H$
and that, for each $k\in\mathbb{R},$ $\eta^{2k}$ can be defined
as a holomorphic function on $\H$. (Cf.~\cite[p.\ 205]{rankin:mod}.)
Note that $\eta(z)^{24}$ is the famous Discriminant function $\Delta(z)$.
Cf.~\cite[pp.\ 196-197]{rankin:mod}. It is clear that $\eta\left(z+1\right)=e\left(\frac{1}{24}\right)$
and as for $\Delta\left(z\right)$ it is also possible to express
$\eta\left(z\right)$ as a lacunary Fourier series (cf. \cite[p.\ 18]{muhlenbruch:03}).

\subsubsection{The multiplier system}

It can be proved (cf. Thm. 3.1 and Thm. 3.4 \cite[p.\ 48 and p.\ 52]{apostol})
that $\eta$ satisfies the following functional equations \begin{eqnarray}
\eta\left(\frac{-1}{z}\right) & = & (-iz)^{\frac{1}{2}}\eta\left(z\right),\,\text{and in general, }\nonumber \\
\eta(Az) & = & v_{\eta}(A)(cz+d)^{\frac{1}{2}}\eta(z),\,\forall A=\left(\begin{array}{cc}
a & b\\
c & d\end{array}\right)\in\SLZ.\label{eq:eta_functional_eq}\end{eqnarray}
This functional equation expresses the fact that $\eta$ is an $\SLZ-$automorphic
form of weight $\frac{1}{2}$ and multiplier system given by $v_{\eta}.$
Accordingly the function $\eta^{2k}$ is an $\SLZ$-automorphic form
of weight $k$ and multiplier system given by $v_{\eta}^{2k},$ and
we can use $\eta^{2k}$ in the context of Proposition \ref{pro:multiplier_exists}
to assure the existence of the multiplier system, $v_{\eta,k}=v_{\eta}^{2k}$,
of weight $k$ on $\SLZ$ (and any of its subgroups, e.g. $\Gamma_{0}(N)$).
We have the following explicit formula for $v=v_{\eta}^{2k}$:\begin{equation}
\frac{1}{2\pi i}\log v\left(\left(\begin{array}{cc}
a & b\\
c & d\end{array}\right)\right)=\left\{ \begin{array}{ccc}
\frac{kb}{12}, &  & a=d=1,c=0,\\
k\left(\frac{a+d-3c}{12c}-s(d,c)\right), &  & c>0,\end{array}\right.\label{eq:eta_multiplier_weightk}\end{equation}
and for $c<0$ we use that $v\left(-A\right)=e^{-k\pi i}v\left(A\right)$
(cf.~Def.~\ref{def:multiplier-system}). Here $s(d,c)$ is the Dedekind
sum, \[
s(d,c)=\sum_{n=1}^{c-1}\frac{n}{c}\left(\left(\frac{dn}{c}\right)\right),\]
where $\left(\left(x\right)\right)$ is the saw-tooth function \[
\left(\left(x\right)\right)=\begin{cases}
x-\left\lfloor x\right\rfloor -\frac{1}{2}, & \text{if }x\notin\mathbb{Z},\text{ }\\
0, & \textrm{if }x\in\mathbb{Z},\end{cases}\]
 and $\left\lfloor x\right\rfloor $ is the greatest integer less
than or equal to $x$. Note that if $x$ is not an integer, then $\left\lfloor -x\right\rfloor =-\left\lfloor x\right\rfloor -1$
so $\left(\left(-x\right)\right)=-\left(\left(x\right)\right),$ and
hence $s(-d,c)=-s(d,c)$ if $gcd(d,c)=1$. 

\begin{rem}
It is also possible to express the eta multiplier explicitly without
Dedekind sums but using extended quadratic residue symbols instead.
We have the following formulas from Knopp \cite[p.\ 51]{knopp:modular}
or van Lint \cite[Thm.\ 3]{van-Lint:dedekindEta:MR0103287}: \begin{equation}
v_{\eta}\left(\left(\begin{array}{cc}
a & b\\
c & d\end{array}\right)\right)=\begin{cases}
\left(\frac{c}{d}\right)e\left(\frac{1}{24}\left[\left(a+d\right)c-bd\left(c^{2}-1\right)+3d-3-3cd\right]\right), & c>0,\,\textrm{even},\\
\left(\frac{d}{c}\right)e\left(\frac{1}{24}\left[\left(a+d\right)c-bd\left(c^{2}-1\right)-3c\right]\right), & c>0,\,\textrm{odd.}\end{cases}\label{eq:Knopps_eta_formula}\end{equation}
(Note that the symbols $\left(\frac{c}{d}\right)_{*}$ and $\left(\frac{d}{c}\right)^{*}$
of \cite{knopp:modular,van-Lint:dedekindEta:MR0103287} agree with
our symbols in these two cases.)
\end{rem}

\begin{rem}
\label{rem:many_mult_system}It is known that, for each $k\in\mathbb{R}$,
there exist exactly $6$ different multiplier systems of weight $k$
on $\PSLZ$ (cf.~\cite[\S 3.4, pp.\ 83, 206 ]{rankin:mod} or \cite[Thm.\ 19, p.\ 132]{maass:modular_functions}).
We will denote these by $v_{\eta,k}^{(r)}=v_{\eta}^{2(k+r)},$ where
$r\in\{0,2,4,6,8,10\}$. Compare \cite[eq.\ (6.4.7)]{rankin:mod};
one knows, of course, that $v_{\eta}^{24}=1$. 

When dealing with the modular group and weight $k$, it is sufficient
to consider only the multiplier system $v_{\eta,k}^{(0)}=v_{\eta,k}=v_{\eta}^{2k}$
(for reasons to be discussed later in Section \ref{sub:Maass-operators-and-symmetry}). 
\end{rem}

\subsection{The $\theta$ multiplier system\label{sub:The-theta-multiplier}}

On any subgroup of $\PSLZ$, we can always use the $\eta$-multiplier
system, but in general, on subgroups of $\PSLZ$, there are also other
multiplier systems available. In particular, on $\Gamma_{0}(4)$,
there is a multiplier system of weight $\frac{1}{2}$ which is interesting
from an arithmetical point of view. 

It is well-known (cf. \cite{shimura:73:half_integral} or \cite{Iwaniec:topics})
that the Jacobi theta function \[
\theta(z)=\sum_{-\infty}^{\infty}e(n^{2}z),\, z\in\H,\]
is automorphic on $\Gamma_{0}(4)$ with weight $k=\frac{1}{2}$ and
can be used to define a multiplier system on $\Gamma_{0}(4).$ Using
the Poisson summation formula one can prove (cf. \cite[pp.\ 72--75]{gunning}
or \cite[pp.\ 167--168]{Iwaniec:topics}) that the theta function
satisfies: \begin{equation}
\theta\left(\frac{-1}{2z}\right)=(-iz)^{\frac{1}{2}}\theta\left(\frac{z}{2}\right),\label{eq:inversion_formula_theta}\end{equation}
 and one can also prove the general formula (cf. \cite[Thm.\ 10.10, p.\ 177]{Iwaniec:topics}
or \cite[p.\ 447]{shimura:73:half_integral}): \begin{eqnarray}
\theta(Az) & = & v_{\theta}\left(A\right)(cz+d)^{\frac{1}{2}}\theta(z),\, A=\left(\begin{array}{cc}
a & b\\
c & d\end{array}\right)\in\Gamma_{0}(4).\label{eq:functional_eq_theta}\end{eqnarray}
 The multiplier $v_{\theta}\left(A\right)$ can be expressed explicitly
as\[
v_{\theta}\left(A\right)=\bar{\epsilon}_{d}\left(\frac{c}{d}\right),\]
where $\epsilon_{d}=1$ if $d\equiv1\mod4$ and $\epsilon_{d}=i$
if $d\equiv-1\mod4,$ and $\left(\frac{c}{d}\right)$ denotes the
extended quadratic residue symbol defined as the traditional Jacobi
symbol if $0<d\equiv1\mod2$ and extended by \[
\left(\frac{c}{d}\right)=\frac{c}{|c|}\left(\frac{c}{-d}\right),\, c\neq0,\]
and \[
\left(\frac{0}{d}\right)=\begin{cases}
1 & \text{if }d=\pm1,\\
0 & \text{otherwise.}\end{cases}\]
 For the sake of completeness we also use the traditional Kronecker
extension, i.e. we define \begin{equation}
\left(\frac{\vphantom{2}c}{2}\right)=\left(\frac{2}{c}\right).\label{eq:kronecker_c2eq2c}\end{equation}
 One can verify that our symbol $\left(\frac{\cdot}{\cdot}\right)$
satisfies reciprocity relations similar to the usual ones:

\begin{prop}
\label{pro:reciprocity} Suppose that $c,d\in\mathbb{Z}$ are odd
and $c\ne0$. Then we have:

\begin{eqnarray*}
\left(\frac{-1}{d}\right) & = & \left(-1\right)^{\left(\frac{d-1}{2}\right)},\\
\left(\frac{c}{d}\right) & = & \begin{cases}
\left(\frac{d}{c}\right)\left(-1\right)^{\left(\frac{d-1}{2}\right)\left(\frac{c-1}{2}\right)}, & d,\,\text{or }c>0,\\
-\left(\frac{d}{c}\right)\left(-1\right)^{\left(\frac{d+1}{2}\right)\left(\frac{c+1}{2}\right)}, & d,\,\text{and }c<0.\end{cases}\end{eqnarray*}

\end{prop}
\begin{rem}
Relations (\ref{eq:inversion_formula_theta}) and (\ref{eq:functional_eq_theta})
can also be proved using the corresponding relations (\ref{eq:eta_functional_eq})
for $\eta$ and the following relation between the $\eta$ and the
$\theta$ functions (cf. \cite[p.\ 177]{Iwaniec:topics} or \cite[Thm.\ 12, p.\ 46]{knopp:modular}):
$\theta(z)=\frac{\eta²\left(\frac{z+1}{2}\right)}{\eta(z+1)}.$
\end{rem}

\section{Maass waveforms}

The slash-operator $f_{|A}(z)=f(Az)$ can be extended to an operator
of weight $k$ as: \[
f_{|[k,A]}(z)=f(Az)j_{A}(z;k)^{-1},\]
 and the natural analog of the Laplace-Beltrami operator, $\Delta,$
which is invariant under this action is the weight-$k$ Laplacian:
\[
\Delta_{k}=\Delta-iyk\frac{\partial}{\partial x}=y^{2}\left(\frac{\partial^{2}}{\partial x^{2}}+\frac{\partial^{2}}{\partial y^{2}}\right)-iyk\frac{\partial}{\partial x}\,.\]
If $\Gamma$ is a Fuchsian group we define the space $\Mas{\Gamma,v,k,\lambda}$
consisting of Maass waveforms on $\Gamma,$ of weight $k$, multiplier
system $v$ and eigenvalue $\lambda,$ as the space of functions which
satisfy the following conditions:

\begin{lyxlist}{00.00.0000}
\item [{1)}] $f_{|[A,k]}(z)=v(A)f(z),$ $\forall A=\left(\begin{array}{cc}
a & b\\
c & d\end{array}\right)\in\overline{\Gamma},$
\item [{2)}] $\Delta_{k}f+\lambda f=0,$ and 
\item [{3)}] $\int_{\FD}|f|^{2}d\mu<\infty.$
\end{lyxlist}
Observe that condition 1) is equivalent to 

\begin{lyxlist}{00.00.0000}
\item [{1')}] $f(Az)=v(A)j_{A}(z;k)f(z),$ $\forall A=\left(\begin{array}{cc}
a & b\\
c & d\end{array}\right)\in\overline{\Gamma}.$
\end{lyxlist}
 \emph{}

\emph{For purposes of the computational work to be described in this
paper, we shall be content to restrict ourselves to cases where $\lambda>\frac{1}{4}$.
(Cf.~also here para. 4 of sect. \ref{ref:small_eigenv} below.)}

Instead of the Bessel equation in the case of weight $0,$ condition
2) above gives us the Whittaker equation, and using the method of
separation of variables gives us Whittaker functions instead of the
K-Bessel functions at weight $0$ (for complete details see \cite[Chap.\ 9]{hejhal:lnm1001}).
Since $f(x+iy)$ is no longer periodic in $x,$ but instead satisfies
$f(z+1)=v(S)f(z)=e(\alpha)f(z),$ with $\alpha\in[0,1)$, the Fourier
series of $f$ can (\cite[pp.\ 26, 348, 420(19)]{hejhal:lnm1001})
be written as \begin{eqnarray}
f(z) & = & \sum_{\underset{{n+\alpha\neq0}}{-\infty}}^{\infty}\frac{c(n)}{\sqrt{|n+\alpha|}}W_{sgn(n+\alpha)\frac{k}{2},iR}(4\pi|n+\alpha|y)e\left((n+\alpha)x\right),\label{eq:fourierseries_f}\end{eqnarray}
where $W_{l,\mu}(x)$ is the Whittaker function in standard notation
(cf.~\cite[vol.\ I, p.\ 264]{erdelyi:53}) and $R$ is the usual
spectral parameter, $\lambda=\frac{1}{4}+R^{2}$. One notes here that
$W_{0,iR}(x)=\pi^{-\frac{1}{2}}x^{\frac{1}{2}}K_{iR}\left(\frac{x}{2}\right).$
For $k=0,$ the expansion above thus reduces to usual Fourier expansion
with $2y^{\frac{1}{2}}K_{iR}(2\pi|n+\alpha|y)$ as in \cite[p.\ 26, prop.\ 4.12]{hejhal:lnm1001}. 

If we have more than one cusp we define functions $f_{j}$ related
to $f$ at each cusp, $p_{j},$ of $\Gamma$ by using the cusp normalizing
maps $\sigma_{j}$ from section \ref{sub:Summary-of-Notation} and
setting $f_{j}(z)=f_{|[\sigma_{j},k]}(z)=j_{\sigma_{j}}(z;k)^{-1}f(\sigma_{j}z).$
It is easy to see that\[
f_{j}(z+1)=v(S_{j})f_{j}(z)=e(\alpha_{j})f_{j}(z),\]
with $\alpha_{j}\in[0,1)$ (cf. \cite[p. 41]{Iwaniec:topics}). Thus
the Fourier series of $f$ at the cusp $j$ can be written as \begin{eqnarray}
f_{j}(z) & = & \sum_{n=-\infty}^{\infty}\frac{c_{j}(n)}{\sqrt{|n+\alpha_{j}|}}W_{sgn(n+\alpha_{j})\frac{k}{2},iR}(4\pi|n+\alpha_{j}|y)e\left((n+\alpha_{j})x\right).\label{eq:fourierseries_fj}\end{eqnarray}
 As in the case of weight $0$ and Dirichlet character, we say that
the cusp number $j$ is \emph{open} or \emph{singular} if $\alpha_{j}=0$
and \emph{closed} if $\alpha_{j}\ne0.$ If all cusps of $\Gamma$
are singular for the multiplier system $v$ we say that $v$ is a
\emph{singular multiplier system for} $\Gamma$. 

\begin{rem}
Observe that, for the eta-multiplier on $\PSLZ$ and weight $k$,
we have $\alpha=\alpha_{1}=\frac{k}{12}.$ 
\end{rem}

\subsection{Decomposition of the discrete spectrum}

It is known (see for example \cite{bruggeman:94} or \cite{hejhal:lnm1001})
that closed cusps (i.e. $v(S_{j})\ne1)$ do not contribute to the
continuous spectrum, and if all cusps are closed there is only the
discrete part of the spectrum left and this is spanned by the Maass
waveforms. We also know (see \cite[p.\ 385]{hejhal:lnm1001}) that
on the modular group with weight $k$ the smallest eigenvalue is \[
\lambda_{min}=\frac{|k|}{2}\left(1-\frac{|k|}{2}\right),\]
or larger. In the case of $\PSLZ$ and $k\ge0,$ $F(z)=y^{\frac{k}{2}}\eta(z)^{2k}$
has eigenvalue equal to $\lambda_{min}$. 

In this paper, any eigenvalues $\lambda\in\left[\lambda_{min},\frac{1}{4}\right]$
will be regarded as exceptional. The non-exceptional eigenvalues thus
satisfy $\frac{1}{4}<\lambda_{0}\le\lambda_{1}\le\cdots\le\lambda_{n}\rightarrow\infty$.
One can obtain lower bounds for the eigenvalue $\lambda_{0}$ (see
\cite[p.\ 183]{bruggeman:94}), but in light of the numerical experiments
in Section \ref{sub:Varying-weight} they are not very effective (cf.~Figure
\ref{fig:Weights_0.1to6_compare}).

\section{Operators }

\subsection{Conjugation and reflection\label{sub:operators_Conjugation}}

Let $J$ and $K$ denote the reflection, $Jz=-\overline{z}$ and conjugation,
$Kz=\overline{z}$. Then $J$ and $K$ act as involutions on the space
of Maass waveforms via the operations \begin{eqnarray*}
Kf & = & f_{|K}(z)=\overline{f(z),}\\
Jf & = & f_{|J}(z)=f(-\overline{z}).\end{eqnarray*}
 It follows from the definition of the action of $\GLR$ on $\H$
(cf.~page \pageref{sub:Summary-of-Notation}) that we can use the
matrix $J=\left(\begin{smallmatrix}1 & 0\\
0 & -1\end{smallmatrix}\right)$ in $\GLR$ to represent the operator $J$. For $A=\left(\begin{smallmatrix}a & b\\
c & d\end{smallmatrix}\right)$ we define \[
A^{*}=JAJ^{-1}=\left(\begin{array}{cc}
a & -b\\
-c & d\end{array}\right),\]
and then $A^{**}=A,$ and $A(z)_{|K}=-A^{*}(z),$ meaning that also
$-\overline{A(z)}=A^{*}(-\bar{z})$. 

\begin{rem}
\label{rem:f_conj_inv} It is easy to verify that if $f\in\Mas{\Gamma_{0}(N),v,k,\lambda},$
then $Kf\in\Mas{\Gamma_{0}(N),\overline{v},-k,\lambda}$ and $Jf\in\Mas{\Gamma_{0}(N),v^{*},-k,\lambda},$
where $v^{*}$ is the multiplier system determined by\[
v^{*}(A)=v(A^{*})\cdot\begin{cases}
1, & c\neq0,\\
e^{\pi ik(1-sgn(d))}, & c=0,\end{cases},\,\textrm{for}\, A=\left(\begin{array}{cc}
a & b\\
c & d\end{array}\right).\]
 Of particular interest is the involution obtained by combining $J$
and $K,$ i.e. \[
KJf(z)=f_{|JK}(z)=\overline{f(-\bar{z}).}\]
It is easily seen that if $f$ has Fourier coefficients $c_{j}(n)$,
then (by \cite[p.\ 265(8)]{erdelyi:53}) $f_{|JK}$ has Fourier coefficients
$\overline{c_{j}(n)}$, and we thus would like to have $f$ and $f_{|JK}$
belonging to the same space (i.e. transform according to the same
multiplier system), since then we can assume that the Fourier coefficients
are real. 
\end{rem}
It is clear that if $f\in\Mas{\Gamma_{0}(N),v,k,\lambda},$ then $KJf\in\Mas{\Gamma_{0}(N),\overline{v^{*}},k,\lambda}$
so we are left to see whether $\overline{v^{*}}=v$ or not. By using
the explicit formulas one can verify that indeed $\overline{v^{*}}=v$
for both $v_{\theta}$ and $v_{\eta}$ (see \cite[pp.\ 66-68]{stromberg:thesis}
for details) and we arrive at the following proposition. 

\begin{prop}
\label{pro:real_fouriercoefficients}If $v$ is either the $\eta$-
or the $\theta$-multiplier system (in the latter case $4|N$) then
a basis $\left\{ g_{1},\ldots,g_{m}\right\} $ of $\Mas{\Gamma_{0}(N),v,k,\lambda}$
can be chosen so that each $g_{j}$ can be expanded in a Fourier series
at $\infty$ with real coefficients. 
\end{prop}
\begin{proof}
We have seen that for both the theta and the eta multiplier systems
the product $KJ$ is a conjugate-linear involution of the space $\Mas{\Gamma_{0}(N),v,k,\lambda},$
and hence we can assume that any $f\in\Mas{\Gamma_{0}(N),v,k,\lambda}$
is an eigenfunction of $KJ$ with eigenvalue $\epsilon,$ where $|\epsilon|=1.$
Note that if $f(z)$ has a Fourier series expansion as above with
Fourier coefficients $c(n)$ then $f_{|KJ}$ has Fourier coefficients
$\overline{c\left(n\right)}$ and hence $\overline{c(n)}=\epsilon c(n).$
Finally we observe that if $\epsilon=e^{i\theta}$we can look at the
function $g=e^{i\frac{\theta}{2}}f$ which then satisfies $KJg=e^{-i\frac{\theta}{2}}KJf=e^{-i\frac{\theta}{2}}e^{i\theta}f=g$.
After proper normalization it is thus no restriction to assume that
the eigenvalue of $KJ$ is $\epsilon=1,$ and that the Fourier coefficients
are real.
\end{proof}
\begin{rem}
For the sake of completeness it should be remarked that in general
one can not simultaneously take Fourier coefficients at cusps other
than $\infty$ to be real (cf. the next subsection where we introduce
the map $\omega_{N}$, which is a cusp normalizing map for the cusp
at $0$ and which has eigenvalues $\pm i^{-k}$). 
\end{rem}

\subsection{The involution $\tau_{N}$ \label{sub:The-involution_omegaN}}

As in the case of zero weight (e.g. \cite{stromberg:04:1}) we define
$\omega_{N}z=\frac{-1}{Nz},$ or equivalently $\omega_{N}=\left(\begin{smallmatrix}0 & \frac{-1}{\sqrt{N}}\\
\sqrt{N} & 0\end{smallmatrix}\right).$ We know that $\omega_{N}$ is an involution of $\Gamma_{0}(N),$
i.e. $\Gamma_{0}(N)=\omega_{N}\Gamma_{0}(N)\omega_{n}^{-1},$ but
the question is how it relates to the weight and multiplier system. 

If $f\in\Mas{\Gamma_{0}(N),v,k,\lambda}$ it is easy to see that $f{}_{|[k,\omega_{N}]}\in\Mas{\Gamma_{0}(N),v^{\omega_{N}},k,\lambda},$
and it is also easy to verify that $v^{\omega_{N}}(T)=v\left(\omega_{N}T\omega_{N}^{-1}\right).$
We also have \begin{eqnarray*}
f_{|[k,\omega_{N}]|[k,\omega_{N}]}(z) & = & j_{\omega_{N}}(z;k)^{-1}f\left(\omega_{N}z\right)_{|[k,\omega_{N}]}=j_{\omega_{N}}(z;k)^{-1}j_{\omega_{N}}(\omega_{N}z;k)^{-1}f\left(\omega_{N}^{2}z\right)\\
 & = & e^{-ikArg(\sqrt{N}z)}e^{-ikArg(-1/\sqrt{N}z)}f\left(z\right)=e^{-i\pi k}f(z),\end{eqnarray*}
and hence if we define the operator $\tau_{N}$ by \begin{eqnarray*}
\tau_{N}f(z) & = & e^{ik\frac{\pi}{2}}f_{|[k,\omega_{N}]}(z),\end{eqnarray*}
we have that \[
\tau_{N}^{2}=Id.\]
 To show that $\tau_{N}$ is a linear involution, we also have to
verify that $v^{\omega_{N}}=v$. This is easily done in the two cases
$N=1$ together with $v=v_{\eta}$ and $N=4$ together with $v=v_{\theta}$.
And we conclude:

\begin{prop}
If $N=1$ and $v=v_{\eta},$ or $N=4$ and $v=v_{\theta}$ the operator
\[
\tau_{N}:\Mas{\Gamma_{0}(N),v,k,\lambda}\rightarrow\Mas{\Gamma_{0}(N),v,k,\lambda},\]
defined by \[
\tau_{N}f(z)=e^{ik\frac{\pi}{2}}f_{|[k,\omega_{N}]}=e^{-ik(Argz-\frac{\pi}{2})}f\left(\frac{-1}{Nz}\right),\]
is a linear involution, i.e. $\tau_{N}\left(af\right)=a\tau_{N}f$
for all $a\in\mathbb{C}$ and $\tau_{N}^{2}f=f$. Hence it has eigenvalues
$\pm1$.
\end{prop}
\begin{rem}
Note that in the case of $\Gamma_{0}\left(4\right)$ (i.e. $N=4$)
if $\tau_{N}f(z)=\pm f(z),$ then \[
f_{2}=f_{|\omega_{N}}=e^{-ik\frac{\pi}{2}}\tau_{N}f=\pm e^{-ik\frac{\pi}{2}}f,\]
which means that the Fourier coefficients at the cusp at $0$ are
proportional to the coefficients at $i\infty$: \begin{equation}
c_{2}(n)=\pm e^{-i\frac{\pi}{2}k}c_{1}(n)=\pm i^{-k}c_{1}(n).\label{eq:coeff_relation_invol_tau}\end{equation}

\end{rem}

\subsection{The operator $L$ \label{sub:The-operator-L} }

\begin{defn}
\label{rem:operatorL}For $N=4$ and the $\theta$-multiplier system
at weight $\frac{1}{2}$, following Kohnen \cite{MR81j:10030} or
Katok-Sarnak \cite{katok-sarnak} we define the operator $L$ acting
on $\Mas{\Gamma_{0}(4),v_{\theta},\frac{1}{2},R}$ by\[
L=\frac{1}{\sqrt{2}}\tau_{4}T_{4,\frac{1}{2}}^{v_{\theta}},\]
where $T_{4,\frac{1}{2}}^{v_{\theta}}$ is the following (exceptional)
Hecke operator \begin{equation}
T_{4,\frac{1}{2}}^{v_{\theta}}f(z)=\frac{1}{2}\sum_{j\mod4}f\left(\frac{z+j}{4}\right).\label{eq:Hecke_op_w_half_T4}\end{equation}
Cf.~Subsection \ref{sub:Heckeop-Theta-multiplier}. It can be shown
that the operator $L$ is self-adjoint, commutes with $\Delta_{\frac{1}{2}}$
and all Hecke operators $T_{p^{2},\frac{1}{2}}^{v_{\theta}},$ $p\ne2$
(cf.~(\ref{eq:def_Heckeop_theta_mult})) and have the eigenvalues
$1$ and $-\frac{1}{2}$. Cf.~\cite{MR81j:10030} and \cite{MR0562506}.
Suppose that $f(z)$ has Fourier expansions à la (\ref{eq:fourierseries_fj})
at the cusps $p_{1}=\infty,$ $p_{2}=0$ and $p_{3}=-\frac{1}{2}$
with Fourier coefficients $a_{j}(n)$ respectively. By calculations
similar to \cite{MR81j:10030} or \cite{MR81j:10030,biro:00} it can
be shown that $Lf$ has Fourier coefficients $b(n)$ (with respect
to $\infty$) given by \begin{eqnarray}
b(n) & = & \frac{1}{2}\begin{cases}
a_{1}(n)+\left(1+i\right)a_{2}\left(\frac{n}{4}\right), & n\equiv0\mod4,\\
a_{1}(n)+\sqrt{2}a_{3}\left(\frac{n-1}{4}\right)\left(-1\right)^{\frac{n-1}{4}}, & n\equiv1\mod4,\\
-a_{1}(n), & n\equiv2,3\mod4.\end{cases}\label{eq:b_for_op_L}\end{eqnarray}

\label{def:def_of_V+}Let $V^{+}\subseteq\Mas{\Gamma_{0}(4),v_{\theta},\frac{1}{2},\lambda}$
denote the subspace introduced by Kohnen in \cite{MR81j:10030}, i.e.
$V^{+}$consists of all $f$ with Fourier coefficients $a_{1}(n)=0$
for $n\equiv2,3\mod4$. By using (\ref{eq:b_for_op_L}) it is easy
to verify that $V^{+}$ is exactly the eigenspace of $L$ corresponding
to the eigenvalue $1$. Unfortunately the eigenspace $V^{-}$ corresponding
to the eigenvalue $-\frac{1}{2}$ is not as simple to characterize.
However, (\ref{eq:b_for_op_L}) can be used to identify $V^{-}$ by
means of certain relations between coefficients at the cusps at $\infty$,
$0$ and $-\frac{1}{2}$. E.g. if $a_{1}(1)=1$ then $a_{3}(0)=-\sqrt{2}$
and if $a_{1}(1)=0$ then $a_{3}(0)=0$. 

Clearly $T_{4,\frac{1}{2}}^{v_{\theta}}$ does not in general commute
with $L$ but in case $Lf=\lambda f$, $T_{4,\frac{1}{2}}f=\lambda_{4}f$
and $\tau_{N}f=\epsilon f$ (with $\epsilon=\pm1$) then $a_{1}\left(4\right)=\epsilon\sqrt{2}\lambda=-\frac{\epsilon}{\sqrt{2}},\,\sqrt{2}\epsilon$.
This should be compared with the results on newforms at weight zero,
e.g.~\cite[p.\ 147]{atkinlehner} and \cite[p.\ 31 ]{stromberg:thesis}. 
\end{defn}

\subsection{Maass operators }

So far, the operators we have seen act on spaces of Maass waveforms
of a given weight and multiplier system. 

We will show that we may limit the range of weights $k$ to investigate
to $k\in[0,6].$ For this we will use the Maass lowering and raising
operators, $E_{k}^{\pm},$which raise or lower the weight of a Maass
waveform by units of $2$. They are defined by \begin{eqnarray*}
E_{k}^{\pm} & = & iy\frac{\partial}{\partial x}\pm y\frac{\partial}{\partial y}+\frac{k}{2},\end{eqnarray*}
 and using the relation between the Whittaker function and the confluent
hypergeometric function together with the transformation formulas
\cite[p.\  258, (10)]{erdelyi:53} (see also: \cite[p.\ 183 (middle)]{maass:modular_functions}
and \cite[p.\ 302 lines -3 and -1]{MR0232968}) we see that for $k>0$
(here we set $Y=4\pi|n+\alpha|y$ and $n_{\alpha}=n+\alpha$)\begin{alignat}{3}
 & E_{k}^{+}\left(W_{\frac{k}{2},iR}(Y)e(n_{\alpha}x)\right) & = & -W_{\frac{k+2}{2},iR}\left(Y\right)e(n_{\alpha}x), & n_{\alpha}>0,\label{eq:Maass_operators_act_1}\\
 & E_{k}^{-}\left(W_{\frac{k}{2},iR}(Y)e(n_{\alpha}x)\right) & = & -\left(\frac{k(k-2)}{4}+\frac{1}{4}+R^{2}\right)W_{\frac{k-2}{2},iR}\left(Y\right)e(n_{\alpha}x), & n_{\alpha}>0,\label{eq:Maass_operators_act_2}\\
 & E_{k}^{+}\left(W_{-\frac{k}{2},iR}(Y)e(n_{\alpha}x)\right) & = & \left(\frac{k(k+2)}{4}+\frac{1}{4}+R^{2}\right)W_{-\frac{k+2}{2},iR}\left(Y\right)e(n_{\alpha}x), & n_{\alpha}<0,\label{eq:Maass_operators_act_3}\\
 & E_{k}^{-}\left(W_{-\frac{k}{2},iR}(Y)e(n_{\alpha}x)\right) & = & W_{-\frac{k-2}{2},iR}\left(Y\right)e(n_{\alpha}x), & n_{\alpha}<0.\label{eq:Maass_operators_act_4}\end{alignat}
To verify that they respect the weight $k$-action it is easiest to
proceed straight forward but to use the following form of the operators
\begin{eqnarray*}
E_{k}^{+} & = & (z-\overline{z})\frac{\partial}{\partial z}+\frac{k}{2},\\
E_{k}^{-} & = & -(z-\overline{z})\frac{\partial}{\partial\overline{z}}+\frac{k}{2},\end{eqnarray*}
and write $e^{-iArg(cz+d)}=\left(\frac{c\overline{z}+d}{cz+d}\right)^{\frac{1}{2}}$.
Cf. ~\cite[p.\ 178]{maass:modular_functions}.

One can then show that $E_{k}^{\pm}$ maps $\Mas{\Gamma_{0}(N),v,k,\lambda}$
into $\Mas{\Gamma_{0}(N),v,k\pm2,\lambda}$, and that the composition
\[
E_{k\pm2}^{\mp}E_{k}^{\pm}:\Mas{\Gamma_{0}(N),v,k,\lambda}\mapsto\Mas{\Gamma_{0}(N),v,k,\lambda}\]
 is just multiplication by a constant, which is readily seen to be
nonzero anytime $\lambda>\frac{1}{4}$. Hence $E_{k}^{\pm}$ acts
bijectively on the spaces corresponding to non-exceptional eigenvalues,
i.e. they are always bijections for $\lambda>\frac{1}{4}$.

\subsubsection{Maass operators and the symmetry about $k=6$ \ref{sub:Maass-operators-and-symmetry}\label{sub:Maass-operators-and-symmetry} }

First of all, observe that $E_{k}^{\pm}$ only change the weight and
not the multiplier system $v,$ but in view of the remark after Def.~\ref{def:multiplier-system}
it is clear that $v=v_{\eta}^{2k}$ is also a multiplier system of
weight $k+r$ for any $r\in2\mathbb{Z},$ and with the notation $v_{\eta,k}^{(r)}=v_{\eta}^{2(k+r)}$
it is clear that \[
E_{k}^{+}:\Mas{\Gamma_{0}(N),v_{\eta,k}^{(r)},k,\lambda}\rightarrow\Mas{\Gamma_{0}(N),v_{\eta,k+2}^{(r-2)},k+2,\lambda},\]
and \[
E_{k}^{-}:\Mas{\Gamma_{0}(N),v_{\eta,k}^{(r)},k,\lambda}\rightarrow\Mas{\Gamma_{0}(N),v_{\eta,k-2}^{(r+2)},k-2,\lambda}.\]

\emph{Our main purpose is to investigate the eigenvalues of Maass
waveforms on the modular group when the weight and multiplier system
are varied.} That is, we would like to investigate the space $\Mas{\Gamma_{0}(1),v_{\eta,k}^{(r)},k,\lambda}$
for all $k\in\mathbb{R}$ and $r\in\{0,2,4,6,8,10\}.$ Suppose that
$\lambda>\frac{1}{4}$ so the the lowering and raising operators act
bijectively, then using the $(k,r)$ to denote that space $\Mas{\Gamma_{0}(1),v_{\eta,k}^{(r)},k,\lambda}$
and using $\approx$ to denote bijectively corresponding spaces we
have:

\begin{itemize}
\item Since all $v_{\eta,k}^{(r)}$ are 24th-roots of unity we have trivially:
$\left(k,r+12\right)=\left(k,r\right)$.
\item A composition of Maass operators which raises or lowers the weight
by $12$ will preserve the multiplier system. Hence $\left(k+12,r\right)\approx\left(k,r\right)$
and we may assume that $k\in[0,12]$.
\item $K$ is a bijection from $\left(k,r\right)$ to $\left(-k,-r\right)=\left(-k,12-r\right)\approx\left(12-k,12-r\right)$
so there is no restriction to assume $k\in[0,6]$ and $r\in\{0,2,4,6,8,10\}$. 
\item By using the raising operator we see that $\left(k,0\right)\approx\left(k+2,-2\right)\approx\left(k+2,10\right)$
and by repeated use we see that with out loss of generality we can
also assume $r=0$. 
\end{itemize}
We are thus justified in our choice of restricting the investigation
to the spaces $\Mas{\Gamma_{0}(1),v,k,\lambda}$ for $k\in[0,6]$
and $v=v_{\eta,k}^{(0)}=v_{\eta}^{2k}.$

\subsection{\label{sub:Hecke-operators-non-trivial}Hecke operators for non-trivial
multiplier systems}

We know that Hecke operators play an important role in the understanding
of the theory of both modular forms and Maass waveforms at integer
weights. 

For general real weights, the Hecke operators are not important (and
may well be non-existent), but we will begin with a general definition
anyway, and then we will consider two special cases: $\Gamma_{0}(4)$
with the theta multiplier system and weight $\frac{1}{2}$ and $\Gamma_{0}(1)$
with the eta multiplier system and weight $1$. 

The general introductory discussion is based on \cite{strombergsson:heckeoperators}
but the specific example of integer weights (as worked out in detail
in \cite{hecke_operators_weight}) is based on ideas from \cite{MR0106888}
and \cite{van-Lint:HeckeOperators:MR0090616}. 

Let $\Gamma\subset\PSLR$ be cofinite and suppose $v:\overline{\Gamma}$$\rightarrow S^{1}$
is a multiplier system of weight $k\in\mathbb{R}.$ The commensurator
of $\Gamma,$ $\textrm{comm}\left(\Gamma\right),$ in $\PSLR$ is
defined as \[
\textrm{comm}\left(\Gamma\right)=\left\{ \alpha\in\PSLR\,|\,\alpha\Gamma\alpha^{-1}\cap\Gamma\,\textrm{ has finite index in }\Gamma\textrm{ and }\alpha\Gamma\alpha^{-1}\right\} .\]
We know that the Hecke operators are associated with the members of
the commensurator, or actually with the double cosets, $\overline{\Gamma}\alpha\overline{\Gamma},$
for $\alpha\in\textrm{comm}\left(\Gamma\right)$. Fix $\alpha\in\textrm{comm}\left(\Gamma\right).$
It can be shown that if the multiplier system $v$ satisfies \begin{equation}
v(g)=v^{\alpha}(g),\,\forall g\in\overline{\Gamma}\cap\alpha^{-1}\overline{\Gamma}\alpha,\label{eq:multiplier_coset_condition}\end{equation}
then we can define an associated multiplier system, $v_{\alpha},$
on the double coset:\[
v_{\alpha}:\overline{\Gamma}\alpha\overline{\Gamma}\rightarrow S^{1},\]
by setting\[
v_{\alpha}\left(g_{1}\alpha g_{2}\right)=\sigma_{k}\left(g_{1}\alpha,g_{2}\right)\sigma_{k}\left(g_{1},\alpha\right)v\left(g_{1}\right)v\left(g_{2}\right),\]
for all $g_{1},g_{2}\in\overline{\Gamma}.$ It might be the case that
there exists an associated multiplier system of $W=\chi v$, where
$\chi$ is a character on $\overline{\Gamma}\alpha\overline{\Gamma},$
even though there does not exist an associated multiplier system of
$v$ itself. Suppose that $v_{\alpha}$ exists as above, and that
we have $\overline{\Gamma}\alpha\overline{\Gamma}=\cup_{i=1}^{d}\overline{\Gamma}\alpha_{i}$.
We then define the Hecke operator $T_{\alpha,k}^{v}:\MAS\left(\overline{\Gamma},v,k,\lambda\right)\rightarrow\MAS\left(\overline{\Gamma},v,k,\lambda\right)$
by \begin{eqnarray*}
\left(T_{\alpha,k}^{v}f\right)(z) & = & \sum_{i=1}^{d}\overline{v_{\alpha}\left(\alpha_{i}\right)}f_{|[\alpha_{i},k]}(z).\end{eqnarray*}
When $\Gamma=\Gamma_{0}(N)$ (or any congruence subgroup of level
$N$), one usually constructs Hecke operators $T_{n,k}^{v}$ corresponding
to positive integers $n$ in which case $\alpha=\left(\begin{smallmatrix}1 & 0\\
0 & n\end{smallmatrix}\right)$ (cf.~e.g.~\cite{atkinlehner}, \cite[ch.\ 5]{gunning}, \cite[ch.\ 9]{rankin:mod},
\cite{miyake} or \cite{shimura} for more details).

\subsubsection{Hecke Operators for the Theta multiplier System\label{sub:Heckeop-Theta-multiplier}}

Consider the case $\Gamma=\Gamma_{0}(4),$ $k=\frac{1}{2}$ and $v=v_{\theta}.$
Let $n$ be a positive integer and $\alpha=\left(\begin{smallmatrix}1 & 0\\
0 & n\end{smallmatrix}\right).$ Then $g\in\overline{\Gamma_{0}(4)}\cap\alpha^{-1}\overline{\Gamma_{0}\left(4\right)}\alpha$
can be written $g=\left(\begin{smallmatrix}a & nb\\
c/n & d\end{smallmatrix}\right),$ and $\alpha g\alpha^{-1}=\left(\begin{smallmatrix}a & b\\
c & d\end{smallmatrix}\right),$ with $ad-bc=1$ and $c\equiv0\mod4n$. It is easy to verify that
$v^{\alpha}(g)=v(\alpha g\alpha^{-1})$ and $v\left(g\right)=v(\alpha g\alpha^{-1})\left(\frac{n}{d}\right)^{-1},$
and hence $v^{\alpha}(g)=v(g)$ if and only if $\left(\frac{n}{d}\right)=1.$
By (\ref{eq:multiplier_coset_condition}) this relation must hold
for all odd integers $d$, and hence it is clear that the extension
$v_{\alpha}$ exists if and only if $n$ is a perfect square. It is
also easy to verify that in this case the multiplier system is given
by $v_{\alpha}\left(g_{1}\alpha g_{2}\right)=v\left(g_{1}\right)v\left(g_{2}\right)$.

Suppose now for simplicity that $n=p^{2}$, with $p\ne2$ a prime
number. The $p^{2}+p$ different coset representatives of $\overline{\Gamma}$
in $\overline{\Gamma}\alpha\overline{\Gamma}$ are given by $\alpha_{b}=\left(\begin{smallmatrix}1 & b\\
0 & p^{2}\end{smallmatrix}\right),\, b=0,\ldots,p^{2}-1,$ $\sigma=\left(\begin{smallmatrix}p^{2} & 0\\
0 & 1\end{smallmatrix}\right)$, and $\beta_{b}=\left(\begin{smallmatrix}p & b\\
0 & p\end{smallmatrix}\right),\, b=1,\ldots,p-1.$ By factorization of these representatives we see that $v_{\alpha}\left(\sigma\right)=v_{\alpha}\left(\alpha_{b}\right)=1$,
$b=0,\ldots,p^{2}-1$ and $v_{\alpha}\left(\beta_{b}\right)=\epsilon_{p}\left(\frac{b}{p}\right)$
for $b=1,\ldots,p-1$. Hence, for $p\ne2$ \begin{eqnarray}
T_{p^{2},\frac{1}{2}}^{v_{\theta}}f(z) & = & \frac{1}{p}\left\{ \sum_{b=0}^{p^{2}-1}\overline{v_{\alpha}}\left(\alpha_{b}\right)f\left(\alpha_{b}z\right)+\overline{v_{\alpha}}\left(\sigma\right)f\left(\sigma z\right)+\sum_{b=1}^{p-1}\overline{v_{\alpha}}\left(\beta_{b}\right)f\left(\beta_{b}z\right)\right\} \nonumber \\
 & = & \frac{1}{p}\left\{ \sum_{b=0}^{p^{2}-1}f\left(\frac{z+b}{p^{2}}\right)+f\left(p^{2}z\right)+\overline{\epsilon_{p}}\sum_{b=0}^{p-1}\left(\frac{b}{p}\right)f\left(z+\frac{b}{p}\right)\right\} .\label{eq:def_Heckeop_theta_mult}\end{eqnarray}
For $p=2$ the $4$ coset representatives are given by $\alpha_{b}=\left(\begin{smallmatrix}1 & b\\
0 & 4\end{smallmatrix}\right),\, b=0,\ldots,3$ and we obtain precisely the operator in (\ref{eq:Hecke_op_w_half_T4}).
The above construction is analogous to \cite[pp.\ 450--451, thm.\ 1.7]{shimura:73:half_integral}.
Suppose that $f\in\MAS\left(\Gamma_{0}(4),v_{\theta},\frac{1}{2},\lambda\right)$
has the following Fourier series at infinity \begin{equation}
f(z)=\sum_{n\ne0}\frac{a(n)}{\sqrt{|n|}}W_{\frac{1}{4}sgn(n),iR}(4\pi|n|y)e(nx)\label{eq:Fourier_exp_weight_half}\end{equation}
 and that $T_{p^{2},\frac{1}{2}}^{v_{\theta}}f(z)$ has a similar
Fourier expansion but with coefficients $b^{(p)}(n).$ Using the formula
for the standard Gauss sum, $\sum_{b=0}^{p^{2}-1}e\left(\frac{nb}{p^{2}}\right)=p^{2}$
if $p^{2}|n$ else $0$, and the twisted version $\sum_{b=1}^{p-1}\left(\frac{b}{p}\right)e\left(\frac{nb}{p}\right)=\sqrt{p}\left(\frac{n}{p}\right)\epsilon_{p}$
(cf.~\cite[Satz 7, p.\ 375]{borevich-shafarevich} or \cite[pp.\ 83-87]{MR1282723})
 we get that for all non-zero integers $n$\begin{eqnarray*}
b^{(p)}(n) & = & \left\{ a(p^{2}n)+a\left(\frac{n}{p^{2}}\right)+p^{-\frac{1}{2}}\left(\frac{n}{p}\right)a(n)\right\} ,\, p\ne2,\,\textrm{and}\\
b^{(2)}(n) & = & a\left(4n\right).\end{eqnarray*}
(We use the standard convention that $a\left(x\right)=0$ if $x\not\in\mathbb{Z}$.)
It is now obvious that our Hecke operator $T_{p^{2},\frac{1}{2}}^{v_{\theta}}$
is equal to the corresponding Hecke operator defined in \cite[p.\ 199]{katok-sarnak}.
Observe that our Fourier coefficients $a(n)=$$\sqrt{n}\times$Katok-Sarnak's
Fourier coefficients $b(n).$

As usual, we consider Hecke eigenforms in $\MAS\left(\Gamma_{0}(4),v_{\theta},\frac{1}{2},\lambda\right)$
which are simultaneous eigenfunctions of all $T_{p^{2},\frac{1}{2}}^{v_{\theta}}$
with $p\ne2$. An additional commuting normal operator can be chosen
as either $T_{4,\frac{1}{2}}^{v_{\theta}}$ \emph{or} $L$ (these
two operators does not commute in general). 

As it turns out, the operator $L$ is particularly useful in connection
with the Shimura correspondence on the Kohnen space, $V^{+}$, where
$L$ has eigenvalue $1$ (cf.~Section \ref{sub:The-Shimura-correspondence}).

Observe that the Hecke eigenvalues in this case does not equal to
the Fourier coefficients. In fact, suppose that $f$ as in (\ref{eq:Fourier_exp_weight_half})
is an eigenfunction of all Hecke operators with $T_{p^{2},\frac{1}{2}}^{v_{\theta}}f=\lambda_{p}f$
and that $a(t)\ne0$ for some integer $t$. It is then easy to see
that \begin{eqnarray}
\lambda_{p} & = & \left\{ \frac{a\left(tp^{2}\right)}{a(t)}+\frac{\left(\frac{t}{p}\right)}{\sqrt{p}}\right\} ,\, p\ne2,\,\textrm{and}\label{eq:Hecke_eigenvalue_p_weight_1_2}\\
\lambda_{2} & = & \frac{a(4t)}{a(t)}.\nonumber \end{eqnarray}
Using multiplicative relations of the Hecke operators one can prove
(cf.~\cite[p.\ 453]{shimura:73:half_integral}) that if $t$ is square
free, then \[
a(tm^{2})a(tn^{2})=a(t)a(tm^{2}n^{2}),\,\textrm{for}\,(m,n)=1.\]
Furthermore, if $f$ is also an eigenfunction of $T_{4,\frac{1}{2}}^{v_{\theta}}$
then \[
a(m)a\left(4n\right)=a(4m)a(n),\, m,n\in\mathbb{Z}.\]

\subsubsection{Hecke operators at integer weights and Fourier coefficients\label{sub:Lifts-at-weight}}

Consider the modular group, $\Gamma=\PSLZ,$ together with the integer
weight $k\not\equiv0\mod12$ and multiplier $v=v_{\eta}^{2k}$. 

It can be shown that for all positive integers $n$ and $m$ with
$kn\equiv-km\equiv k\mod12$ we can construct Hecke operators $\T{n}$
and Hecke-like operators $\Th{m}$ acting on $\Mas{\Gamma,v,k,\lambda}$.
Using these operators one can obtain multiplicative relations for
Fourier coefficients similar to the weight zero case. This is shown
in detail in \cite{hecke_operators_weight} and here we will only
state the theorem and give a brief outline of the ideas of the proof. 

\begin{thm}
\label{thm:hecke_operators_weight}Let $k$ be an integer $k\not\equiv0\mod12$
and $R>0$ then there exist a basis of $\Mas{\Gamma,v,k,R}$ consisting
of Maass wave forms $f$ with Fourier expansions at infinity \[
f(z)=\sum_{n=-\infty}^{\infty}\frac{c\left(n\right)}{\sqrt{\left|n+\frac{k}{12}\right|}}W_{\frac{k}{2}\sgn\left(n+\frac{k}{12}\right),iR}\left(4\pi\left|n+\frac{k}{12}\right|y\right)e\left(\left(n+\frac{k}{12}\right)x\right)\]
where the coefficients $c\left(n\right)$ satisfies the following
multiplicative relations if $c\left(0\right)\ne0$. For positive integers
$m,n$ with $12m\equiv12n\equiv0\mod k$ set $m_{1}=\frac{12m}{k}$,
$n_{1}=\frac{12n}{k}$ and $D=\frac{k}{\left(12,k\right)}$. If $\left(m_{1}+1,D\right)=\left(n_{1}+1,D\right)=1$
\begin{equation}
c(m)c(n)=c(0)\sum_{0<d|(m_{1}+1,n_{1}+1)}\chi_{k}\left(d\right)\, c\left(\frac{k}{12}\left(\frac{\left(m_{1}+1\right)\left(n_{1}+1\right)}{d^{2}}-1\right)\right),\label{eq:w1_coeff_basic_mult_rel}\end{equation}
and if $\left(m_{1}-1,D\right)=\left(n_{1}-1,D\right)=1$ then \begin{equation}
c\left(-m\right)c\left(-n\right)=\Lambda_{k,R}\; c\left(0\right)\sum_{0<d|\left(m_{1}-1,n_{1}-1\right)}\chi_{k}\left(d\right)\, c\left(\frac{k}{12}\left(\frac{\left(m_{1}-1\right)\left(n_{1}-1\right)}{d^{2}}-1\right)\right),\label{eq:mult_rel_w1_neg_coef}\end{equation}
where $\chi_{k}\left(d\right)=i^{k\left(d-1\right)}$ ($=\left(\frac{-1}{d}\right)^{k}$
for odd $d$) and \begin{equation}
\Lambda_{k,R}=\begin{cases}
\prod_{j=1}^{l}\left(j\left(j-1\right)+\frac{1}{4}+R^{2}\right)^{2}, & k=2l,\\
-R^{2}\prod_{j=1}^{l}\left(j^{2}+R^{2}\right)^{2}, & k=2l+1.\end{cases}\label{eq:Lambda_kR}\end{equation}

\end{thm}
In particular, we see that if $k|12$ we have $D=1$ and the multiplicative
relations (\ref{eq:w1_coeff_basic_mult_rel}) and (\ref{eq:mult_rel_w1_neg_coef})
are valid for all positive integers. 

\begin{rem}
As in the weight zero case and the coefficient $c\left(1\right)$
one can show that if $D=1$ and $f$ is an eigenfunction of all Hecke
operators (defined below) then $c\left(0\right)\ne0$ unless $f\left(z\right)$
is identically $0$. In the case $D>0$, if $c\left(0\right)=0$ and
$f\left(z\right)\not\equiv0$ we can choose an integer $n_{0}$ such
that $c\left(n_{0}\right)\ne0$ and obtain a similar set of multiplicative
relations. 
\end{rem}
The proof of the {}``positive part'' of the theorem, i.e. equation
(\ref{eq:w1_coeff_basic_mult_rel}) relies on the construction of
a family of Hecke operators $\T{m}$ with $km\equiv k\mod12$. It
is shown that this family consist of self-adjoint operators commuting
with each other and the weight $k$ Laplacian. An explicit form of
$\T{m}$ is \begin{equation}
\T{m}f\left(z\right)=\frac{1}{\sqrt{m}}\sum_{ad=m,d>0}\chi_{k}\left(d\right)\sum_{b=0}^{d-1}\overline{v\left(T\right)^{bd}}f\left(\frac{az+b}{d}\right).\label{eq:Hecke_weight_k_Tm_formula}\end{equation}
It is not hard to show that if $f\left(z\right)$ has Fourier coefficients
$c\left(n\right)$ then $\T{m}$ has coefficients \[
b\left(n\right)=\sum_{0<d|\left(m,n-\frac{k\left(m-1\right)}{12}\right)}\chi_{k}\left(d\right)c\left(\frac{nm}{d^{2}}+\frac{k\left(m-d^{2}\right)}{12d^{2}}\right)\]
from which we see that if $\T{m}f=\lambda_{m}f$ and $c\left(0\right)\ne0$
then\[
\lambda_{m}=\frac{1}{c\left(0\right)}\left[c\left(\frac{k\left(m-1\right)}{12}\right)+\chi_{k}\left(D\right)c\left(\frac{k\left(m-D^{2}\right)}{12D^{2}}\right)\right].\]
The proof of (\ref{eq:w1_coeff_basic_mult_rel}) is concluded with
a proof of the following multiplication rule (using straight-forward
calculations and induction). If $km\equiv kn\equiv k\mod12$ then
\[
\T{m}\T{n}=\sum_{0<d|\left(m,n\right)}\chi_{k}\left(d\right)\T{\frac{mn}{d^{2}}}.\]

To obtain the {}``negative part'' of theorem, i.e. (\ref{eq:mult_rel_w1_neg_coef})
we have to consider another family of operators, $\Th{m}$ with $km\equiv-k\mod12$.
These operators are given as a combination of a bijection $\Theta=J\circ\mathcal{E}_{k}^{-}$
mapping $\Mas{\Gamma,v,k,\lambda}$ to $\Mas{\Gamma,\overline{v},k,\lambda}$
and an Hecke operator $\T[\overline{v}]{m}$ mapping $\Mas{\Gamma,\overline{v},k,\lambda}$
back to $\Mas{\Gamma,v,k,\lambda}$. Here $\mathcal{E}_{k}^{-}=E_{2-k}^{-}\circ\cdots\circ E_{k-2}^{-}\circ E_{k}^{-}$
maps $\Mas{\Gamma,v,k,\lambda}$ to $\Mas{\Gamma,v,-k,\lambda}$ bijectively
since $\lambda>\frac{1}{4}$ and $J$ reflects this space back to
$\Mas{\Gamma,\overline{v},k,\lambda}$. $\T[\overline{v}]{m}$ is
given by (\ref{eq:Hecke_weight_k_Tm_formula}) with $v$ interchanged
with $\overline{v}$. The operator $\Theta$ is similar to the operator
defined by Maass in \cite[p.\ 181]{maass:modular_functions}. Using
(\ref{eq:Maass_operators_act_1})-(\ref{eq:Maass_operators_act_4})
it is easy to show that $\Theta^{2}f=\Lambda_{k,R}f$ for all $f\in\Mas{\Gamma,k,v,\lambda}$
and that if $f\left(z\right)$ has Fourier coefficients $c\left(n\right)$
then $\Th{m}f\left(z\right)$ has coefficients \[
d\left(n\right)=\sum_{0<d|\left(m,n-\frac{k\left(m+1\right)}{12}\right)}\chi_{k}\left(d\right)c'\left(\frac{nm}{d^{2}}+\frac{k\left(m+d^{2}\right)}{12d^{2}}\right)\]
where $c'\left(n\right)=c\left(-n\right)\begin{cases}
1, & n\ge1,\\
\Lambda_{k,R}, & n\le0.\end{cases}$ If $\Th{m}f=\mu_{m}f$ and $c\left(0\right)\ne0$ then \[
\mu_{m}=\frac{1}{c\left(0\right)}\left[c\left(\frac{k\left(m+1\right)}{12}\right)+\chi_{k}\left(D\right)c\left(\frac{k\left(m+D^{2}\right)}{12D^{2}}\right)\right].\]
The multiplication rule for the operators $\Th{m}$ is that for $km\equiv kn\equiv k\mod12$
we have \[
\Th{m}\Th{n}=\Lambda_{k,R}\sum_{0<d|\left(m,n\right)}\chi_{k}\left(d\right)\Th{\frac{mn}{d^{2}}}.\]
This in combination with the expressions for $\mu_{m}$ and $\lambda_{m}$
concludes the proof of (\ref{eq:mult_rel_w1_neg_coef}).

\begin{example}
Look at the specific case $k=1$ and a function $f\in\Mas{\Gamma,v,1,\lambda}$
then by by setting $m=1$ in (\ref{eq:w1_coeff_basic_mult_rel}) and
using the normalization $c(1)=1$ we see that \[
c\left(n\right)=\sum_{0<d|\left(12n+1,13\right)}\chi_{1}\left(d\right)c\left(\frac{\left(13\left(12n+1\right)-d^{2}\right)}{12d^{2}}\phantom{}\right)\]
 and hence if $(12n+1,13)=1$ then we get a striking proportionality
relation:

\begin{align}
c\left(n\right)= & c(0)c\left(\frac{13(12n+1)-1}{12}\right)=c(0)c\left(13n+1\right).\label{eq:w1_coeff_prop1}\\
\prefixtext{\mbox{For}\,\ensuremath{12n+1=13l}\,\mbox{we get}}\nonumber \\
c\left(n\right)= & c(0)\left(c\left(13n+1\right)+c\left(\frac{n-1}{13}\right)\right),\label{eq:w1_coeff_prop2}\\
\prefixtext{\mbox{and if\,}\ensuremath{\left(l,13\right)=1}\mbox{\, we can combine these two equations}\mbox{ and obtain:}}\nonumber \\
c\left(n\right)= & c(0)\left(c\left(13n+1\right)+c(0)c\left(n\right)\right)\label{eq:w1_coeff_prop3}\end{align}
and hence \begin{equation}
c(n)=\frac{c(0)}{1-c(0)^{2}}c\left(13n+1\right).\label{eq:wt_1_hecke_rel_pos}\end{equation}
Consider now $m=-1$ in (\ref{eq:mult_rel_w1_neg_coef}) and note
that $\Lambda_{1,R}=-R^{2}$. If $(12n-1,11)=1$ then \begin{equation}
c(-n)=\frac{-R^{2}c(0)}{c(-1)}c(11n-1),\label{eq:wt_1_hecke_rel__neg_1}\end{equation}
and if $(12n-1,11)=11$ then \begin{equation}
c(-n)=\frac{-R^{2}c(0)}{c(-1)}\left[c(11n-1)-c\left(\frac{n-1}{11}\right)\right].\label{eq:wt_1_hecke_rel__neg_2}\end{equation}

\end{example}
We have not seen relations of type (\ref{eq:w1_coeff_prop1})-(\ref{eq:wt_1_hecke_rel__neg_2})
earlier in the literature. For some numerical examples see Tables
\ref{tab:N1R47709} and \ref{tab:N1R36624}. Examples of relations
for higher weights can be found in \cite{hecke_operators_weight}. 

\begin{rem}
An alternative approach to the above coefficient relations in the
case of weight $1$ is to identify $\Mas{\Gamma_{0}(1),v_{\eta}^{2},1,\lambda}$
with a subspace of $\Mas{\Gamma_{0}(144),\left(\frac{-1}{d}\right),1,\lambda}$
via the map $f(z)\mapsto g(z)=f(12z)$. Cf. e.g.~\cite[sec.\ 2.4.7]{stromberg:thesis}. 

This identification also provides an explanation for the occurrence
of CM-type forms found (numerically) in $\Mas{\Gamma_{0}(1),v_{\eta}^{2},1,\lambda}$.
These forms have eigenvalues in an arithmetic progression: $R_{k}=\frac{2\pi k}{\ln\left(7+2\sqrt{12}\right)},$
for $k\in\mathbb{Z}.$ See Table \ref{tab:egenvN1_w1}. 
\end{rem}

\section{The Eisenstein series for $\PSLZ$ at weight zero}

In case there is a cusp $p_{j}$ at which the multiplier system is
singular (i.e. $v(\Tm_{j})=1$) we have a continuous spectrum: $[\frac{1}{4},\infty)$
(with multiplicity equal to the number of singular cusps), and in
general we can not say much about the embedded discrete spectrum in
$[\frac{1}{4},\infty)$. 

Examples of singular cusps are the cusp at infinity for the eta multiplier
and weight $k\equiv0\mod12$ on $\PSLZ$ and the cusps at $0$ and
$i\infty$ for the theta multiplier and weight $\frac{1}{2}$ on $\Gamma_{0}\left(4\right)$. 

Remember that Maass waveforms are part of the discrete spectrum, but
as we continuously {}``turn off'' the multiplier system, i.e. for
$v=v_{\eta}^{2k}$ we let $k\rightarrow0$, the continuous spectrum
will emerge in the limit. For this reason we want to review some details
concerning the Eisenstein series on the modular group.

At weight $0$ and singular character $\chi$, the continuous spectrum
of $\Delta=\Delta_{0}$ is the interval $[\frac{1}{4},\infty)$ and
the eigenpacket is given by the Eisenstein series $E(z;s;\chi)$ defined
by

\[
E(z;s;\chi)=\sum_{T\in\Gamma_{\infty}\backslash\Gamma}\chi(T^{-1})\left(\Im(Tz)\right)^{s},\]
where $\Gamma_{\infty}=[S].$ For the trivial character we have the
Fourier series expansion (cf.~\cite[p.\ 65 and p.\ 76]{hejhal:lnm1001})\[
E(z;s;\chi)=y^{s}+\varphi(s)y^{1-s}\sum_{n\neq0}\varphi_{n}(s)\sqrt{y}K_{s-\frac{1}{2}}(2\pi|n|y)e(nx),\]
 where \begin{eqnarray*}
\varphi(s) & = & \sqrt{\pi}\frac{\Gamma\left(s-\frac{1}{2}\right)}{\Gamma(s)}\frac{\zeta(2s-1)}{\zeta(2s)},\,\textrm{and}\\
\varphi_{n}(s) & = & \frac{2\pi^{s}|n|^{s-\frac{1}{2}}}{\Gamma(s)}\frac{\sigma_{1-2s}(|n|)}{\zeta(2s)}.\end{eqnarray*}
 Hence we can see that for $s=\frac{1}{2}+iR$ , the $n$th Fourier
coefficient of $E(z;s)$ is given by \begin{eqnarray}
c(n) & = & \varphi_{n}\left(\frac{1}{2}+iR\right)=K\cdot|n|^{iR}\sigma_{-2iR}(|n|)\label{eq:eisenstein_coeff}\\
 & = & K\cdot|n|^{iR}\sum_{d||n|,d>0}d^{-2iR},\nonumber \end{eqnarray}
where $K=K(R)$ is a constant dependent on $R$. For a prime $p$
we get \begin{eqnarray*}
c(p) & = & K\cdot p^{iR}(1+p^{-2iR})=2K\cdot\cos(R\ln p).\end{eqnarray*}
Using this formula we can compute quotients of various $c(p)$ (e.g.~$\frac{c(2)}{c(3)}$)
and compare this with corresponding quotients for the experimentally
obtained forms of weight $k\approx0$.

\section{\label{sub:The-Shimura-correspondence}The Shimura correspondence}

We know that the $\theta-$function is an automorphic form (not a
cusp form) of weight $\frac{1}{2}$ on $\Gamma_{0}(4),$ hence we
consider $\Gamma_{0}(4)$ together with the $\theta-$multiplier system
(cf. section \ref{sub:The-theta-multiplier}).

\subsection{Introduction -- the holomorphic case}

We will consider the Shimura correspondence only in the particular
case of trivial character and level a square-free multiple of $4$.
Let $S_{k}(N)$ denote the space of holomorphic cusp forms of weight
$k\in\mathbb{Z}$ (and trivial multiplier) on $\Gamma_{0}(N),$ and
let $S_{k+\frac{1}{2}}(4N)$ denote the space of holomorphic cusp
forms of weight $k+\frac{1}{2},\, k\in\mathbb{Z}$ and multiplier
$v_{\theta}$, on $\Gamma_{0}(4N)$. The \emph{Shimura correspondence}
is a correspondence between the space $S_{k+\frac{1}{2}}(4N)$ and
spaces $S_{2k}(N')$ for certain integers $N'|4N$ (e.g. $N'=2N$
or $N$). 

The map from $S_{k+\frac{1}{2}}$ to $S_{2k}$ was first constructed
by Shimura \cite{shimura:73:half_integral} and later an adjoint map
from $S_{2k}$ to $S_{k+\frac{1}{2}}$ was constructed by Shintani
\cite{shintani:75:onconstruction}. The former uses a Dirichlet-series
and the latter uses an integral against a theta-function. Both these
maps commute with the Hecke operators that are acting on $S_{2k}(N)$
and $S_{k+\frac{1}{2}}(4N)$ respectively. Kohnen \cite{MR81j:10030,MR84b:10038}
proved that for $N$ odd and square-free, the correspondence is a
bijection between the newforms on $S_{2k}(N)$ and a certain subspace,
$V^{+}\subseteq S_{k+\frac{1}{2}}(4N)$. The subspace $V^{+}$ is
composed of Hecke eigenfunctions whose Fourier coefficients, $c(n)$,
satisfy certain vanishing properties; namely, $c(n)=0$ for $n\equiv2,3\mod4$
(see also Section \ref{sub:The-operator-L}). 

Following from the Shimura correspondence is also a connection between
certain Fourier coefficients of the half integral weight forms and
critical values of twisted L-series for the corresponding integral
weight form. Cf.~e.g.~\cite{MR577010,MR646366}, \cite{MR629468,MR783554}
and \cite{shimura:fc_of_hilbert}.

\subsection{The Shimura correspondence for Maass waveforms \label{sub:Extension-to-the-non-holo}}

The extension of the Shimura correspondence and Kohnen's result to
spaces of Maass waveforms has been investigated by e.g. Sarnak \cite{sarnak:82:MR750670},
Hejhal \cite{MR726196}, Duke \cite{duke:88}, Katok-Sarnak \cite{katok-sarnak},
Khuri-Makdisi \cite{khuri:fc_of_hilbert}, Kojima \cite{kojima:95,kojima:00,kojima:04},
Biró \cite{biro:00} and Arakawa \cite{arakawa:2}. Of these, \cite{khuri:fc_of_hilbert}
and \cite{kojima:00,kojima:04} are written in the more general setting
of Hilbert modular forms for number fields. Reading \cite{khuri:fc_of_hilbert}
together with \cite{shimura:fc_of_hilbert} and \cite{katok-sarnak}
gives a good picture of the current state of affairs. 

Throughout this section let $\MAS\left(N,R\right)^{+}$ denote the
space of even (with respect to $J:z\mapsto-\overline{z}$) Hecke normalized
Maass waveforms in $\MAS\left(\Gamma_{0}(N),1,0,R\right)$ and let
$\MAS_{\frac{1}{2}}\left(4,R\right)$ denote the space of Hecke normalized
(with respect to all $T_{p^{2}},$ $p\ne2$) weight $\frac{1}{2}$
Maass waveforms in $\MAS\left(\Gamma_{0}(4),v_{\theta},\frac{1}{2},R\right)$.
Also let $V^{+}\subseteq\MAS_{\frac{1}{2}}\left(4,R\right)$ be defined
as in Section \ref{sub:The-operator-L}. 

For $f\in\MAS\left(N,R\right)^{+}$ and $\phi\in\MAS_{\frac{1}{2}}\left(4,R\right)$
we will use the following notation: \begin{align*}
f(z) & =\sum_{\underset{n\ne0}{n=-\infty}}^{\infty}A(n)\sqrt{y}K_{iR}\left(2\pi|n|y\right)e(nx),\\
\prefixtext{and}\\
\phi\left(z\right) & =\sum_{\underset{n\ne0}{n=-\infty}}^{\infty}\frac{a(n)}{\sqrt{|n|}}W_{\frac{1}{4}sgn(n),iR}\left(4\pi|n|y\right)e\left(nx\right).\end{align*}
 (Both expansions are given with respect to the cusp at $\infty$).

The existence of a Shimura correspondence and an inverse for Maass
waveforms is expressed by the following proposition (essentially \cite[thm. 5.1 and 5.2]{khuri:fc_of_hilbert}): 

\begin{prop}
\label{thm:existence_of_shim_corr}~
\begin{lyxlist}{00.00.0000}
\item [{a)}] The Shimura correspondence gives a map $\Phi:\MAS_{\frac{1}{2}}\left(4,R\right)\rightarrow\MAS\left(2,2R\right)^{+}$.
\item [{b)}] Conversely, let $f\in\MAS\left(2,2R\right)^{+}$. Then there
exists a $\phi\in\Phi^{-1}(f)\in\MAS_{\frac{1}{2}}\left(4,R\right)$. 
\item [{c)}] If $A(p)$ is the Hecke eigenvalue of $f$ with respect to
the operator $T_{p}$ and $\lambda_{p}$ is the eigenvalue of $\phi\in\Phi^{-1}(f)$
with respect to $T_{p^{2},\frac{1}{2}}^{v_{\theta}}$ then we actually
have \[
A(p)=\lambda_{p}.\]

\end{lyxlist}
\end{prop}
\begin{proof}
See the proof of \cite[thm. 5.1 and 5.2]{khuri:fc_of_hilbert}. For
c) observe the difference in normalization of the Hecke operators.
Cf. also~\cite[thm.\ 2]{kojima:95}, Proposition \ref{thm:Shimura_corr_prop}
a) and (\ref{eq:Hecke_eigenvalue_p_weight_1_2}).
\end{proof}
\begin{rem}
From \cite[prop. 4.1 and 2.3]{katok-sarnak} we know that the above
correspondence $\Phi$ actually restricts to a map between $\MAS\left(1,2R\right)^{+}$
and the subspace $V^{+}$. And hence a new form in $\MAS\left(2,2R\right)^{+}$
will be mapped to $V^{-}$, the orthogonal complement of $V^{+}$.
Vice versa, a form in $V^{-}$ will be mapped to a new form in $\MAS\left(2,2R\right)^{+}$. 
\end{rem}
To generalize the results mentioned at the end of the previous subsection
to Maass waveforms we need the following definition. For $f\in\MAS(2,R)^{+}$
with Fourier coefficients $\left\{ A(n)\right\} $ and a given Dirichlet
character $\chi$ we define the $\chi$-twisted L-series of $f$ by
\[
L(f,\chi,s)=\sum_{n\ne0}^{\infty}A(n)\chi(n)|n|^{-s-\frac{1}{2}}.\]

\begin{prop}
\label{thm:Shimura_corr_prop}Let $\phi\in\MAS_{\frac{1}{2}}\left(4,R\right)$
have Fourier coefficients $\left\{ a(n)\right\} $ and correspond
(via Prop. \ref{thm:existence_of_shim_corr}) to $f=\Phi(\phi)\in\MAS\left(N,2R\right)^{+}$
(where $N=1$ if and only if $\phi\in V^{+},$ otherwise $N=2$) with
Fourier coefficients $\left\{ A(n)\right\} $. Let $t\in\mathbb{Z}^{+}$
be square free and let $\chi_{t}'$ be the quadratic residue symbol
$\left(\frac{t}{\cdot}\right)$ considered mod $Nt$ (i.e. we have
$\chi'_{t}(n)=\left(\frac{t}{n}\right)$ if $\left(n,N\right)=1$
and $\chi'_{t}(n)=0$ otherwise. 

Then the following properties hold: 
\begin{lyxlist}{00.00.0000}
\item [{a)}] If $a(t)\ne0$ then $A(n)$ can be expressed by: 
\end{lyxlist}
\begin{eqnarray}
A(n) & = & \sum_{\underset{k>0}{kd=n}}\frac{\chi_{t}'\left(k\right)}{\sqrt{k}}\frac{a\left(td^{2}\right)}{a\left(t\right)},\, n\in\mathbb{Z}^{+}.\label{eq:shimura_corr_coeff}\end{eqnarray}

To express $a(t)$ in terms of $A(n)$'s we get three cases.

\begin{lyxlist}{00.00.0000}
\item [{b)}] If $f$ is an oldform. i.e. $f\in\MAS(1,2R)^{+}$ then $\phi\in V^{+}$
and hence \[
a(t)=0,\, t\equiv2,3\mod4.\]

\item [{c)}] If $f$ is a newform with eigenvalue $\epsilon$ with respect
to the involution $z\mapsto\frac{-1}{2z}$ then \[
a(t)=0\]
 for $t\equiv5\mod8$ if $\epsilon=1$ and for $t\equiv1\mod8$ if
$\epsilon=-1$. 
\item [{d)}] For all other (square free) values of $t$ the following formula
holds\[
\left|a(t)\right|^{2}=Q\frac{\left\langle \phi,\phi\right\rangle }{\left\langle f,f\right\rangle }L\left(f,\chi_{t},0\right)\]
where $Q$ is a constant independent of $t$. 
\item [{e)}] If $\phi$ is a normalized Hecke eigenform (for all $T_{p^{2},\frac{1}{2}}^{v_{\theta}}$,
$p$ an odd prime) then $f$ is also a normalized Hecke eigenform. 
\end{lyxlist}
\end{prop}
\begin{proof}
Cf.\emph{~\cite[pp.\ 448(prop.\ 1.3), 458 (main theorem), 474 (line -11)]{shimura:73:half_integral},}
and \cite[p.\ 422]{khuri:fc_of_hilbert} for the choice of $\chi_{t}'$
(the $N=1$ case is simply to incorporate the results of \cite{katok-sarnak}).
Relation a) follows from \cite[thm.\ 5.1 (2)]{khuri:fc_of_hilbert},
which in our case can be written as \begin{equation}
\sum_{n=1}^{\infty}A(n)n^{-s}=c\cdot L\left(s+\frac{1}{2},\chi_{t}'\right)\sum_{n=1}^{\infty}a\left(tn^{2}\right)n^{-s}.\label{eq:Lfun_id_shimura}\end{equation}
 (Cf.~e.g.~also \cite[p.\ 458]{shimura:73:half_integral}, \cite[p.\ 159]{niwa:half_integral},
\cite[prop. 3.1]{shimura:fc_of_hilbert}.) The relation b) follows
from \cite[prop.\ 2.3]{katok-sarnak} (see also \cite{MR81j:10030})
and relation c) is implicit in \cite[thm.\ 8.1 (8.6)]{khuri:fc_of_hilbert}.
Simply observe that in our normalization the relevant factor in the
product is given by \[
\left(\sqrt{2}A(2)-\left(\frac{2}{t}\right)\right)=-\left(\epsilon+\left(\frac{2}{t}\right)\right)\]
 since $A(2)=\frac{-\epsilon}{\sqrt{2}}$ for a newform (analogous
to \cite[thm.\ 3]{atkinlehner}). (Alternatively consider the sign
of the functional equation for $L\left(f,\chi_{t},s\right)$.) Relation
d) is the Maass waveform-analogue of \cite[thm.\ 1]{MR629468} and
follows from \cite[thm.\ 8.1]{khuri:fc_of_hilbert} (in the notation
of \cite{khuri:fc_of_hilbert} we have $a(n)=|n|^{\frac{1}{4}}\mu(n;\phi)$).
\end{proof}
\begin{rem}
That the map $\Phi$ is well-defined and surjective from $\MAS_{\frac{1}{2}}\left(4,R\right)$
onto $\MAS\left(2,2R\right)^{+}$ follows from \cite[thm.\ 5.1 and 5.2]{khuri:fc_of_hilbert}
and {}``multiplicity one'' for the Hecke operators $T_{p}$ on $\MAS\left(2,2R\right)^{+}$
(cf.~\cite[thm.\ 4.6]{Andreas:thesis}). Experimentally we have observed
that $\Phi$ restricted to $V^{+}$ also seems to be injective. Theoretically,
this is still an open problem that might be possible to resolve using
the trace formula for the Hecke operators $T_{p^{2}}$ on $V^{+}$.

\end{rem}

\begin{rem}
Note that the same argument as for \ref{thm:Shimura_corr_prop} c)
implies that to any $f\in\MAS\left(2,2R\right)^{+}$ with Fourier
coefficient $A(2)=\frac{\pm1}{\sqrt{2}}$ there corresponds a function
$\phi\in\Phi^{-1}\left(f\right)\subseteq\MAS_{\frac{1}{2}}\left(4,R\right)$
with coefficients $a(n)=0$ for either $n\equiv1$ or $5$ $\mod8$
respectively. Hence, since oldforms on $\Gamma_{0}(2)$ occur in pairs,
we can choose two forms $f_{1},f_{2}$ with $A_{1}(2)=\frac{1}{\sqrt{2}}$
respectively $A_{2}(2)=\frac{-1}{\sqrt{2}}$. These two functions
thus correspond to two linearly independent non-zero functions in
$\MAS_{\frac{1}{2}}\left(4,R\right)$ which hence is at least two
dimensional when $2R$ is an even eigenvalue for $\PSLZ$. 
\end{rem}
The relations a)--c) in Proposition \ref{thm:Shimura_corr_prop} appeared
in \cite[ch.\ 2]{stromberg:thesis} as experimental observations (cf.~\ref{sub:Half-integer-weight},
in particular Tables \ref{tab:N4R6889}-\ref{tab:fourier_coeff_comparison}).

See also \cite[\S4.1]{MR993311} and \cite[p.\ 633]{MR904946}, and
\cite[pp.\ 502 (bottom) -- 504 (top)]{shimura:fc_of_hilbert} for
some additional perspectives.

\section{Some computational remarks\label{sec:Some-Computational-Remarks}}

We recall the key ingredients in the standard Hejhal's algorithm (cf.~e.g.~\cite{hejhal:92,hejhal:99_eigenf,hejhal:calc_of_maass_cusp_forms},
\cite{stromberg:04:1,stromberg:thesis} and \cite{andreas:effective_comp})
to compute weight zero Maass waveforms on cofinite Fuchsian groups.
First of all the asymptotic properties of the K-Bessel function are
used to obtain an approximation the Maass waveform by a truncated
Fourier series. By viewing this as a discrete Fourier transform one
can use the inverse transform to express any coefficient as a linear
combination of the other coefficients. Finally, by using the assumed
automorphy properties of the function one obtains a non-trivial linear
system that is satisfied by the Fourier coefficients. 

To use this algorithm to also \emph{locate} eigenvalues the most general
method is to use two different sets of sampling points for the inverse
transform and try to minimize the difference between the correspondingly
computed solution vectors. 

The following two modifications are needed in order to make the weight
zero algorithm work in the general case: 

\begin{itemize}
\item the $K$-Bessel function needs to be replaced with the Whittaker function,
and
\item the automorphy condition needs to incorporate the multiplier system.
\end{itemize}
The first modification, although trivial in theory is the most complex
in terms of the numerics involved. There was no efficient algorithm
for the Whittaker function in the literature and a new algorithm had
to be developed. We used the integral representation (cf.~\cite[p.\ 375 (top)]{hejhal:lnm1001})\[
W_{k,iR}(2x)=\sqrt{\frac{2x}{\pi}}\int_{0}^{\infty}e^{-x\cosh t}\Psi\left[-k;\frac{1}{2};x(1+\cosh t)\right]\cosh\left[iRt\right]dt,\]
where $\Psi$ is a confluent hypergeometric function together with
a stationary phase method to develop a robust and efficient algorithm.
This algorithm is in essence similar to Hejhal's algorithm (cf.~\cite{hejhal:92})
for the K-Bessel function, $K_{iR}(x)$, and the generalization of
it made by Avelin (cf.~\cite[Avelin]{helen:deform_published}) to
a complex argument, $K_{s}(x),$ $s\in\mathbb{C}.$ The details of
this algorithm is described in \cite[Ch.\ 4]{stromberg:thesis}.

Let $\Gamma$ be a Fuchsian group with fundamental domain $\mathcal{F}$
and let $p_{j},\,\sigma_{j},\,1\le j\le\kappa$, $q_{i}$ and $U_{i},\,1\le i\le\kappa_{0}$
be as in Section \ref{sub:Summary-of-Notation}. Let $I\left(w\right)=i$
if $q_{i}$ is the closest (in the hyperbolic metric) parabolic vertex
to $w\in\mathcal{F}$. 

Consider $z\in\H$ and let $T_{j}\in\Gamma$ be the pull-back map
of $\sigma_{j}z$, i.e. $T_{j}\sigma_{j}z\in\mathcal{F}$. Observe
that $f_{j}(z)=f_{|\sigma_{j}}(z)=j_{\sigma_{j}}(z;k)^{-1}f(\sigma_{j}z),$
and with $z_{j}^{*}=\sigma_{I(j)}^{-1}U_{I(j)}T_{j}\sigma_{j}z$ where
we set $I\left(j\right)=I\left(T_{j}\sigma_{j}z\right)$ (cf.~\cite[p.\ 23]{stromberg:04:1})
the automorphy condition now becomes:\begin{eqnarray}
f_{j}(z) & = & j_{\sigma_{j}}(z;k)^{-1}f(\sigma_{j}z)=j_{\sigma_{j}}(z;k)^{-1}f(T_{j}^{-1}U_{I(j)}^{-1}\sigma_{I(j)}z_{j}^{*})\nonumber \\
 & = & j_{\sigma_{j}}(z;k)^{-1}v(T_{j}^{-1}U_{I(j)}^{-1},k)j_{T_{j}^{-1}U_{I(j)}^{-1}}(\sigma_{I(j)}z_{j}^{*};k)f(\sigma_{I(j)}z_{j}^{*})\label{eq:autmorphy_cond_weight}\\
 & = & j_{\sigma_{j}}(z;k)^{-1}j_{T_{j}^{-1}U_{I(j)}^{-1}}(\sigma_{I(j)}z_{j}^{*};k)j_{\sigma_{I(j)}}(z_{j}^{*};k)v(T_{j}^{-1}U_{I(j)}^{-1},k)f_{I(j)}(z_{j}^{*}).\nonumber \end{eqnarray}
The entire setup of the Phase 1 algorithm, i.e. locating eigenvalues
and computing the first few Fourier coefficients, goes through exactly
as in the case of weight $0$ (cf.~\cite{stromberg:04:1} or \cite[Ch.\ 1]{stromberg:thesis})
with some trivial modifications. For the sake of completeness we will
give an outline of the modified algorithm. Consider $f\in\Mas{\Gamma,v,k,\lambda}$
and using the notation $n_{i}=n+\alpha_{i}$ $f$ has a Fourier series
expansion at the cusp $p_{i}$: \[
f_{i}(z)=f_{|\left[k,\sigma_{j}\right]}\left(z\right)=\sum_{\underset{{n+\alpha_{i}\neq0}}{-\infty}}^{\infty}\frac{c_{i}(n)}{\sqrt{|n_{i}|}}W_{\frac{k}{2}sgn(n_{i}),iR}(4\pi|n_{i}|y)e\left(n_{i}x\right),\]
and since the Whittaker function is ultimately exponentially decaying,
given an $\epsilon>0,$ there exists a constant (depending on $y$
and $\epsilon$), $M(y),$ such that \[
f_{i}(z)=\hat{f}_{i}(z)+[[\epsilon]],\]
where we use $[[\epsilon]]$ to denote a constant of magnitude less
than $\epsilon.$ The truncated Fourier series \[
\hat{f}_{i}(z)=\sum_{\underset{{n+\alpha_{i}\neq0}}{-M(Y)}}^{M(Y)}\frac{c_{i}(n)}{\sqrt{|n_{i}|}}W_{\frac{k}{2}sgn(n_{i}),iR}(4\pi|n_{i}|y)e\left(n_{i}x\right),\]
is now viewed as a discrete Fourier transform, and if we take the
inverse transform over the horocyclic points: $z_{m}=x_{m}+iY$, $1-Q\le m\le Q$,
where $x_{m}=\frac{1}{2Q}(\frac{1}{2}-m)$ we get: \begin{align*}
\frac{1}{2Q}\sum_{1-Q}^{Q}\hat{f}_{i}(z_{m})e(-n_{i}x_{m})= & \frac{1}{2Q}\sum_{m=1-Q}^{Q}\,\sum_{l=\underset{{l_{i}\neq0}}{-M(Y)}}^{M(Y)}\frac{c_{i}(l)}{\sqrt{|l_{i}|}}W_{\frac{k}{2}sgn(l_{i}),iR}(4\pi|l_{i}|Y)e\left(l_{i}x_{m}-n_{i}x_{m}\right)\\
= & \sum_{l=\underset{{l_{i}\neq0}}{-M(Y)}}^{M(Y)}\frac{c_{i}(l)}{\sqrt{|l_{i}|}}W_{\frac{k}{2}sgn(l_{i}),iR}\left(4\pi|l_{i}|Y\right)\frac{1}{2Q}\sum_{1-Q}^{Q}e\left(l_{i}x_{m}-n_{i}x_{m}\right)\\
= & \frac{c_{i}(n)}{\sqrt{|n_{i}|}}W_{\frac{k}{2}sgn(n_{i}),iR}\left(4\pi|n_{i}|Y\right),\end{align*}
where we used the fact that \begin{eqnarray*}
\frac{1}{2Q}\sum_{1-Q}^{Q}e\left(l_{i}x_{m}-n_{i}x_{n}\right) & = & e\left(\frac{l_{i}-n_{i}}{4Q}\right)\frac{1}{2Q}\sum_{1-Q}^{Q}e\left(-\left(l_{i}-n_{i}\right)\frac{m}{2Q}\right)=\delta_{nl}.\end{eqnarray*}
Now we also have $f_{i}(z_{m})=\chi_{mi}f_{I(m,i)}(z_{mi}^{*}),$
where (by (\ref{eq:autmorphy_cond_weight})) \[
\chi_{mi}=j_{\sigma_{i}}(z_{m};k)^{-1}j_{T_{i}^{-1}U_{I(m,i)}^{-1}\sigma_{I(m,i)}}(z_{mi}^{*};k)w(T_{i}^{-1}U_{I(m,i)}^{-1},\sigma_{I(m,i)})v(T_{i}^{-1}U_{I(m,i)}^{-1},k).\]
 Hence\begin{eqnarray*}
\frac{c_{i}(n)}{\sqrt{|n_{i}|}}W_{sgn(n_{i})\frac{k}{2},iR}(4\pi|n_{i}|y) & = & \frac{1}{2Q}\sum_{1-Q}^{Q}f_{i}(z_{m})e(-n_{i}x_{m})+[[\epsilon]]\\
 & = & \frac{1}{2Q}\sum_{1-Q}^{Q}\chi_{mi}f_{I(m,i)}(z_{mi}^{*})e(-n_{i}x_{m})+[[\epsilon]]\\
 & = & \sum_{j=1}^{\kappa}\sum_{\underset{{n+\alpha_{i}\neq0}}{l=-M_{0}}}^{M_{0}}c_{j}(l)V_{nl}^{ij}+2[[\epsilon]],\end{eqnarray*}
where \begin{eqnarray*}
V_{nl}^{ij} & = & \frac{1}{\sqrt{|l_{j}|}}\frac{1}{2Q}\sum_{\underset{{I(m,i)=j}}{1-Q}}^{Q}\chi_{mi}W_{sgn(l_{j})\frac{k}{2},iR}(4\pi|l_{j}|y_{mj}^{*})\times e\left(l_{j}x_{mj}^{*}\right)e(-n_{i}x_{m}).\end{eqnarray*}
We then define $\tilde{V}_{nl}^{ij}$ by \[
\tilde{V}_{nl}^{ij}=V_{nl}^{ij}-\delta_{nl}\delta_{ji}\frac{1}{\sqrt{|n_{i}|}}W_{sgn(n_{i})\frac{k}{2},iR}(4\pi|n_{i}|Y),\]
and if we neglect the error of magnitude $\epsilon$ we finally obtain
a (well-conditioned) linear system \[
CV=0,\tag{*}\]
which must be satisfied by the Fourier coefficients of $f.$ Here
$V$ is the matrix $\tilde{V}_{nl}^{ij}$ and $C$ is the vector $c_{i}(n),$
both depending on $R$ and $Y$. The basic idea of Phase 1 is now
to locate eigenvalues $R$ together with sets of corresponding Fourier
coefficients, $c_{i}(n)$, by solving ({*}) repeatedly for different
$R$'s, and seeking those values of $R$ for which the solution vectors
$C=C(R,Y)$ are stable under changes of $Y$. That is, if $C(R,Y_{1})\approx C(R,Y_{2})$
(for $Y_{1}$and $Y_{2}<Y_{0}$, for some suitable constant $Y_{0}$)
we take it as an indication that the $R$ is close to an eigenvalue
and that the components of $C$ are close to the corresponding Fourier
coefficients. For more details and justifications see \cite[sect.\ 3.2]{stromberg:04:1}.

We note that the Phase 2 algorithm (i.e. computation of a larger set
of Fourier coefficients) also generalizes to non-zero weight in a
similar manner (cf.~\cite[sect.\ 3.3]{stromberg:04:1}).

\section{Numerical results\label{sec:Numerical-results}}

The experimental excursions have been directed towards three essentially
different subjects, but, in each, we have worked in an \emph{exploratory}
spirit.

\begin{itemize}
\item First we tried to get an over-all picture of the distribution of small
to middle-range sized eigenvalues on $\PSLZ$ (and the eta multiplier)
for {}``large'' weights, e.g. $R\in[0,14]$ and $k\in[0.1,6].$ 
\item Second, we continuously {}``turned off'' the multiplier system $v_{\eta}^{2k}$
on $\PSLZ$ by letting the weight $k\rightarrow0$ and studied the
varying distribution of eigenvalues, $N_{k}(T),$ as well as the formation
of the continuous spectrum.
\item There are some cases where arithmeticity plays a role even in nonzero
weight. We studied the Shimura correspondence between weight zero
forms on $\Gamma_{0}(2)$ and weight one half forms on $\Gamma_{0}(4)$.
And we also studied weight one forms on $\PSLZ$ which correspond
to Hecke eigenforms on $\Gamma_{0}(144)$ with a Dirichlet character. 
\end{itemize}

\subsection{Varying weight\label{sub:Varying-weight}}

\begin{figure}

\caption{\label{fig:Weights_0.1to6}Section of eigenvalues with $0<R\le14,$
and weight $0.1\le k\le6.$}

\includegraphics[scale=0.75]{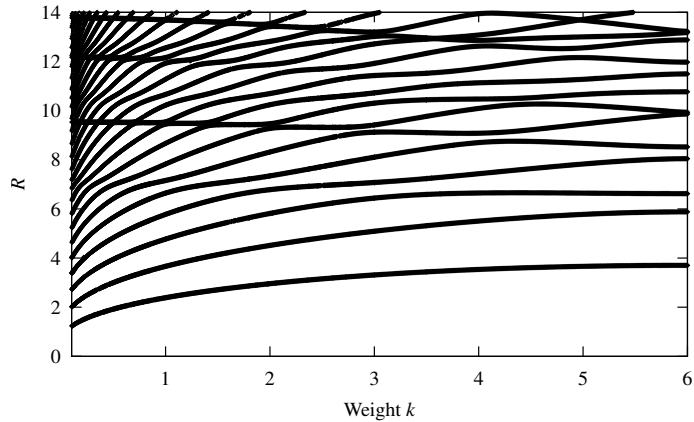}
\end{figure}

The first experiment considered was to tabulate the first few eigenvalues
(up to $R=14$) for $\PSLZ$ and the multiplier system $v_{\eta}^{2k}$,
of weight $k\in(0,6]$ (cf.~Section \ref{sub:Maass-operators-and-symmetry}).
We made the computations for $k\in[0.1,\,6]$ with a grid size of
$0.01,$ and the results are presented in Figure \ref{fig:Weights_0.1to6}.
We stress here, that the arcs in Figure \ref{fig:Weights_0.1to6}
terminate at $k=0.1$; it is not excluded (and actually expected)
that they might go lower.

For some examples of eigenvalues for {}``large'' weights see Table
\ref{tab:large_wt}. This data should be compared with data obtained
by M\"uhlenbruch, \cite[p.\ 160]{muhlenbruch:03}, who used a completely
different method (with much less accuracy). We note here that as $R$
increases, the negative Fourier coefficients seem to grow rapidly
in magnitude (as compared to the positive ones, with the normalization
$c(1)=1$) for {}``large'' weights. 

We believe that we have found all eigenvalues with $(R,k)\in[0,14]\times[0.1,\,6]$.
This belief is founded on the {}``continuity'' of the resulting
graphs $R_{j}(k)$ (cf.~Figure \ref{fig:Weights_0.1to6}), where
$R_{j}(k)$ is the $j$-th eigenvalue at weight $k$, considered as
a function of $k$. By standard results (e.g.~\cite[p.\ 149]{bruggeman:varying_weightIII})
$R_{j}(k)$ should be a real analytic function in this interval. 

\begin{figure}

\caption{\label{fig:Weights_0.1to6_compare}Comparison with the theoretical
lower bounds in $k\in[0.1,6]$. }

\includegraphics[scale=0.75]{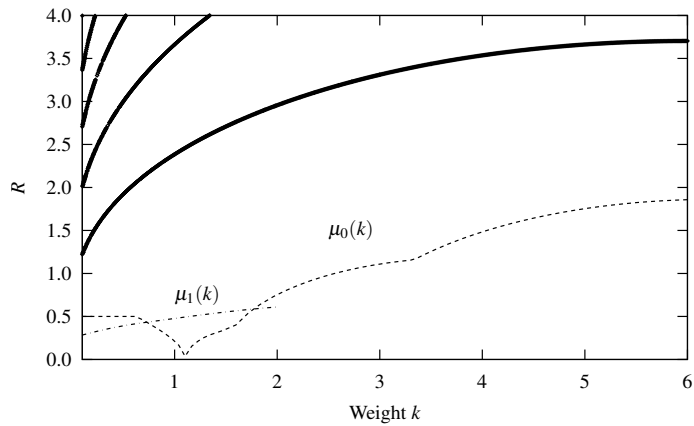}
\end{figure}

\label{ref:small_eigenv}Remember that, for $k\ge0$, the smallest
eigenvalue, $\lambda_{min},$ corresponds to the function $F(z)=y^{\frac{k}{2}}\eta(z)^{2k}$.
A lower bound for the \emph{second smallest} eigenvalue, $\lambda_{0}(k)$,
is discussed in \cite[p.\ 183]{bruggeman:varying_weightIII}. Bruggeman
finds two such bounds, both positive, which he calls $\mu_{0}(k)$
and $\mu_{1}(k)$ ($\mu_{1}$ is better than $\mu_{0}$ in a certain
interval $I\subset[0,2]$.) Figure \ref{fig:Weights_0.1to6_compare}
shows a comparison between the $R$-values corresponding to these
bounds and the smallest experimentally found eigenvalues in the interval
$k\in[0.1,\,6]$; we see that Bruggeman's bounds can be improved quite
a bit. 

\begin{table}[H]
\begin{threeparttable}
\centering
\caption{Eigenvalues for $N=1$, $v=v_{\eta}$}
\label{tab:large_wt}
\begin{tabular}{l|l}
\begin{tabular}{d{12}@{\hspace{2mm}}c@{\hspace{2mm}}c}%
\multicolumn{3}{l}{ $k=5.0$ }\\
\multicolumn{1}{c}{$R$}&
\multicolumn{1}{c}{$|c(-1)|$\tnote{a}}&
\multicolumn{1}{c}{$H(y1,y2)$}\\
\hline\noalign{\smallskip}
   3.66240686691           & 2E+3      & 8E-09 \\
   5.77698688079           & 6E+3      & 1E-09 \\
   6.64285171609           & 1E+4      & 2E-09 \\
   7.82634704661           & 7E+4      & 8E-07 \\
   8.66620831839           & 8E+4      & 1E-08 \\
   9.45156176778           & 4E+4      & 4E-09 \\
    10.21802876612           & 9E+4      & 2E-08 \\
    10.65897262925           & 2E+5      & 5E-09 \\
    11.27526358329           & 2E+5      & 2E-08 \\
    12.15792337439           & 5E+6      & 2E-06 \\
    12.55403510011           & 3E+5      & 3E-09 \\
    13.00123950671           & 4E+4      & 2E-09 \\
    13.67542640619           & 8E+5      & 8E-09 \\
    13.71353384347           & 4E+6      & 4E-07 \\
    14.47039277248           & 6E+5      & 1E-09 \\
    15.03845367363           & 1E+6      & 6E-09 \\
    15.39856858318           & 1E+6      & 1E-09 \\
    15.85705128856           & 7E+4      & 7E-09 \\
    16.14536205683           & 5E+6      & 2E-07 \\
    16.45061260141           & 2E+6      & 1E-08 \\
    16.93043847901           & 2E+7      & 2E-08 \\
    17.51562192888           & 2E+5      & 4E-09 \\
    17.59022138300           & 1E+6      & 6E-10 \\
    18.13826107361           & 7E+6      & 2E-08 \\
   18.32637702289           & 5E+5      & 5E-09 \\
    18.76341585136           & 2E+6      & 5E-09 \\
    19.16629116326           & 4E+6      & 8E-10 \\
    19.67214438521           & 3E+7      & 2E-07 \\
    19.68520099819           & 1E+6      & 3E-10 \\
    20.00524829746           & 2E+6      & 4E-09 \\
    20.38266630653           & 3E+6      & 3E-10 \\
    20.67297062056           & 6E+6      & 4E-08 \\
    20.97339376061           & 6E+6      & 8E-10 \\
\end{tabular} 
&
\begin{tabular}{d{12}@{\hspace{2mm}}c@{\hspace{2mm}}c}%
\multicolumn{3}{l}{$k=5.25$}\\
\multicolumn{1}{c}{$R$}&
\multicolumn{1}{c}{$|c(-1)|$\tnote{a}}&
\multicolumn{1}{c}{$H(y1,y2)$}\\
\hline\noalign{\smallskip}
    3.68037312372           & 3E+3      & 3E-08 \\
    5.82067054942           & 9E+3      & 6E-09 \\
    6.63460520751           & 2E+4      & 3E-09 \\
    7.90867228426           & 2E+5      & 9E-09 \\
    8.61646891946           & 1E+5      & 8E-09 \\
    9.56930344151           & 7E+4      & 5E-09 \\
    10.15656706121           & 2E+5      & 2E-09 \\
    10.70911890024           & 2E+5      & 4E-10 \\
    11.34046324165           & 4E+5      & 2E-08 \\
    12.11839521329           & 4E+6      & 2E-06 \\
    12.65021958486           & 4E+5      & 9E-09 \\
    13.02622821839           & 8E+4      & 2E-09 \\
    13.56022943627           & 1E+6      & 1E-07 \\
    13.87057635696           & 7E+5      & 4E-07 \\
    14.48204838116           & 1E+6      & 7E-05 \\
    15.09966704087           & 4E+7      & 1E-07 \\
    15.38981845044           & 1E+6      & 4E-09 \\
    15.94059443942           & 1E+4      & 3E-09 \\
    16.09999759486           & 8E+6      & 7E-06 \\
    16.52671557073           & 5E+6      & 1E-09 \\
    16.90856097808           & 2E+7      & 5E-08 \\
    17.53730159778           & 1E+6      & 3E-10 \\
    17.74142373355           & 7E+6      & 4E-09 \\
    18.02022951826           & 6E+6      & 9E-10 \\
    18.37970066644           & 8E+5      & 2E-09 \\
    18.90587158951           & 1E+7      & 8E-09 \\
    19.09726131554           & 1E+7      & 5E-10 \\
    19.66894569593           & 5E+6      & 3E-10 \\
    19.73195996101           & 7E+6      & 4E-10 \\
    20.12609436572           & 2E+6      & 3E-10 \\
    20.35571778301           & 1E+7      & 7E-10 \\
    20.71020380483           & 6E+6      & 8E-09 \\
    20.88321504381           & 1E+7     & 2E-09 \\
\end{tabular}
\end{tabular}

\begin{tablenotes}
\item[a] The normalization we have used here is the usual $c(1)=1$.
\end{tablenotes}

\end{threeparttable}
\end{table}

%

\subsection{Small weights}

The investigation of eigenvalues for small weights $k$ has been done
in the interval $R\in[0,20]$, and $k\in\left\{ 10^{-j}\,|1\le j\le12\right\} \subset\left[10^{-12},0.1\right].$
We believe that most cusp forms were found. Let us use the notation
\[
\lambda_{j}(k)\]
for the $j$-th discrete eigenvalue at weight $k,$ and $\phi_{j}(k)$
for the corresponding cusp form. It is then a basic fact that $\lambda_{j}(k)$
depends continuously on $k,$ but it can also be shown that for $k\in(0,12)$
$\lambda_{j}$ is even real analytic in $k.$ That is, $\phi_{j}(k)$
belong to an {}``analytic family'' in the terminology of Bruggeman
(see \cite{bruggeman:varying_weightIII,bruggeman:94}). In connection
with this, it should also be noted that our experiments support the
statement in Observation 71 in \cite{muhlenbruch:03} that the first
few cusp forms at weight $0$ do not belong to an analytic family
of cusp forms defined in the interval $(-12,12)$. Indeed we find
that in the range considered we actually seem to have $\lambda_{j}(k)\rightarrow0$
as $k\rightarrow0$ which is consistent with Bruggeman's result (cf.
case ii) b) in Prop. 2.17 in \cite[p. 149]{bruggeman:varying_weightIII}).

Our experiments indicate that for fixed small $k$, the {}``generic''
cusp forms $\phi_{j}(k)$ can be divided into two classes: 

\begin{itemize}
\item $C(k),$ and
\item $E(k).$ 
\end{itemize}
Here $C(k)$ consists of functions $\phi_{j}(k)$ such that $\lambda_{j}(k)$
is close to an eigenvalue of a cusp form at weight $0$ and the Fourier
coefficients are close to the corresponding coefficients of the weight-zero
cusp form. 

$E(k)$ on the other hand consists of functions $\phi_{j}(k)$ such
that the Fourier coefficients are close to the Fourier coefficients
of the Eisenstein series $E(z,s)$ where $\lambda_{j}(k)=s(1-s).$ 

The typical difference between the Fourier coefficients at weight
$k$ and weight $0$ are in both cases basically of order $k$; for
the forms in $C(k)$, the distance between $\lambda_{j}(0)$ and the
corresponding discrete eigenvalue at weight $0$ is much smaller than
$k$. 

The {}``generic'' in the statement above means that we exclude certain
places where the families $\phi_{j}(k)$ changes character between
$C(k)$ and $E(k).$ In these cases we have a situation of \emph{almost}
multiplicity 2, and in too coarse resolution it actually looks like
the analytic families intersect.

Weyl's law For non-trivial $\eta$-multiplier and a \emph{fixed} non-zero
weight $k\in(0,12)$ on $\PSLZ$ is \begin{equation}
N_{k}(T)=\sharp\{ R\le T,\,\textrm{weight }k\}=\frac{T^{2}}{12}-\frac{T}{\pi}\ln\left|1-e^{\frac{k\pi}{6}i}\right|+S(T)+O(1),\label{eq:weyls_law_fixedk}\end{equation}
 for $T\ge1$ (cf.~\cite[p. 466]{hejhal:lnm1001} with $\kappa_{0}=0$).
Our experiments seem to suggest that as $k\rightarrow0$, the main
contribution is proportional to the factor $\ln\left|1-e^{\frac{k\pi}{6}i}\right|$,
and indeed it is easy show that the $O(1)$ term is even uniformly
bounded in $k$. 

To obtain asymptotics for $S(T)$ when $k\rightarrow0$ (and $T$
is bounded) is a bit more involved. We want to generalize \cite[thm.\ 2.29, p.\ 468]{hejhal:lnm1001}
by following the pointers provided in \cite[proof of thm.\ 2.29, p.\ 468]{hejhal:lnm1001}.
Basically, our goal is to single out any terms that grow as $k\rightarrow0$
in the estimates for $S(T)$. 

By careful bookkeeping we see that these terms are all of the form
$\ln\left|1-e^{\frac{k\pi}{6}i}\right|$ and this kind of terms only
show up when we apply the functional equation for the logarithmic
derivative of the Selberg zeta function, $Z\left(s\right)$ (\cite[thm.\ 2.18, p.\ 441]{hejhal:lnm1001}).
To obtain the necessary estimate of $\frac{Z'}{Z}\left(s\right)$
we study the integral in the left hand side of \cite[prop. 8.6, p.\ 121]{hejhal:lnm548}.
Let the assumptions of \cite[ass. 8.5 p.\ 120]{hejhal:lnm548} and
\cite[p.\ 125 row 5]{hejhal:lnm548} apply, e.g.~$s=\sigma+iT$,
$U=T+10$ and $\frac{1}{2}\le\frac{1}{2}+\epsilon=\sigma_{1}<\sigma$.
Through the integral estimates \cite[prop.\ 8.7-8, 8.10-11]{hejhal:lnm548}
we consider either $\frac{Z'}{Z}\!\left(\xi\right)-2\ln\left|1-e^{\frac{k\pi}{6}i}\right|$
or $\frac{Z'}{Z}\!\left(\xi\right)$, depending on whether the functional
equation is applied along this part of the path or not and then collect
the extra terms at the end. It turns out that we consider the modified
integrand in all parts to the left of $\Re\left(\xi\right)=\frac{1}{2}-\epsilon$.
With this modification, all estimates in \cite[prop.\ 8.7-8, 8.10-14, thm.\ 8.15 and 8.17]{hejhal:lnm548}
go through except for the addition of a term $R_{m}(s)=\frac{2}{\ln x}\ln\left|1-e^{\frac{k\pi}{6}i}\right|\, I_{m}(s)$
where $I_{m}(s)$ is the integral\begin{eqnarray*}
I_{m}(s) & = & \frac{1}{2\pi i}\int_{\gamma}\frac{x^{\xi-s}-x^{2\left(\xi-s\right)}}{\left(\xi-s\right)^{2}}d\xi,\end{eqnarray*}
where the path is ${\gamma=\gamma}_{U}=\left[\frac{1}{2}-\epsilon-iU,-A-iU\right]\cup\left[-A-iU,-A+iU\right]\cup\left[-A+iU,\frac{1}{2}-\epsilon+iU\right].$
The integrand has no poles to the right of $\Re\left(\xi\right)=\frac{1}{2}-\epsilon$
and is hence equal to the integral from $\frac{1}{2}-\epsilon-iU$
to $\frac{1}{2}-\epsilon+iU$. Elementary estimates now show that
\[
\left|I_{m}\left(s\right)\right|\le\frac{2x^{\left(\frac{1}{2}-\epsilon-\sigma\right)}}{2\pi}\int_{-U}^{U}\frac{1}{\epsilon^{2}+\left(t-T\right)^{2}}dt\le\frac{x^{-2\epsilon}}{\pi\epsilon^{2}}\int_{-U}^{U}\frac{1}{1+\left(\frac{t-T}{\epsilon}\right)^{2}}dt\le\frac{x^{-2\epsilon}}{\epsilon}.\]
Using $x=T^{\frac{1}{4}}$ (cf.~\cite[p.\ 138, row -6]{hejhal:lnm548})
together with $\epsilon=\frac{1}{\ln x}$ (cf.~\cite[ass.\ 8.16(b), p.\ 130]{hejhal:lnm548})
we see that \[
\left|R_{m}\left(s\right)\right|\le x^{-\frac{2}{\ln x}}\ln\left|1-e^{\frac{k\pi}{6}i}\right|\ln x\,\frac{1}{\ln x}=e^{-2}\ln\left|1-e^{\frac{k\pi}{6}i}\right|.\]
This error term is now to be added to $\frac{Z'}{Z}$ in \cite[thm. 8.15]{hejhal:lnm548}
and subtracted from the right hand side of \cite[thm. 8.17]{hejhal:lnm548}
and finally in \cite[p.\ 135 (8.18)]{hejhal:lnm548} we need to add
$-\int_{\frac{1}{2}}^{2}\Im\left[R_{m}\left(s\right)\right]$ to the
expression of $\pi S\left(T\right)$. Thus the analogue of \cite[thm. 8.1, p.\ 119]{hejhal:lnm548}
and \cite[thm. 2.29, p.\ 468]{hejhal:lnm1001} is $S\left(T\right)=O\left(\frac{T}{\ln T}\right)+O\left(1\right)\ln\left|1-e^{\frac{k\pi}{6}i}\right|$
(with implied constants independent of $k$) and inserted into (\ref{eq:weyls_law_fixedk})
we conclude the following Weyl's law:\begin{alignat}{1}
N_{k}(T) & =\frac{T^{2}}{12}-\ln\left|1-e^{\frac{k\pi i}{6}}\right|\left[\frac{T}{\pi}+O\left(1\right)\right]+O\left(\frac{T}{\ln T}\right),\,\label{eq:Weyls_law_wto0}\end{alignat}
uniformly for $0<k\le0.1$ (say) and $T\ge2$ and the $O\left(1\right)$
term is bounded in magnitude by ${\frac{3}{2\pi}e}^{-2}\approx0.065$.

\begin{figure}

\caption{Plot of Weyl's law with constant $T=20$ and weight $k\rightarrow0$\label{fig:Plot-of-Weyls}}

\includegraphics[scale=0.75]{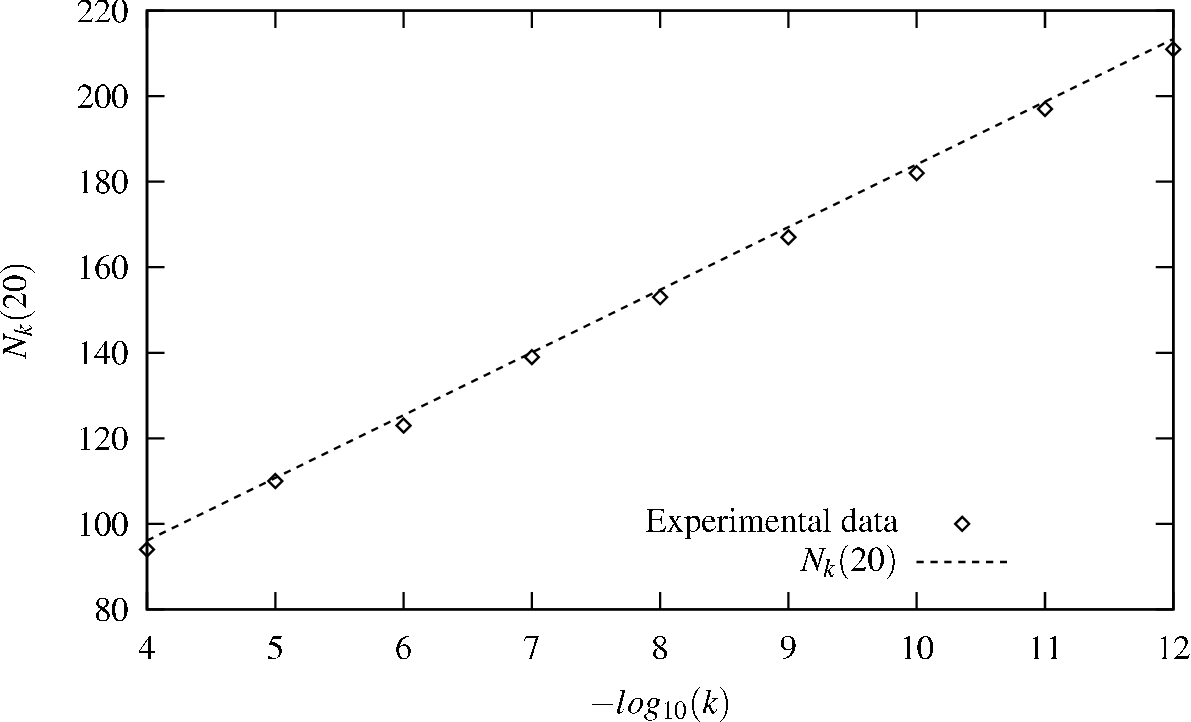}
\end{figure}

We have computed $N_{k}(T)$ experimentally for $k=10^{-j},$$j=4,\ldots,12,$
and $0\le T\le20$. Figure \ref{fig:Plot-of-Weyls} shows a picture
of the experimental values compared to the theoretical values obtained
by (\ref{eq:Weyls_law_wto0}) (with the $O$-terms neglected) and
the difference indeed seems to be constant over this range. 

\begin{figure}

\caption{Section of eigenvalues with $9\le R\le14,$ and $1E-9\le k\le1E-7.$
\label{fig:eigenvalues-go-down}}

\includegraphics[scale=0.75]{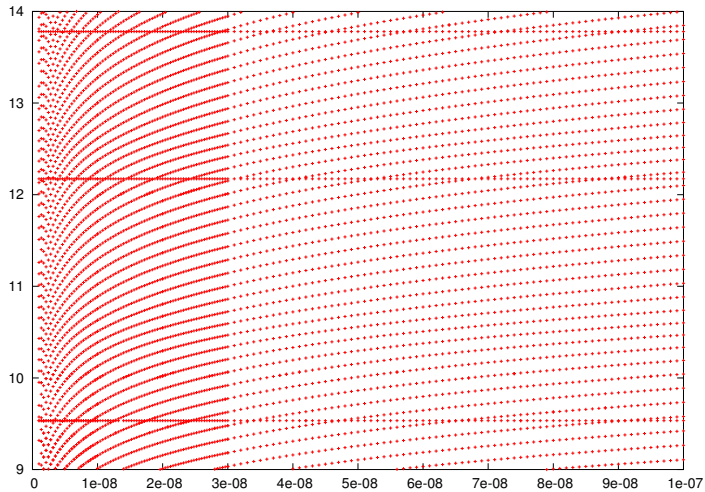}
\end{figure}

From the form of the Weyl's law above, we also see that the successive
spacings $\Delta_{n}(k)=R_{n+1}(k)-R_{n}(k)$ should look about like
$\frac{1}{\frac{dN_{k}}{dT}}$, i.e. $\Delta_{n}(k)\approx\frac{1}{\frac{T}{6}+\frac{1}{\pi}\left|\ln\frac{k\pi}{6}\right|}\sim\frac{\pi}{\left|\ln\frac{k\pi}{6}\right|}$
as $k\rightarrow0$. Figure \ref{fig:eigenvalues-go-down} provides
a nice illustration of this fact, where it is clearly seen that the
average spacings are almost constant for small $k$ and this constant
is roughly proportional to $\frac{1}{|\ln k|}$. 

The mean gap over 100 typical cases ($R_{n}$, $5\le n\le155$) turned
out to be $0.117$ at $k=10^{-9},$ $0.107$ at $k=10^{-10}$, $0.101$
at $k=10^{-11}$ and $0.095$ at $k=10^{-12}$. The relative quotients,
$\Delta_{n}(k)\left[\frac{R_{n}}{6}+\frac{1}{\pi}\left|\ln\frac{k\pi}{6}\right|\right]$,
are $1.019,$ $0.990$, $1.001$ and $0.9934,$ respectively and it
is not inconceivable that one obtains $1$ in the (logarithmic) limit.

\subsubsection{Level repulsion}

\begin{figure}

\caption{\label{fig:Level-repulsion-atR13}Level repulsion at $R=13.779\ldots$}

\includegraphics[scale=0.5]{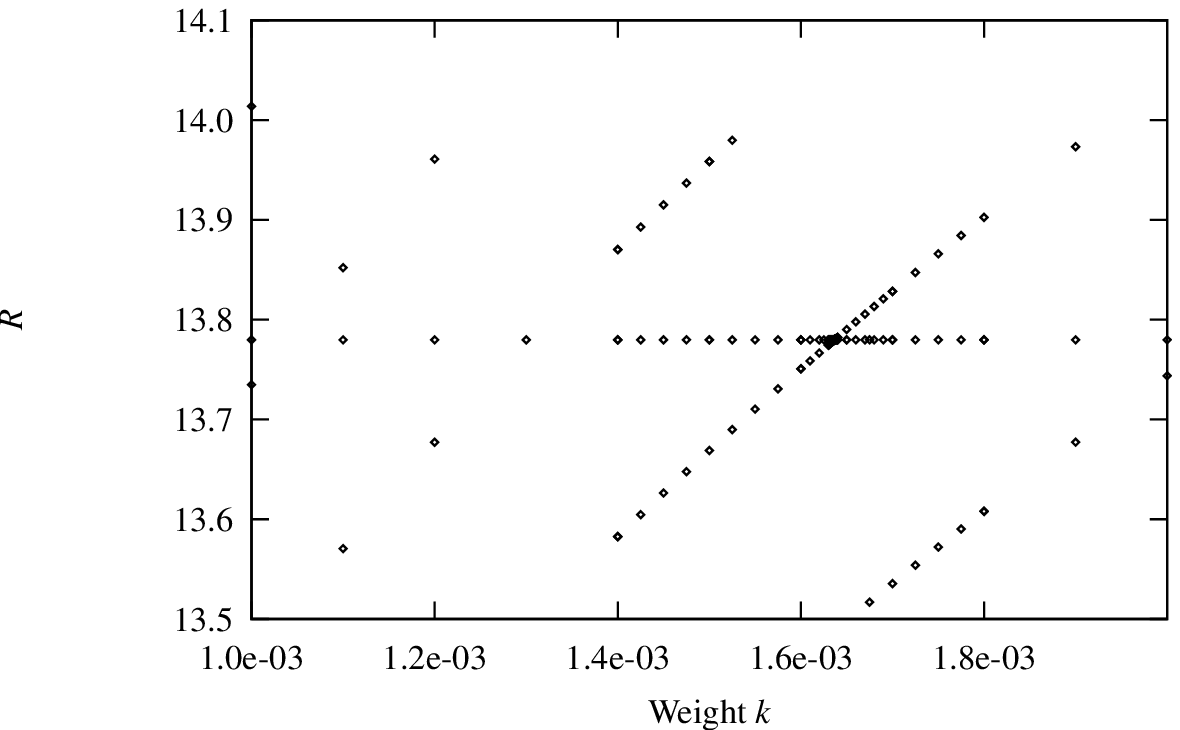}~~\includegraphics[scale=0.5]{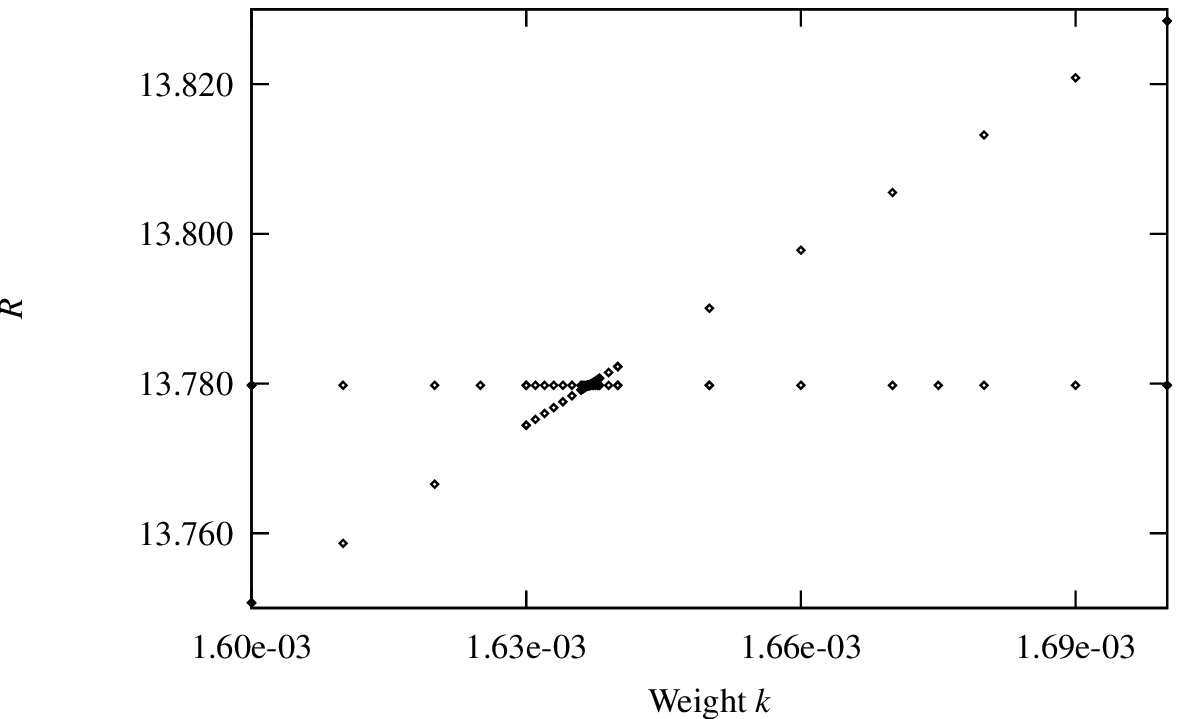}\\
\includegraphics[scale=0.5]{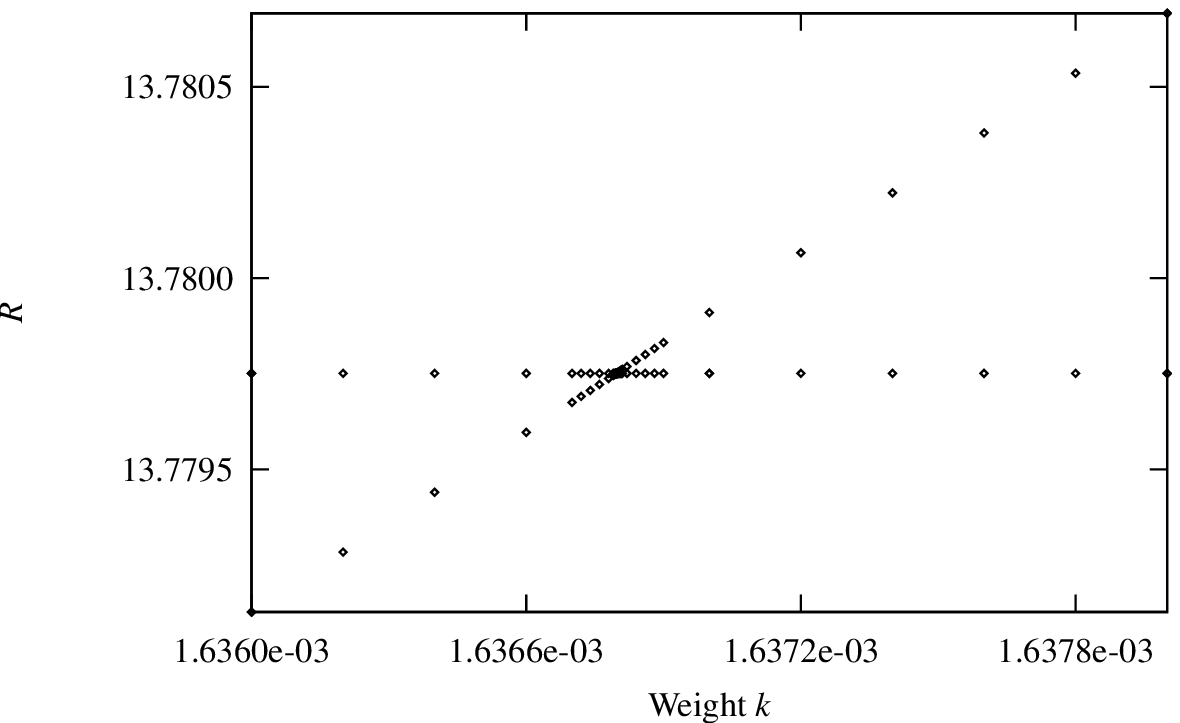}~~\includegraphics[scale=0.5]{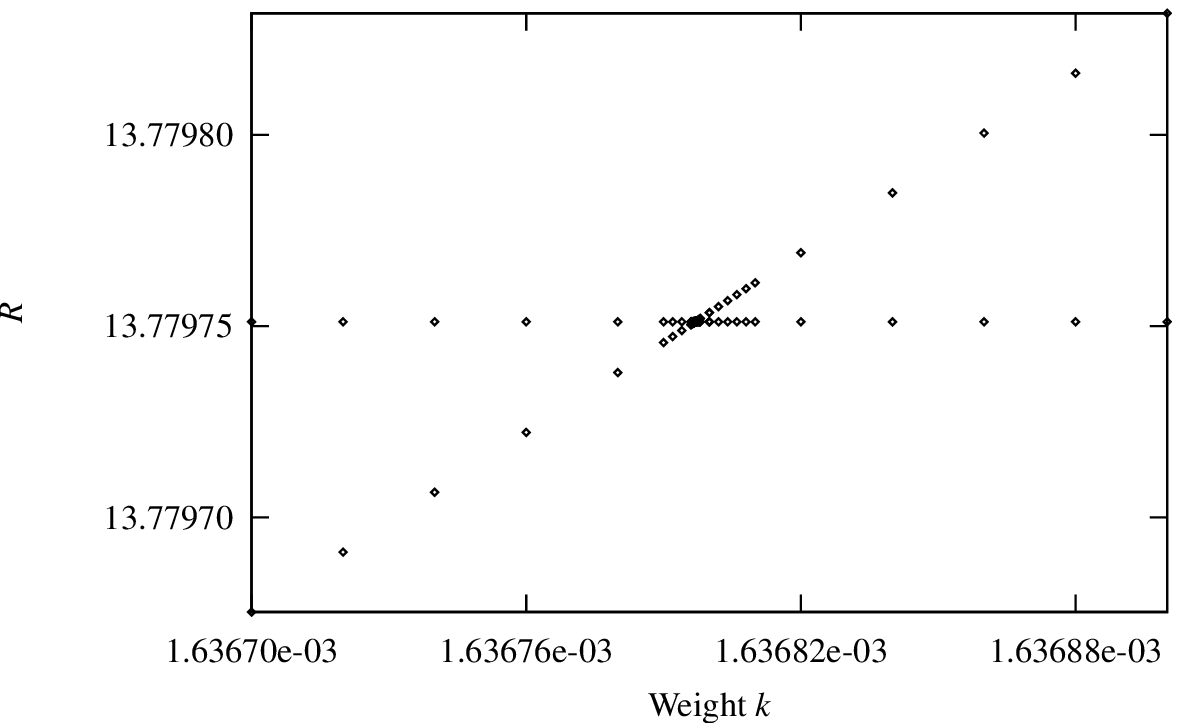}\\
\includegraphics[scale=0.5]{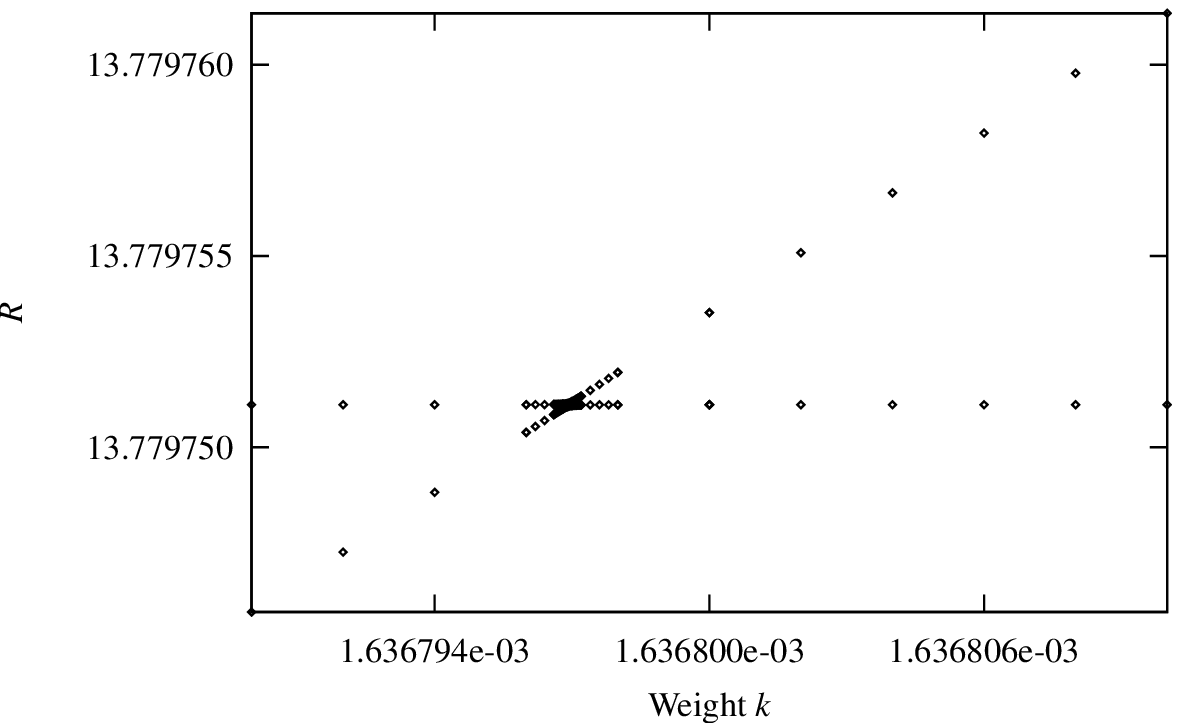}~~\includegraphics[scale=0.5]{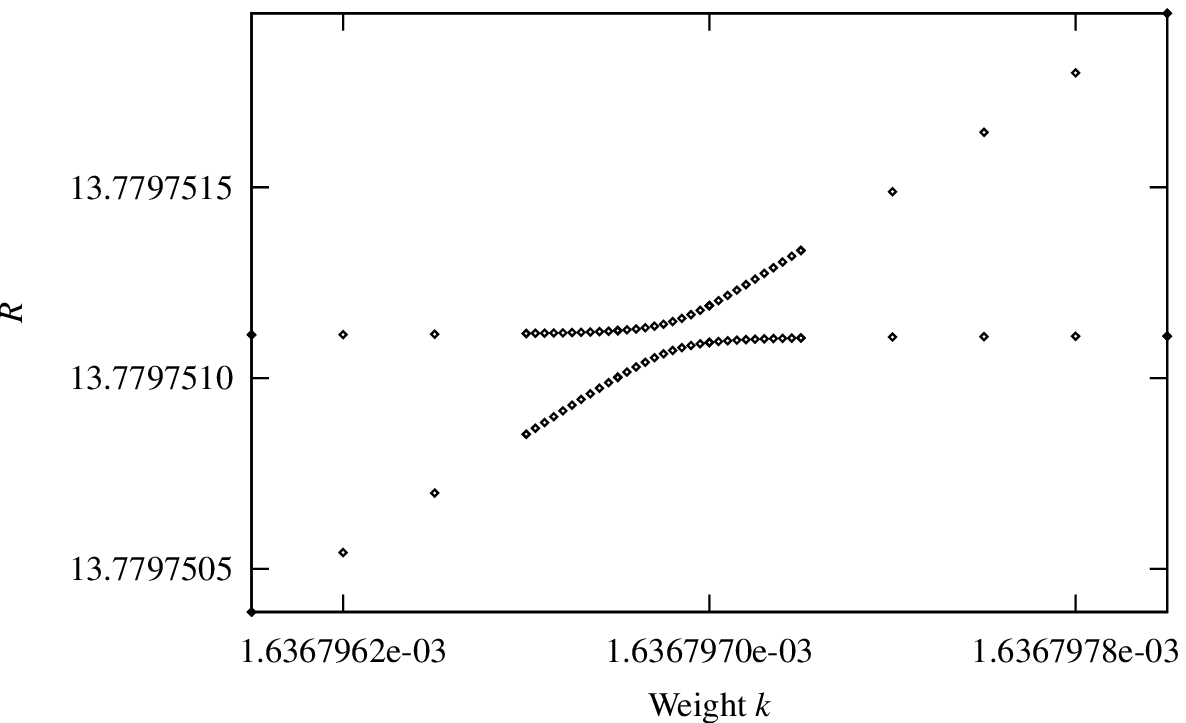}
\end{figure}

\begin{figure}

\caption{\label{fig:Level-repulsion-atR9}Level repulsion at $R=9.533\ldots$}

\includegraphics[scale=0.5]{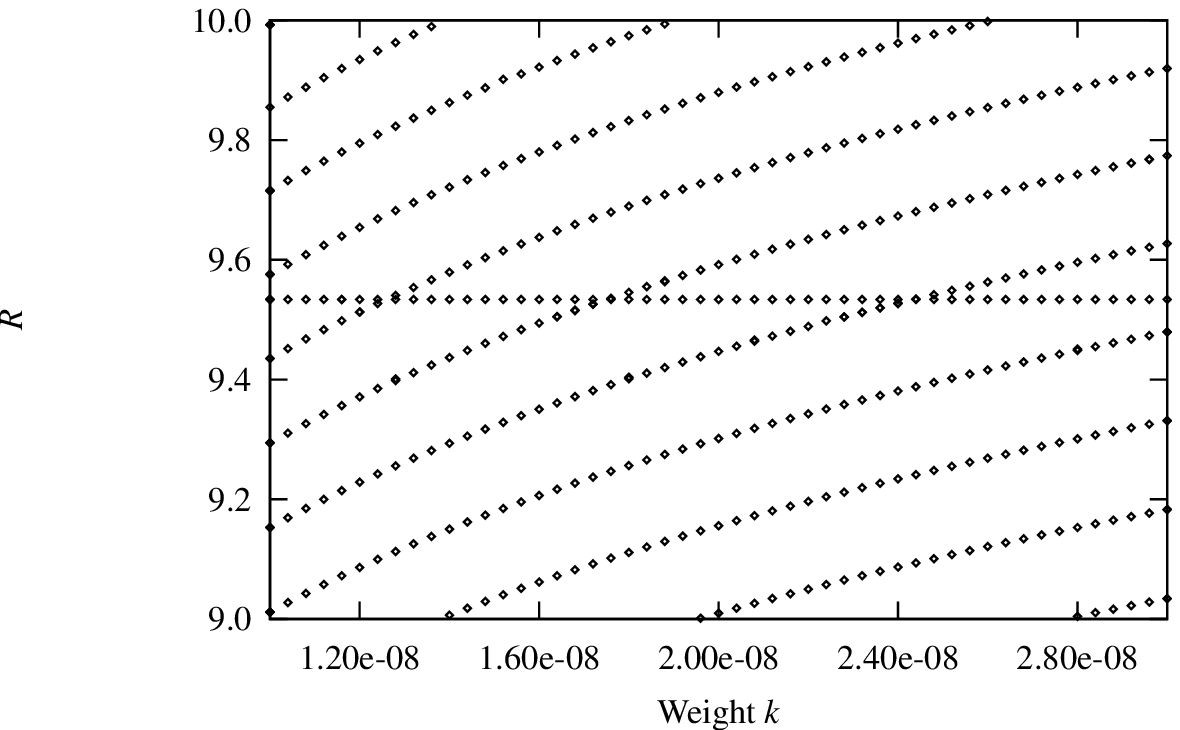}~~\includegraphics[scale=0.5]{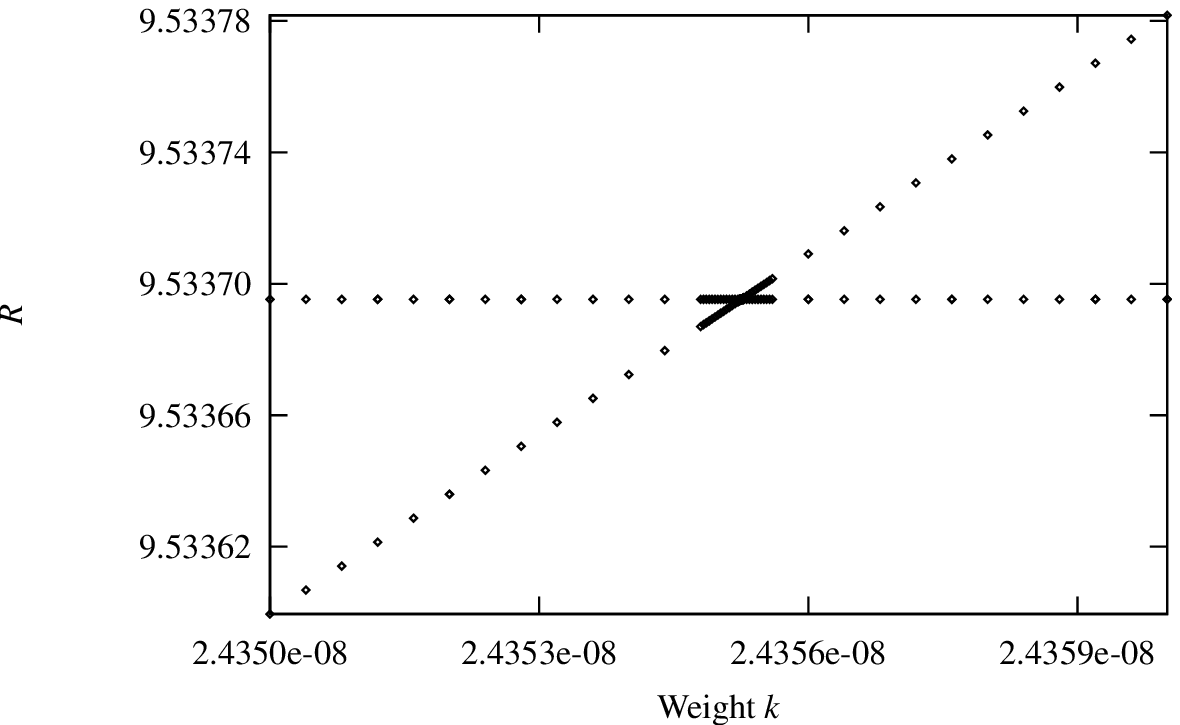}\\
\includegraphics[scale=0.5]{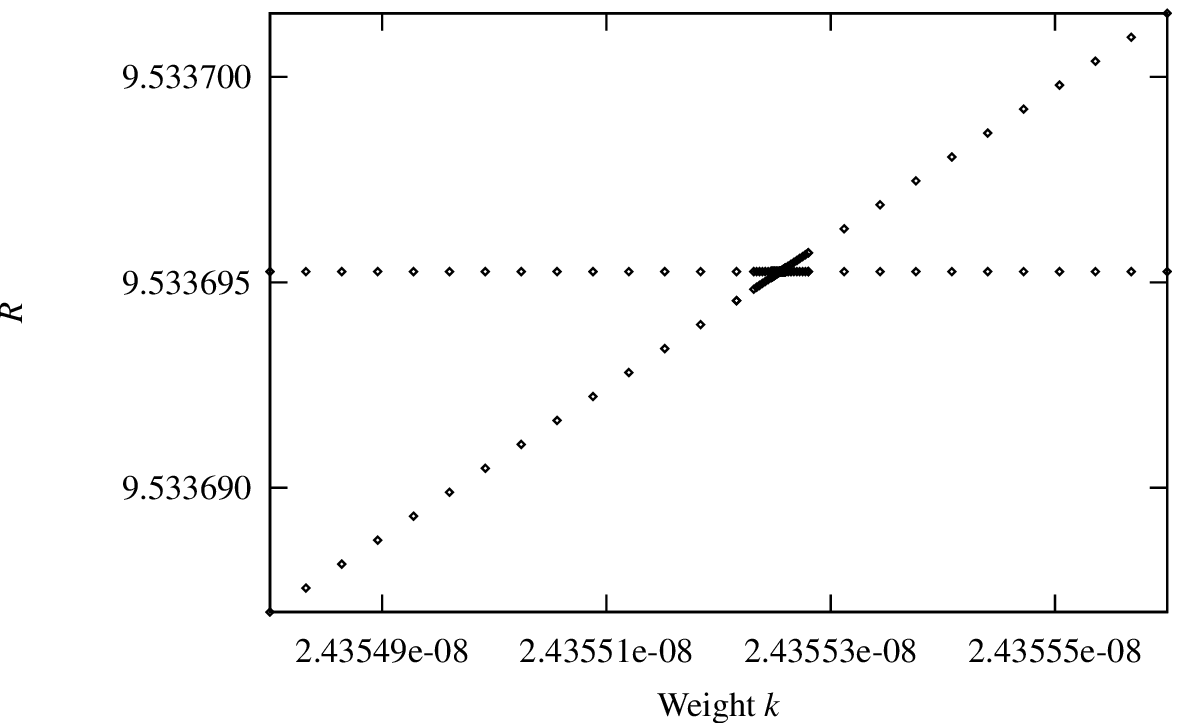}~~\includegraphics[scale=0.5]{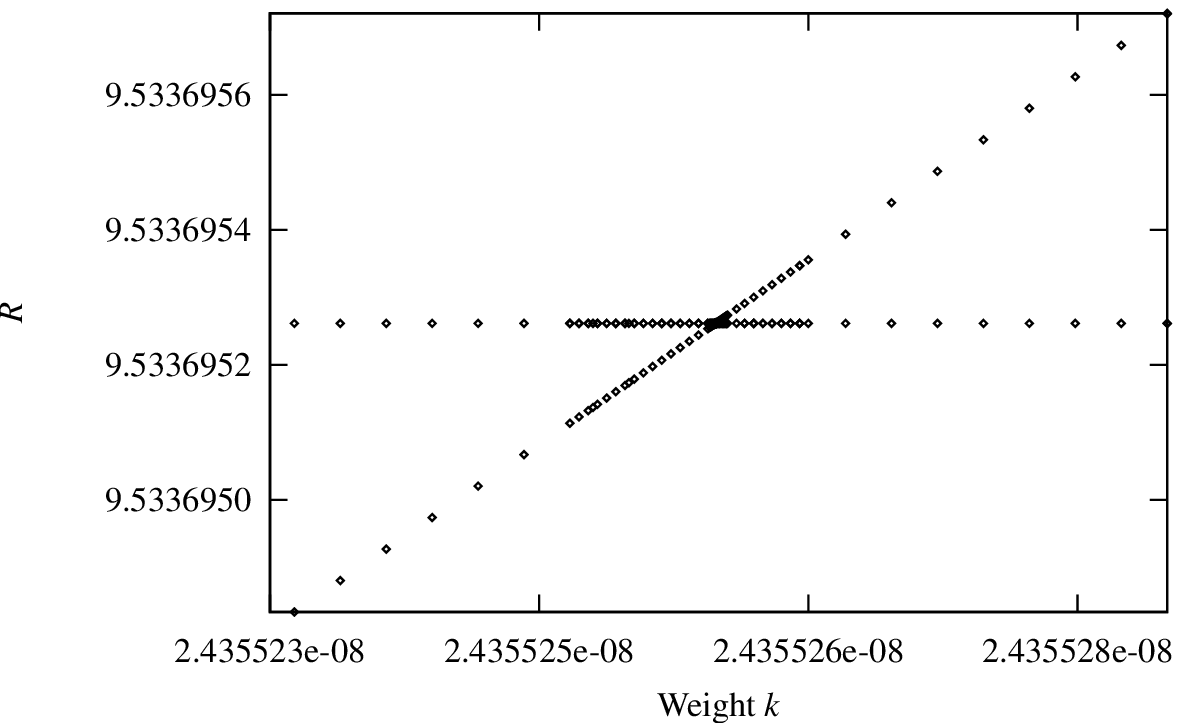}\\
\includegraphics[scale=0.5]{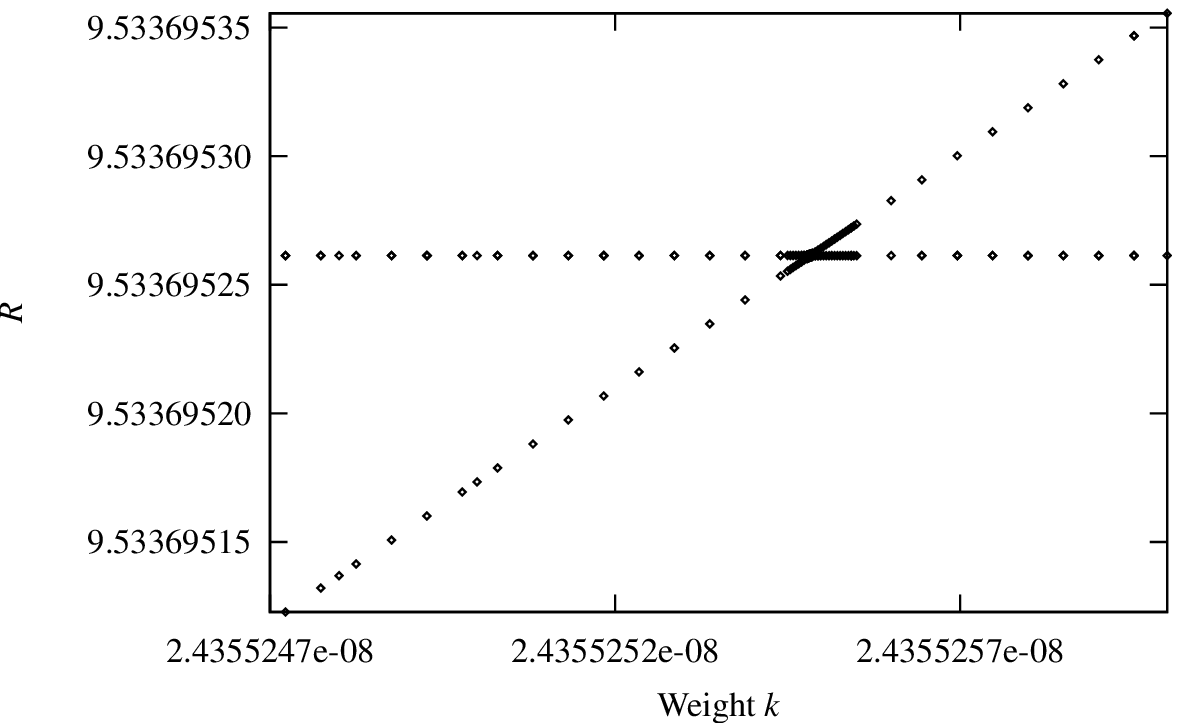}~~\includegraphics[scale=0.5]{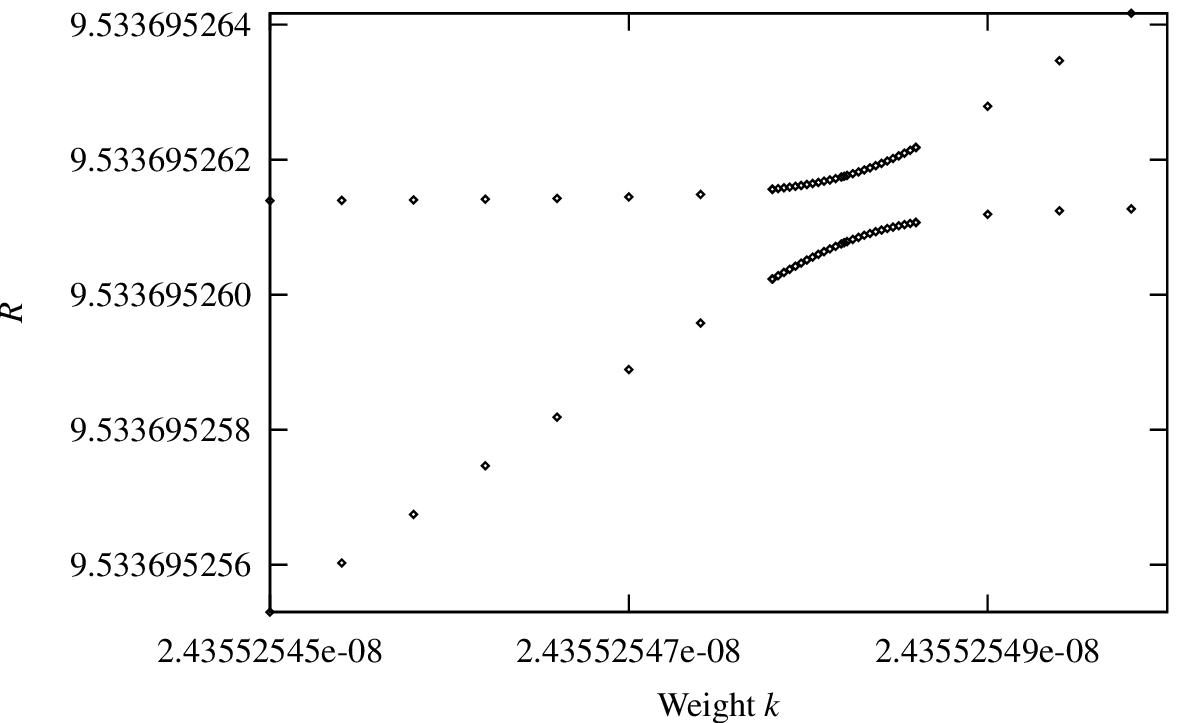}
\end{figure}

From figures like \ref{fig:eigenvalues-go-down} one may be tempted
to think that there are horizontal lines corresponding to cusp forms
at weight $0$ which crosses the lines that are going down (i.e. corresponding
to the Eisenstein series at weight $0$). This is not the case! If
we look closer we will see that there is actually {}``level repulsion''
here, i.e. the horizontal {}``cusp-form-line'' turns down before
the {}``near crossing'' and becomes an {}``Eisenstein-series-line''
and the previous {}``Eisenstein-series-line'' turns into a {}``cusp-form-line''.
See also Figures \ref{fig:Level-repulsion-atR13} and \ref{fig:Level-repulsion-atR9}.
More precisely formulated: if there is a {}``near crossing'' at
the weight $k_{0}$ close to the eigenvalue $R_{0}\approx R_{j}(k_{0})\approx R_{j+1}(k_{0}),$
then there are two analytic families $\phi_{j}(k)$ and $\phi_{j+1}(k)$
such that for some $\delta>\epsilon>0$:\[
\phi_{j}(k)\in\begin{cases}
C(k), & k\in[k_{0}+\epsilon,k_{0}+\delta],\\
E(k), & k\in[k_{0}-\delta,k_{0}-\epsilon],\end{cases}\]
and \[
\phi_{j+1}(k)\in\begin{cases}
E(k), & k\in[k_{0}+\epsilon,k_{0}+\delta],\\
C(k), & k\in[k_{0}-\delta,k_{0}-\epsilon],\end{cases}\]
and in the interval $(k_{0}-\epsilon,k_{0}+\epsilon)$ both families
display a mixing between the two types $E(k)$ and $C(k).$ In fact,
the Fourier coefficients of $\phi_{j+1}$ converge (as $k\rightarrow k_{0}-\epsilon$)
to values close to the Fourier coefficients of $\phi_{j}$ for $k>k_{0}+\epsilon$
and vice versa. Since the two functions also need to be orthogonal
it is clear that the Fourier coefficients exhibit {}``wild'' behavior
in the small interval surrounding the {}``near crossing''. Note
also that, as $k\rightarrow0$, \emph{all} $\phi_{j}(k)\in E(k)$
and $R_{j}(k)\rightarrow0.$ 

See Table \ref{tab:cuspf_and_eisen_comp} for examples of Fourier
coefficients corresponding to eigenfunctions of types $E(k)$ and
$C(k)$, close to an avoided crossing at weight $k=9.044605824E-8$.
Table \ref{tab:cuspf_and_eisen_comp2} illustrates the agreement between
the Fourier coefficients of a more generic cusp form in $E(k)$ and
the corresponding coefficients of the Eisenstein series (appropriately
normalized) at weight $0$. The level of agreement is striking to
put it mildly; likewise in Table \ref{tab:cuspf_and_eisen_comp} for
the $C(k)$ eigenfunction. The {}``1 for 1'' nature of this convergence
\emph{in the presence of a limiting continuous spectrum} seems not
to have been suspected earlier. Cf. \cite[thm.\ 6.6 and cor.\ 6.9]{MR1052555}
and \cite{MR803365}.

The fact that the system seems to avoid accidental degeneracies by
means of level repulsion and avoided crossings is in agreement with
the Wigner-von Neumann theorem, cf.~\cite{wigner-vn}. 

\subsection{Lifts at weight 1\label{sub:Lifts-at-weight-data}}

As we saw in Section \ref{sub:Lifts-at-weight} we could prove the
existence of certain Hecke relations at weight $1$ (e.g. (\ref{eq:wt_1_hecke_rel_pos}),
(\ref{eq:wt_1_hecke_rel__neg_1}) and \ref{eq:wt_1_hecke_rel__neg_2})).
Tables \ref{tab:N1R47709} and \ref{tab:N1R36624} contain numerical
verifications of these relations. Table \ref{tab:egenvN1_w1} contains
a list of computed eigenvalues on $\MAS\left(\Gamma_{0}(1),v_{\eta}^{2},1,\lambda\right)$,
and the eigenvalues corresponding to cosine CM-forms are indicated.
In these cases, we have computed the actual error since we know the
exact eigenvalues:\[
R_{k}=\frac{2\pi k}{\ln\left(7+2\sqrt{12}\right)},\, k\in\mathbb{Z}^{+}.\]

Note that the actual error is in general much smaller than the error-parameter
which is basically $H(Y_{1},Y_{2})=|c(2)-c'(2)|+|c(3)-c'(3)|+|c(4)-c'(4)|,$
where $c(n)$ is computed with $Y_{1}$ and $c'(n)$ with $Y_{2}.$

\begin{table}
\begin{threeparttable}
\centering
\caption[]{Fourier coefficients for a CM-form \\
$f \in  \MAS (\Gamma_0(1),v_{\eta}^{2},1,4.770984191561)$}
\label{tab:N1R47709}
	\begin{tabular}[t]{d{1}d{12}cl}
$n$ &                  & $c(n)/c(0)$\tnote{a}   & Error    \\
\hline\noalign{\smallskip}
 0&   1.755930576575    &   &       \\
 1&   1.000000000000    &   &       \\
 2&  -1.755930576574    & $c(27)$  &  0.4E-08     \\
 3&   1.571810322167    & $c(40)$  &  0.1E-08     \\
 4&   1.755930576575    & $c(53)$  &  0.7E-08     \\
 5&  -1.770268323978    & $c(66)$  &  0.3E-09     \\
 6&  -2.474798320759    & $c(79)$  &  0.1E-08     \\
 7&   0.000000000000    & $c(92)$ &   0.4E-08    \\
 8&   0.346240855507    & $c(105)$  & 0.5E-08      \\
 9&   3.510179255561    & $c(118)$  & 0.3E-08      \\
10&   1.116593241680    & $c(131)$  & 0.4E-08       \\
11&  -0.000000000001    & $c(144)$  & 0.3E-08       \\
12&   0.000000000001    & $c(157)$  & 0.1E-07      \\
13&  -3.019229958496    & $c(170)$  & 0.8E-08       \\
14&  -1.186431979458    & $c(183)+c(1)$  & 0.3E-08      \\
15&  -3.079783541463    & $c(196)$  & 0.5E-08       \\
  &   &              $-c(n)c(-1)/R^{2}/c(0)$    &   \\
-1&    5.055064268188        & \tnote{b}          &  0.1E-11          \\
-2&  -11.180067729976        & $c(21)$         &  0.4E-07             \\
-3&    0.000000000001        & $c(32)$          & 0.2E-07              \\
-4&   16.472675660354        & $c(43)$          & 0.2E-08              \\
-5&   16.729030199659        &  $c(54)$          &0.5E-07               \\
-6&   13.098490835617        &   $c(65)$         &0.5E-08               \\
-7&   16.046414105740        &   $c(76)$         &0.1E-07               \\
-8&   -0.000000000023        &  $c(87)$          &0.1E-07               \\
-9&   13.340740291248        &   $c(98)$         &0.3E-07               \\
-10&    0.000000000006        &   $c(109)$         &0.4E-07               \\
-11 &    2.151215034503        &  $c(120)$          &0.7E-08               \\
-12 &    2.878852009050        &  $c(131)-c(1)$          & 0.2E-07              \\
-13 &    0.000000000039        &  $c(142)$          &0.6E-08               \\
-14 &    7.496795049955        &  $c(153)$         & 0.4E-07              \\
-15 &  -12.830881007618        &  $c(164)$         & 0.4E-07              \\
\end{tabular} 
\begin{tablenotes}
\item[a] This quotient is deduced from formula (\ref{eq:w1_coeff_prop1}) or (\ref{eq:w1_coeff_prop2}) on p.\ \pageref{eq:w1_coeff_prop1}.
\item[b]  $c(-1)^{2}=-R^{2}c(0)(c(10)-c(0))$
\end{tablenotes}
\end{threeparttable}
\end{table}
\begin{table}
\begin{threeparttable}
\centering
\caption[]{Fourier coefficients for a non-CM-form \\
$f \in  \MAS (\Gamma_0(1),v_{\eta}^{2},1,3.66240686698667)$}
\label{tab:N1R36624}
	\begin{tabular}[t]{d{1}d{12}cl}
$n$ &                  & $c(n)/c(0)$\tnote{a}   & Error    \\
\hline\noalign{\smallskip}
 0&  -1.352193685534     &   &       \\
 1&   1.000000000000    &   &       \\
 2&  -1.697113317091  & $c(27)$  &  0.5E-08     \\
 3&  -0.057989599353   & $c(40)$  &  0.4E-07     \\
 4&   2.461764397786    & $c(53)$  &  0.1E-07     \\
 5&   0.510856433057   & $c(66)$  &  0.6E-08     \\
 6&  -0.952325903762  & $c(79)$  &  0.8E-08     \\
 7&  -1.660630683908  & $c(92)$ &   0.2E-07    \\
 8&  -2.343382246022   & $c(105)$  & 0.9E-08      \\
 9&   1.271097206907   & $c(118)$  & 0.1E-07      \\
10&  -0.203512820511   & $c(131)$  & 0.6E-08       \\
11&   2.110622602834   & $c(144)$  & 0.5E-08       \\
12&   2.170616908700     & $c(157)$  & 0.4E-08      \\
13&   0.449799127363   & $c(170)$  & 0.4E-07       \\
14&   0.612654661780   & $c(183)+c(1)$  & 0.4E-07      \\
15&  -1.684453740441   & $c(196)$  & 0.1E-07       \\
16&   0.400312170289   & $c(209)$  & 0.2E-07       \\
17&  -2.868395110060   & $c(222)$  & 0.1E-07       \\
18&  -1.931595991172   & $c(235)$  & 0.2E-07      \\
19&  -0.591212766919    & $c(248)$  & 0.1E-07       \\
20&   0.792151138999    & $c(261)$  & 0.9E-08       \\
21&  -1.717242193922  & $c(274)$    & 0.3E-07      \\
22&   1.369169138277   & $c(287)$    & 0.2E-07      \\
23&   2.007854712832    & $c(300)$    & 0.8E-09       \\
24&   0.447826147902  & $c(313)$    & 0.1E-07      \\
25&   3.051006373828  & $c(326)$    & 0.4E-08      \\
26&  -0.032419986064      & $c(339)$    & 0.6E-08       \\
27&   1.255081531820     & $c(352)+a(2)$  & 0.2E-07      \\
28&  -0.707047087424   & $c(365)$    & 0.3E-07       \\
29&   1.272283355260    & $c(378)$    & 0.1E-07       \\
30&  -0.184187214400  & $c(391)$    & 0.2E-07     \\
\end{tabular} 
\begin{tablenotes}
\item[a] This quotient is deduced from formula (\ref{eq:w1_coeff_prop1}) or (\ref{eq:w1_coeff_prop2}) on p.\ \pageref{eq:w1_coeff_prop1}.
\end{tablenotes}
\end{threeparttable}
\end{table}

\begin{table}
\begin{threeparttable}
\centering
\caption[]{Eigenvalues for $\Mas{\PSLZ,1,v^2_{\eta}}$} 
\label{tab:egenvN1_w1}
\begin{tabular}[t]{d{14}cc}
$R$ &   $H(y_1,y_2)$ & True error \tnote{a}      \\
\hline\noalign{\smallskip}
2.38549209578045 & 1E-12 & 1E-15 \tnote{b} \\
3.66240686698667 & 1E-11\\
4.77098419156091 & 1E-13 & 1E-14 \tnote{b} \\
5.77698688078694 & 6E-12\\
6.64285171613711 & 1E-11\\
7.15647628734173 & 1E-11 & 4E-13 \tnote{b} \\
7.82634704540775 & 1E-13\\
8.66620831896793 & 7E-14\\
9.45156176783224 & 7E-12\\
9.54196838312186 & 2E-13 & 6E-14 \tnote{b} \\
10.21802876776059 & 3E-13\\
10.65897262920241 & 2E-12\\
11.27526358349387 & 2E-13\\
11.92746047890219 & 7E-11 & 5E-14 \tnote{b} \\
12.15792337422149 & 7E-13\\
12.55403509998720 & 1E-13\\
13.00123950642372 & 9E-13\\
13.67542640643589 & 2E-13\\
13.71353384358095 & 1E-12\\
14.31295257468268 & 1E-12 & 1E-14 \tnote{b} \\
14.47039277253940 & 1E-12\\
15.03845367358721 & 2E-13\\
15.39856858348441 & 1E-12\\
15.85705128717333 & 3E-10\\
16.14536205734475 & 3E-11\\
16.45061260131967 & 4E-11\\
16.69844467046304 & 4E-12 & 1E-13 \tnote{a} \\
16.93043847896222 & 4E-15\\
17.51562192885174 & 2E-10\\
17.59022138305996 & 1E-10\\
18.13826107340244 & 3E-11\\
18.32637702205910 & 4E-11\\
18.76341585146817 & 1E-11\\
\end{tabular}
\begin{tablenotes}
\item[a] For CM-forms, the true error is computed with respect to 
the eigenvalue $R_k=\frac{2\pi k}{\ln(\eta_0)}$, where $\eta_0=7+2\sqrt{12}$.
\item[b] These forms correspond to CM-forms.
\end{tablenotes}
\end{threeparttable}
\end{table}

\subsection{Half integer weight\label{sub:Half-integer-weight}}

We now consider the case of $\Gamma_{0}(4)$ and the $\theta$-multiplier
system for weight $k=\frac{1}{2}$. The aim of our investigation in
this case was to study the Shimura lift, and in particular to investigate
the dimensions of the spaces of half integer weight forms. As remarked
at the end of section \ref{sub:The-Shimura-correspondence} several
properties of the Shimura correspondence were observed numerically,
and in the original version of these notes, \cite[ch.\ 2]{stromberg:thesis},
we formulated a number of experimentally inspired conjectures. Except
for the question of injectivity these conjectures have now been resolved
with Propositions \ref{thm:existence_of_shim_corr} and \ref{thm:Shimura_corr_prop}.
(With the obvious correction in the first conjecture: the dimension
of $\MAS(\Gamma_{0}(4),v_{\theta},\frac{1}{2},R)$ corresponding to
an oldspace is \emph{at least} two dimensional.)

Tables \ref{tab:N4R6889} and \ref{tab:N4R46143} contain examples
of Fourier coefficients at weight $\frac{1}{2}$. Note especially
in table \ref{tab:N4R46143} the agreement with Proposition \ref{thm:Shimura_corr_prop}
c) displayed by the Fourier coefficients $c(n)$ for $n\equiv5,1\mod8$
respectively. Here, it is known that $8.92287648699174$ and $12.09299487507860$
on $\Gamma_{0}(2)$ correspond to the eigenvalues $1$ and $-1$,
respectively, with respect to the involution $\omega_{2}$. 

Table \ref{tab:fourier_coeff_comparison} contains a comparison of
Fourier coefficients computed both from forms on $\Gamma_{0}(2)$
via (\ref{eq:shimura_corr_coeff}) and computed directly. Additional
Fourier coefficients for the weight $0$ forms are available in Table
\ref{tab:fourier_coeff_w_half_supp}. 

\begin{acknowledgement*}
This paper is based on Chapter 2 of my Ph.D. thesis \cite{stromberg:thesis}
but with additional theoretical discussions in e.g. Section \ref{sub:The-Shimura-correspondence}
and \ref{sub:Hecke-operators-non-trivial}. I am grateful to my advisor
Dennis Hejhal and my second advisor Andreas Strömbergsson for many
giving discussions and lots of valuable comments and suggestions. 
\end{acknowledgement*}
\bibliographystyle{amsplain}
\bibliography{/home/fredrik/Documents/matematik/refs}

\providecommand{\bysame}{\leavevmode\hbox to3em{\hrulefill}\thinspace}
\providecommand{\MR}{\relax\ifhmode\unskip\space\fi MR }
\providecommand{\MRhref}[2]{%
  \href{http://www.ams.org/mathscinet-getitem?mr=#1}{#2}
}
\providecommand{\href}[2]{#2}
\begin{thebibliography}{10}

\bibitem{apostol}
T.~M. Apostol, \emph{Modular {F}unctions and {D}irichlet {S}eries in {N}umber
  {T}heory}, Springer-Verlag, 1976.

\bibitem{arakawa:2}
T.~Arakawa, \emph{Shimura correspondence for {M}aass wave forms and {S}elberg
  zeta functions}, Proceedings of the conference "automorphic forms and
  representations of algebraic groups over local fields" (Hiroshi Saito and
  Tetusya Takahashi, eds.), RIMS Kyoto Univ., 2003, pp.~1--14.

\bibitem{atkinlehner}
A.~O.~L. Atkin and J.~Lehner, \emph{Hecke operators on
  \uppercase{$\Gamma_{0}(m)$}}, Math. Ann. \textbf{185} (1970), 134--160.

\bibitem{helen:deform_published}
H.~Avelin, \emph{Deformation of {$\Gamma_{0}(5)$} {C}usp {F}orms}, Math. Comp.
  (2006), posted on October 4, 2006, PII S 0025-5718(06)01911-9 (to appear in
  print).

\bibitem{biro:00}
A.~Biró, \emph{Cycle integrals of {M}aass forms of weight 0 and fourier
  coefficients of {M}aass forms of weight 1/2}, Acta Arithmetica \textbf{94}
  (2000), no.~2, 103--152.

\bibitem{andreas:effective_comp}
A.~Booker, A.~Strömbergsson, and A.~Venkatesh, \emph{{E}ffective {C}omputations
  with {M}aass {F}orms}, In preparation.

\bibitem{borevich-shafarevich}
Z.I. Borevich and I.R. Shafarevich, \emph{Zahlentheorie}, Birkhäuser, 1966.

\bibitem{MR808915}
R.~W. Bruggeman, \emph{Modular forms of varying weight. {I}}, Math. Z.
  \textbf{190} (1985), no.~4, 477--495. \MR{MR808915 (87c:11051)}

\bibitem{MR840831}
\bysame, \emph{Modular forms of varying weight. {II}}, Math. Z. \textbf{192}
  (1986), no.~2, 297--328. \MR{MR840831 (87k:11059)}

\bibitem{bruggeman:varying_weightIII}
\bysame, \emph{Modular forms of varying weight. {III}}, J. Reine Angew. Math.
  \textbf{371} (1986), 144--190.

\bibitem{bruggeman:94}
\bysame, \emph{Families of automorphic forms}, Birkhäuser, 1994.

\bibitem{MR993311}
D.~Bump, \emph{The {R}ankin-{S}elberg method: a survey}, Number {T}heory,
  {T}race {F}ormulas and {D}iscrete {G}roups (K.E.~Aubert et~al, ed.), Academic
  Press, 1989, pp.~49--109.

\bibitem{MR904946}
D.~Bump and J.~Hoffstein, \emph{On {S}himura's correspondence}, Duke Math. J.
  \textbf{55} (1987), no.~3, 661--691.

\bibitem{duke:88}
W.~Duke, \emph{Hyperbolic distribution problems and half-integral weight
  {M}aass forms}, Invent. Math. \textbf{92} (1988), 73--90.

\bibitem{erdelyi:53}
A.~Erd{\'e}lyi, W.~Magnus, F.~Oberhettinger, and F.~G. Tricomi, \emph{Higher
  {T}ranscendental {F}unctions. {V}ols. {I}, {II}}, McGraw-Hill, 1953.

\bibitem{farmer-lemurell}
D.~W. Farmer and S.~Lemurell, \emph{Deformations of {M}aass forms}, Math. Comp.
  \textbf{74} (2005), no.~252, 1967--1982 (electronic). \MR{MR2164106}

\bibitem{farmer-lemurell:2}
D.W. Farmer and S.~Lemurell, \emph{{M}aass forms and their {L}-functions},
  arXiv:math.NT/0506102.

\bibitem{gunning}
R.~C. Gunning, \emph{Lectures on {M}odular {F}orms}, Princeton University
  Press, 1962.

\bibitem{hardy:on_toe_representaation}
G.H. Hardy, \emph{On the representation of a number as a sum of any number of
  squares and in particular five}, Trans. Am. Math soc. \textbf{21} (1920),
  255--284.

\bibitem{MR726196}
D.~A. Hejhal, \emph{Some {D}irichlet series with coefficients related to
  periods of automorphic eigenforms. {I,II}}, Proc. Japan Acad. Ser. A, {\bf
  58} (1982), pp.\ 413--417; {\bf 59} (1983), pp.\ 335-338.

\bibitem{hejhal:lnm548}
\bysame, \emph{The {S}elberg {T}race {F}ormula for \rm {PSL}(2,$\mathbb{R}$),
  \em {V}ol.1}, Lecture Notes in Mathematics, vol. 548, Springer-Verlag, 1976.

\bibitem{hejhal:lnm1001}
\bysame, \emph{The {S}elberg {T}race {F}ormula for \rm {PSL}(2,$\mathbb{R}$),
  \em {V}ol.2}, Lecture Notes in Mathematics, vol. 1001, Springer-Verlag, 1983.

\bibitem{MR803365}
\bysame, \emph{A continuity method for spectral theory on {F}uchsian groups},
  Modular forms (R.A. Rankin, ed.), Ellis-Horwood, Chichester, 1984,
  pp.~107--140.

\bibitem{MR1052555}
\bysame, \emph{Regular {$b$}-groups, degenerating {R}iemann surfaces, and
  spectral theory}, Mem. Amer. Math. Soc. \textbf{88} (1990), no.~437, iv+138.

\bibitem{hejhal:92}
\bysame, \emph{Eigenvalues of the {L}aplacian for {H}ecke triangle groups},
  Mem. Amer. Math. Soc. \textbf{97} (1992), no.~469, vi+165.

\bibitem{hejhal:99_eigenf}
\bysame, \emph{On eigenfunctions of the {L}aplacian for {H}ecke triangle
  groups}, Emerging applications of number theory (D.~Hejhal, J.~Friedman,
  et~al., eds.), IMA Vol. Math. Appl., vol. 109, Springer, New York, 1999,
  pp.~291--315.

\bibitem{hejhal:calc_of_maass_cusp_forms}
\bysame, \emph{On the {C}alculation of {M}aass {C}usp {F}orms}, Proceedings of
  the "International School on Mathematical Aspects of Quantum Chaos II",
  Lecture Notes in Physics, Springer, 2004, to appear.

\bibitem{Iwaniec:topics}
H.~Iwaniec, \emph{Topics in \uppercase{C}lassical \uppercase{A}utomorphic
  \uppercase{F}orms}, American Mathematical Society, 1997.

\bibitem{katok-sarnak}
S.~Katok and P.~Sarnak, \emph{{H}eegner points, cycles and {M}aass forms},
  Israel Journal of Mathematics \textbf{84} (1993), 193--227.

\bibitem{khuri:fc_of_hilbert}
K.~Khuri-Makdisi, \emph{On the {F}ourier coefficients of nonholomorphic
  {H}ilbert modular forms of half-integral weight}, Duke Math. J. \textbf{84}
  (1996), 399--452.

\bibitem{knopp:modular}
M.I. Knopp, \emph{Modular {F}unctions in {A}nalytic {N}umber {T}heory}, Markham
  Publishing Company, 1970.

\bibitem{MR81j:10030}
W.~Kohnen, \emph{Modular forms of half-integral weight on {$\Gamma
  \sb{0}(4)$}}, Math. Ann. \textbf{248} (1980), no.~3, 249--266.

\bibitem{MR84b:10038}
\bysame, \emph{Newforms of half-integral weight}, J. Reine Angew. Math.
  \textbf{333} (1982), 32--72.

\bibitem{MR783554}
\bysame, \emph{Fourier coefficients of modular forms of half-integral weight},
  Math. Ann. \textbf{271} (1985), no.~2, 237--268. \MR{MR783554 (86i:11018)}

\bibitem{MR629468}
W.~Kohnen and D.~Zagier, \emph{Values of {$L$}-series of modular forms at the
  center of the critical strip}, Invent. Math. \textbf{64} (1981), no.~2,
  175--198. \MR{MR629468 (83b:10029)}

\bibitem{kojima:95}
H.~Kojima, \emph{Shimura correspondence of {M}aass wave forms with half
  integral weight}, Acta Arithmetica \textbf{69} (1995), no.~4, 367--385.

\bibitem{kojima:00}
\bysame, \emph{On the fourier coefficients of {M}aass wave forms of half
  integral weight over an imaginary quadratic field}, J. reine angew. Math.
  \textbf{526} (2000), 155--179.

\bibitem{kojima:04}
Hisashi Kojima, \emph{On the {F}ourier coefficients of {H}ilbert-{M}aass wave
  forms of half integral weight over arbitrary algebraic number fields}, J.
  Number Theory \textbf{107} (2004), no.~1, 25--62. \MR{MR2059949
  (2005c:11055)}

\bibitem{MR1282723}
S.~Lang, \emph{Algebraic {N}umber {T}heory}, second ed., Graduate Texts in
  Mathematics, vol. 110, Springer-Verlag, New York, 1994.

\bibitem{maass:49}
H.~Maass, \emph{Über eine neue {A}rt von nichtanalytischen automorphen
  {F}unktionen und die {B}estimmung {D}irichletscher {R}eihen durch
  {F}unktionalgleichungen}, Math. Ann. \textbf{121} (1949), 141--183.

\bibitem{MR0065583}
\bysame, \emph{Die {D}ifferentialgleichungen in der {T}heorie der elliptischen
  {M}odulfunktionen}, Math. Ann. \textbf{125} (1952), 235--263 (1953).
  \MR{MR0065583 (16,449c)}

\bibitem{maass:modular_functions}
\bysame, \emph{Lectures on {M}odular {F}unctions of {O}ne {C}omplex
  {V}ariable}, second ed., Tata Institute of Fundamental Research Lectures on
  Mathematics and Physics, vol.~29, Tata Institute of Fundamental Research,
  Bombay, 1983.

\bibitem{MR0232968}
W.~Magnus and R.~Oberhettinger, F.~Soni, \emph{Formulas and {T}heorems for the
  {S}pecial {F}unctions of {M}athematical {P}hysics}, Third enlarged edition,
  Springer-Verlag, 1966.

\bibitem{miyake}
T.~Miyake, \emph{Modular {F}orms}, Springer-Verlag, 1997.

\bibitem{mordell:representations}
L.~J. Mordell, \emph{On the representations of a number as a sum of an odd
  number of squares}, Trans. Cambr. Phil. Soc. \textbf{22} (1919), 361--372.

\bibitem{muhlenbruch:03}
T.~M\"uhlenbruch, \emph{Systems of {A}utomorphic {F}orms and {P}eriod
  {F}unctions}, Ph.D. thesis, Universiteit Utrecht, 2003.

\bibitem{niwa:half_integral}
S.~Niwa, \emph{Modular forms of half integral weight and the integral of
  certain theta-functions}, Nagoya Math. J. \textbf{56} (1975), 147--161.

\bibitem{MR0562506}
Shinji Niwa, \emph{On {S}himura's trace formula}, Nagoya Math. J. \textbf{66}
  (1977), 183--202. \MR{MR0562506 (58 \#27781)}

\bibitem{MR0384702}
S.~J. Patterson, \emph{The {L}aplacian operator on a {R}iemann surface},
  Compositio Math. \textbf{31} (1975), no.~1, 83--107. \MR{MR0384702 (52
  \#5575)}

\bibitem{MR0419364}
\bysame, \emph{The {L}aplacian operator on a {R}iemann surface. {II}},
  Compositio Math. \textbf{32} (1976), no.~1, 71--112. \MR{MR0419364 (54
  \#7385)}

\bibitem{MR0491511}
\bysame, \emph{The {L}aplacian operator on a {R}iemann surface. {III}},
  Compositio Math. \textbf{33} (1976), no.~3, 227--259. \MR{MR0491511 (58
  \#10750)}

\bibitem{petersson:30}
H.~Petersson, \emph{Theorie der automorphen {F}ormen beliebiger reeller
  {D}imension und ihre {D}arstellung durch eine neue {A}rt {P}oincaréscher
  {R}eihen}, Math. Ann. \textbf{103} (1930), 369--436.

\bibitem{petersson:37:analytischen}
\bysame, \emph{Zur analytischen {T}heorie der {G}renzkreisgruppen. {I}.
  {G}renzkreisgruppen und {R}iemannsche {F}lächen; {T}heorie der {F}aktoren-
  und {M}ultiplikatorsysteme komplexer {D}imension}, Math. Ann. \textbf{115}
  (1937), 23--67.

\bibitem{MR0028964}
\bysame, \emph{Automorphe {F}ormen als metrische {I}nvarianten. {II}.
  {M}ultiplikative {D}ifferentiale als {G}renzwerte metrischer {I}nvarianten
  von stetig ver\"anderlicher reeller {D}imension}, Math. Nachr. \textbf{1}
  (1948), 218--257. \MR{MR0028964 (10,525g)}

\bibitem{rankin:mod}
R.~A. Rankin, \emph{Modular {F}orms and {F}unctions}, Cambridge University
  Press, 1976.

\bibitem{MR0243062}
W.~Roelcke, \emph{Das {E}igenwertproblem der automorphen {F}ormen in der
  hyperbolischen {E}bene. {I}, {II}}, Math. Ann. 167 (1966), 292--337; ibid.
  \textbf{168} (1966), 261--324. \MR{MR0243062 (39 \#4386)}

\bibitem{sarnak:82:MR750670}
P.~Sarnak, \emph{Additive number theory and {M}aass forms}, Number theory (New
  York, 1982), Lecture Notes in Math., vol. 1052, Springer, Berlin, 1984,
  pp.~286--309.

\bibitem{MR0088511}
A.~Selberg, \emph{Harmonic analysis and discontinuous groups in weakly
  symmetric {R}iemannian spaces with applications to {D}irichlet series}, J.
  Indian Math. Soc. (N.S.) \textbf{20} (1956), 47--87. \MR{MR0088511 (19,531g)}

\bibitem{shimura}
G.~Shimura, \emph{Introduction to the \uppercase{A}rithmetic \uppercase{T}heory
  of \uppercase{A}utomorphic \uppercase{F}orms}, Princeton Univ Press, 1971.

\bibitem{shimura:73:half_integral}
\bysame, \emph{On modular forms of half integral weight}, Ann. of Math.
  \textbf{97} (1973), 440--481.

\bibitem{shimura:fc_of_hilbert}
\bysame, \emph{On the {F}ourier coefficients of {H}ilbert modular forms of
  half-integral weight}, Duke Math. J. \textbf{71} (1993), 501--557.

\bibitem{shintani:75:onconstruction}
T.~Shintani, \emph{On construction of holomorphic cusp forms of half integral
  weight}, Nagoya Math. J. \textbf{58} (1975), 83--126.

\bibitem{MR803370}
H.~M. Stark, \emph{Fourier coefficients of {M}aass waveforms}, Modular forms
  (Durham, 1983), Ellis Horwood Ser. Math. Appl.: Statist. Oper. Res., Horwood,
  Chichester, 1984, pp.~263--269. \MR{MR803370 (87h:11128)}

\bibitem{stromberg:thesis}
F.~Strömberg, \emph{Computational aspects of maass waveforms}, Ph.D. thesis,
  Uppsala University, 2004.

\bibitem{stromberg:04:1}
\bysame, \emph{{M}aass waveforms on {$(\Gamma_0(N),\chi)$}, computational
  aspects}, Proceedings of the "International School on Mathematical Aspects of
  Quantum Chaos II", 2005, to appear.

\bibitem{hecke_operators_weight}
\bysame, \emph{Hecke {O}perators for {M}aass {W}aveforms on $psl(2,\mathbb{Z})$
  with {I}nteger {W}eight and {E}ta {M}ultiplier}, Preprint, 2006.

\bibitem{strombergsson:heckeoperators}
A.~Strömbergsson, \emph{The {S}elberg {T}race {F}ormula for
  ${SL}(2,\mathbb{R})$ and arbitrary real weight}, unpublished manuscript.

\bibitem{Andreas:thesis}
\bysame, \emph{Studies in the analytical and spectral theory of automorphic
  forms}, Ph.D. thesis, Uppsala University, Dept. of Math., 2001.

\bibitem{then:02}
H.~Then, \emph{{M}aass cusp forms for large eigenvalues}, Math. Comp.
  \textbf{74} (2005), no.~249, 363--381.

\bibitem{van-Lint:HeckeOperators:MR0090616}
J.~H. van Lint, \emph{Hecke operators and {E}uler products}, Drukkerij ``Luctor
  et Emergo'', Leiden, 1957. \MR{MR0090616 (19,839f)}

\bibitem{van-Lint:dedekindEta:MR0103287}
\bysame, \emph{On the multiplier system of the {R}iemann-{D}edekind function
  {$\eta $}}, Nederl. Akad. Wetensch. Proc. Ser. A. 61 = Indag. Math.
  \textbf{20} (1958), 522--527. \MR{MR0103287 (21 \#2065)}

\bibitem{wigner-vn}
J.~Von~Neumann and E.~Wigner, \emph{Über das {V}erhalten von {E}igenwerten bei
  adiabatischen {P}rocessen}, Phys. Zeit. \textbf{30} (1929), 467--470.

\bibitem{MR577010}
J.-L. Waldspurger, \emph{Correspondance de {S}himura}, J. Math. Pures Appl. (9)
  \textbf{59} (1980), no.~1, 1--132. \MR{MR577010 (83f:10029)}

\bibitem{MR646366}
\bysame, \emph{Sur les coefficients de {F}ourier des formes modulaires de poids
  demi-entier}, J. Math. Pures Appl. (9) \textbf{60} (1981), no.~4, 375--484.
  \MR{MR646366 (83h:10061)}

\bibitem{MR936998}
A.~M. Winkler, \emph{Cusp forms and {H}ecke groups}, J. Reine Angew. Math.
  \textbf{386} (1988), 187--204. \MR{MR936998 (90g:11067)}

\bibitem{MR0106888}
K.~Wohlfahrt, \emph{\"{U}ber {O}peratoren {H}eckescher {A}rt bei {M}odulformen
  reeller {D}imension}, Math. Nachr. \textbf{16} (1957), 233--256.

\end{thebibliography}

\begin{table}
\begin{threeparttable}
\centering
\caption[]{Fourier coefficients for \\
$f_{1,2} \in \MAS (\Gamma_0(4),\frac12,6.889875675945)$}
\label{tab:N4R6889}
\begin{tabular}{c@{\hspace{5mm}}c}
\multicolumn{2}{c}{Fourier coefficients for $f_1\in V^{+}$, observe that $a(n)=0$ for $n\equiv2,3 \mod 4$}\\
	\begin{tabular}[t]{d{1}d{14}}
\multicolumn{1}{c}{$n$} &
\multicolumn{1}{c}{$a(n)$}            \\
	\hline\noalign{\smallskip}
	4	&0.84219769675471           \\
	5	&0.18355821406443         \\
	8	&0.56907998524429         \\
	9	&-0.33045049673565        \\
	12      &-0.41296169213831        \\
	13      &0.60153537988057         \\
	16      &0.30482066289397      \\
	17      &-0.88689598690620      \\
	20      &0.41418282092059      \\
	21      &0.50212917023175      \\
	23      &-0.00000000000110      \\
	24      &-1.04429548341249      \\
	25      &0.28984678984391      \\
	\end{tabular}
&
	\begin{tabular}[t]{d{1}d{14}}
\multicolumn{1}{c}{$n$} &
\multicolumn{1}{c}{$a(-n)$}            \\
	\hline\noalign{\smallskip}
	 3   &   1.01825299171456   \\
	 4   &   2.18968040385979   \\
	 7   &   2.06305218095270   \\
	 8   &  -1.17157116610978   \\
	11   &  -1.02718719694121   \\
	12   &   2.29759751514531   \\
	15   &  -4.08935474990375   \\
	16   &   3.39248165496032   \\
	19   &  -1.31804633824673   \\
	20   &   1.87725570455517   \\
	23   &  -2.22265197818114   \\
	24   &   0.58816620330140   \\
	\end{tabular}
\\
\multicolumn{2}{c}{$|c(4)c(9)-c(36)|=0.2E-08$} \\
\end{tabular}
\begin{tabular}{c@{\hspace{5mm}}c}
\multicolumn{2}{c}{Fourier coefficients for $f_2 \not\in V^{+}$, observe that
  $a(n)=0$ for $n\equiv1 \mod 8$}\\
	\begin{tabular}[t]{d{1}d{14}d{2}d{14}}
\multicolumn{1}{c}{$n$} &
\multicolumn{1}{c}{$a(n)$}&
\multicolumn{1}{c}{$n$} &
\multicolumn{1}{c}{$a(n)$}  \\ 
	\hline\noalign{\smallskip}
   0&   0.000000000000  &   19&  -0.936283350934  \\
   1&   0.000000000000  &   20&   0.192087703124  \\
   2&   1.000000000000  &   21&   1.247834928335  \\
   3&  -0.725665465042  &   22&   0.411443268212  \\
   4&  -0.710754741008  &   23&  -0.018564690418  \\
   5&   0.456158224759  &   24&  -1.297582830232  \\
   6&  -1.835059236817  &   25&   0.000000000000  \\
   7&   0.481289183972  &   26&  -0.610656711780  \\
   8&   0.707106781187  &   27&  -0.179166638197  \\
   9&   0.000000000000  &   28&   0.340322845698  \\
  10&  -0.651049038069  &   29&   0.001563849679  \\
  11&  -0.470914040036  &  30&   1.462745418892  \\
  12&  -0.513122971204  &  31&   0.400763051265  \\
  13&   1.494868058150  &  32&   0.095523702475  \\
  14&   0.968484734380  &  33&   0.000000000000  \\
  15&  -0.945563574521  &  34&  -0.879782299796  \\
  16&  -1.101175502961  &  35&  -1.833603623066  \\
  17&   0.000000000000  &  36&   0.234869257223  \\
  18&   0.824250041644  &    & \\
\end{tabular}
&
\end{tabular}
\end{threeparttable}
\end{table}

\begin{table}
\begin{threeparttable}
\centering
\caption{Fourier coefficients for $f \in \MAS \left(\Gamma_0(4),\frac12,R\right)$}
\label{tab:N4R46143}
\begin{tabular}{cc}
\begin{tabular}[t]{lr}
\begin{tabular}[t]{d{1}d{14}}
\multicolumn{2}{l}{$R=4.461438243496$}\\
\multicolumn{1}{c}{$n$} &
\multicolumn{1}{c}{$a(n)$}            \\
\hline\noalign{\smallskip}

 0 &      0.00000000000000\\
 1 &      1.00000000000000\\
 2 &      0.63334968449036\\
 3 &      0.63517832947402\\
 4 &     -0.70710678118667\\
 5 &     -0.00000000000003\\
 6 &      1.28035706400142\\
 7 &     -0.90756258916698\\
 8 &     -0.44784585676555\\
 9 &      0.52643872643776\\
10 &     -0.57763498966972\\
11 &      1.12485377915641\\
12 &     -0.44913890403389\\
13 &      0.00000000000001\\
14 &      0.48078071833327\\
15 &     -1.58539012005784\\
16 &      0.50000000000015\\
17 &     -0.14882069214483\\
18 &      1.06474902295605\\
19 &      0.21268916863632\\
20 &      0.00000000000003\\
21 &      0.00000000000002\\
22 &     -1.56248056455209\\
23 &      0.85478501960318\\
24 &     -0.90534916229569\\
25 &      0.45696099733973\\
36 &     -0.37224839333969\\
\end{tabular} \\
\hline\noalign{\smallskip}
\multicolumn{2}{c}{$|c(4)c(9)-c(36)|=4E-12$}\\
\multicolumn{2}{c}{$|c(21)|=2E-14$}\\
\end{tabular} &
\begin{tabular}[t]{lr}
\begin{tabular}[t]{d{1}d{12}}
\multicolumn{2}{l}{$R=6.046497437542$}\\
\multicolumn{1}{c}{$n$} &
\multicolumn{1}{c}{$a(n)$}            \\
\hline\noalign{\smallskip}
   0 &   0.000000000000   \\
   1 &   0.000000000000   \\
   2 &   1.000000000000   \\
   3 &   1.770795863371   \\
   4 &   0.000000000029   \\
   5 &   1.331470494003   \\
   6 &  -0.749395507674   \\
   7 &   0.074369313625   \\
   8 &   0.707106781209   \\
   9 &  -0.000000000016   \\
  10 &  -1.214451367832   \\
  11 &   0.620756243833   \\
  12 &   1.252141763204   \\
  13 &  -1.123686021586   \\
  14 &  -0.711286031630   \\
  15 &   0.964706032413   \\
  16 &  -0.000000000112   \\
  17 &  -0.000000000108   \\
  18 &  -0.128644484609   \\
  19 &  -0.499664871594   \\
  20 &   0.941491814843   \\
  21 &  -0.447011950327   \\
  22 &  -0.805501416928   \\
  23 &   1.236716356721   \\
  24 &  -0.529902643583   \\
  25 &   0.000000008756   \\
  36 &  -0.000000000168   \\
 288 &  -0.064322242377   \\
\end{tabular} \\
\hline\noalign{\smallskip}
\multicolumn{2}{c}{$|c(2\cdot3^{2})c(2\cdot4^{2})-c(2\cdot12^{2})|=8E-11$}\\
\multicolumn{2}{c}{$|c(17)|=1E-10$}\\
\end{tabular} 
\end{tabular}
\end{threeparttable}
\end{table}

\begin{table}
\caption{Supplemental table of Fourier coefficients for  
$\MAS \left(\Gamma_0(2),0,R\right)$ }
\label{tab:fourier_coeff_w_half_supp}
\begin{tabular}[t]{ll}
\begin{tabular}[t]{d{1}d{14}}
\multicolumn{1}{c}{$R$}
 &  8.92287648699174 \\
\multicolumn{1}{c}{$n$} 
&	\multicolumn{1}{c}{$A(n)$}\\
\hline\noalign{\smallskip}
2&   -0.70710678118654 \\
3&  1.10378899562734 \\
4&    0.49999999999993 \\
5&   0.90417459283958 \\
6&   -0.78049668380711 \\
7&    0.82934246755499 \\
8&  -0.35355339059330 \\
9&    0.21835014686776 \\
10&  -0.63934798597339 \\
\end{tabular}
&
\begin{tabular}[t]{d{1}d{14}}
\multicolumn{1}{c}{$R$}
 &  12.0929948750786  \\
\multicolumn{1}{c}{$n$} 
 &	\multicolumn{1}{c}{$A(n)$}\\
\hline\noalign{\smallskip}
2&    0.70710678118655 \\
3&  -0.70599475399569 \\
4&   0.49999999999999 \\
5&  -0.79974825694039 \\
6&   -0.49921367803249 \\
7&   -1.71337067862845 \\
8&    0.35355339059328 \\
9&   -0.50157140733054 \\
10&   -0.56550741572467 \\
\end{tabular}
\\
\noalign{\smallskip}
\hline\noalign{\smallskip}
\noalign{\smallskip}
\begin{tabular}[t]{d{1}d{15}r}
\multicolumn{1}{c}{$R$}
 &  13.77975135189073 \\
& \multicolumn{1}{c}{(even wrt $ z \mapsto - \frac{1}{2z} $)}\\
\multicolumn{1}{c}{$n$} 
&	\multicolumn{1}{c}{$A(n)$}\\
\hline\noalign{\smallskip}
2&    2.96351804031448 \\
3&   0.24689977245401 \\
4&    3.59139177031902 \\
5&  0.73706038534834 \\
6&    0.73169192981688 \\
7&   -0.26142007576500 \\
8&   2.60064131148226 \\
9&   -0.93904050235826 \\
10&    2.18429174878080 \\
\end{tabular}
&
\begin{tabular}[t]{d{1}d{15}}
\multicolumn{1}{c}{$R$} 
&  13.77975135189073 \\
& \multicolumn{1}{c}{(odd wrt $ z \mapsto - \frac{1}{2z} $)}\\
\multicolumn{1}{c}{$n$}  
&	\multicolumn{1}{c}{$A(n)$}\\
\hline\noalign{\smallskip}
2&   0.13509091556820 \\
3&    0.24689977245398 \\
4&  -0.79070303958101 \\
5&   0.73706038534830 \\
6&   0.03335391631439 \\
7&  -0.26142007576538 \\
8&  -1.36013067551284 \\
9&  -0.93904050236089 \\
10&   0.09957016228575 \\
\end{tabular} \\
\end{tabular}
\end{table}

\begin{table}
\centering
\caption{Comparison of Fourier coefficients $A(n)$ computed directly for  
$\MAS \left(\Gamma_0(N),0,R\right)$ 
 $(N=1,2)$, vs $\hat{A}(n)$ computed on 
$\MAS \left(\Gamma_0(4),\frac12,\frac12 R\right)$ and using (\ref{eq:shimura_corr_coeff}).}
\label{tab:fourier_coeff_comparison}
\begin{tabular}[t]{d{1}d{15}d{15}d{17}}
\multicolumn{1}{c}{$n$}            &
\multicolumn{1}{c}{$\hat{A}(n)$}   &
\multicolumn{1}{c}{$A(n)$}           &
\multicolumn{1}{c}{$|A(n)-\hat{A}(n)|$} \\
\hline\noalign{\smallskip}
\multicolumn{4}{c}{
$f \in \MAS (\Gamma_0(2),0,8.922876486992)$ ($t=1$)} \\
\noalign{\smallskip}
2  & -0.70710678118665 & -0.707106781186 & 0.6E-12 \\
3  &  1.10378899562739 & 1.103788995627 & 0.4E-12 \\
5  &  0.90417459283969 & 0.904174592840 & 0.3E-12 \\
\hline\noalign{\smallskip}
\multicolumn{4}{c}{
$f  \subseteq \MAS (\Gamma_0(1),0,13.779751351891)$ ($t=1$)}\\
\noalign{\smallskip}   
2 & 1.54930447794126 & 1.549304477941 & 0.7E-12 \\
3 & 0.24689977245398 & 0.246899772454 & 0.3E-12 \\
4 & 1.40034436536892 & 1.400344365369 & 0.2E-12 \\
5 & 0.73706038534387 & 0.737060385348 & 0.4E-11 \\
6 & 0.38252292109716 & 0.382522923066 & 0.2E-08 \\
\hline\noalign{\smallskip}
\multicolumn{4}{c}{
$f  \in \MAS (\Gamma_0(1),0,13.779751351891)$ ($t=2$)}\\
\noalign{\smallskip}
3 & 0.24689977245437 & 0.246899772454 & 0.8E-13 \\
5 & 0.73706038535004 & 0.737060385348 & 0.2E-11 \\
7 &-0.26142007624377 &-0.261420075765 & 0.5E-09 \\
9 &-0.93904050238904 &-0.939040502362 & 0.3E-10 \\
\hline\noalign{\smallskip}
\multicolumn{4}{c}{
$f \in \MAS (\Gamma_0(2),0,12.092994875079)$ ($t=2$)}\\
\noalign{\smallskip}
3 & -0.70599475379863  & -0.705994753996 & 0.2E-09 \\
5 & -0.79974825694696  & -0.799748256940 & 0.7E-11 \\
7 & -1.71337067860377  & -1.713370678628 & 0.2E-10 \\
9 & -0.50157140750090  & -0.501571407330 & 0.2E-09 \\
\noalign{\smallskip}
\multicolumn{4}{c}{In segment 3, the calculation is based on (\ref{eq:shimura_corr_coeff})
and the second portion of Table \ref{tab:N4R6889}. }
\end{tabular}
\end{table}

\begin{table}
\centering
\caption{ Comparison of Fourier coefficients for weights $k=9.044605824E-08$ and $k=0$
  near an ``avoided crossing''.}
\label{tab:cuspf_and_eisen_comp}
\begin{tabular}[h]{l}
\multicolumn{1}{l}{Corresponds to the cusp form} \\
\begin{tabular}[h]{lrrl}
	\begin{tabular}[t]{c}
	\multicolumn{1}{c}{$k$} \\
	\multicolumn{1}{c}{$R$} \\
		$c(2)$ \\
		$c(3)$ \\
		$c(4)$ \\
		$c(5)$ \\
		$c(6)$ \\
	\end{tabular}
&
	\begin{tabular}[t]{d{14}r}
           9.0446058240E-08 \\
	  13.77975135189074 \\
   	  1.54930480559976 \\
	  0.24689988546553 \\
	  1.40034433555250 \\
	  0.73706067260516 \\
	  0.38252272069428 \\
       \end{tabular}
&
	\begin{tabular}[t]{d{14}r}
		0 \\
		13.77975135189074 \\
		1.54930447794069 \\
	   	0.24689977245411 \\
	  	1.40034436536841 \\
		0.73706038534787 \\
		0.38252292306557 \\
	\end{tabular}
&
	\begin{tabular}[t]{d{7}}
\multicolumn{1}{c}{Difference}\\
		\\
	0.3E-06\\
	0.1E-06\\
	0.1E-06\\
	0.2E-06\\
	0.2E-06\\
	\end{tabular}
\end{tabular}
\\
\noalign{\smallskip}
\multicolumn{1}{l}{Corresponds to the Eisenstein series} \\
\begin{tabular}[h]{lrrl}
	\begin{tabular}[t]{c}
		$k$\\
		$R$\\
		$c(2)$ \\
		$c(3)$ \\
		$c(4)$ \\
		$c(5)$ \\
		$c(6)$ \\
	\end{tabular}
&
	\begin{tabular}[t]{d{14}}
   9.0446058240E-08 \\
  13.77975135138225\\
  -2.06525760334129 \\
  -1.72891679536648 \\
   2.97153986917404 \\
  -2.02747287754385 \\
   3.41008729668221 \\

 	\end{tabular}
&
	\begin{tabular}[t]{d{14}}
		0  \\
	13.77975135138225  \\
	-1.98398933080188 \\
	-1.68449330640991  \\
	2.93621366473571\\
	-1.96531634618530  \\
	 3.34201674772446 \\
	\end{tabular}
&
	\begin{tabular}[t]{d{7}}
\multicolumn{1}{c}{Difference}\\
		\\
	0.8E-01	\\
	0.4E-01	\\
	0.4E-01	\\
	0.6E-01	\\
	0.7E-01	\\
	\end{tabular}
\end{tabular}
\\
\noalign{\smallskip}
\multicolumn{1}{c}{The coefficients for the Eisenstein series at $k=0$
   were computed using (\ref{eq:eisenstein_coeff}), i.e.:} \\
\begin{tabular}[t]{l}
	$c(2)=2\cos (R\ln2)$ \\   
	$c(3)=2\cos (R\ln3)$ \\   
	$c(4)=1+2\cos (R\ln4)$ \\   
	$c(5)=2\cos (R\ln5)$ \\   
	$c(6)=2\cos (R\ln6)+2\cos (R(\ln 3 - \ln 2))$ \\   
\end{tabular}
\end{tabular}
\end{table}

\begin{table}
\centering
\caption{Comparison of Fourier coefficients for weights $k=9.044605824E-08$ and $k=0$
  ``far'' from an ''avoided crossing''.
The weight $0$ coefficients were computed using the formulas in Table \ref{tab:cuspf_and_eisen_comp}.}

\label{tab:cuspf_and_eisen_comp2}
\begin{tabular}[h]{lrrl}
  \begin{tabular}[t]{c}
	$k$\\
	$R$\\
	$c(2)$ \\
	$c(3)$ \\
	$c(4)$ \\
	$c(5)$ \\
	$c(6)$ \\
\end{tabular}
&
  \begin{tabular}[t]{d{14}}
9.0446058240E-08 \\
13.62696884857618   \\         
-1.99957085683552 \\
-1.48069687587703 \\
2.99828354611637\\
-1.99647405201235\\
2.96075820067617\\
\end{tabular}
&
  \begin{tabular}[t]{d{14}}
0 \\
13.62696884857618\\
-1.99957081810438 \\
 -1.48069680342062\\
 2.99828345661464\\
 -1.99647406885962\\
 2.96075811858031\\
\end{tabular}
&
  \begin{tabular}[t]{d{7}}
\multicolumn{1}{c}{Difference}\\
\\
0.4E-07\\ 
0.7E-07\\
0.7E-07\\
0.2E-07\\
0.8E-07\\
\end{tabular}
\\
\end{tabular}
\end{table}


%

\end{document}